%% file: diffeo.tex
\input dpmac

\ifx\documentclass\undefined
\catcode`^=7
\mathcode`@="8000 \catcode`@=13 \def@{\grabalpha{cat}\sf\empty} \catcode`@=12
\mathcode`"="8000 \catcode`"=13 \def"{\mathop\grabalpha{fun}\eurm\nolimits} \catcode`"=12
\mathcode`\`="8000 \catcode`\`=13 \def`{\grabalpha{simp}\bf\empty} \catcode"60=12 
\def\mathspecials{\catcode`@=13 \catcode`"=13 \catcode`\`=13\relax}
\catcode`^=\active

\font\sevenss=cmss8 at 7pt
\font\fivess=cmss8 at 5pt
\scriptfont\ssfam\sevenss
\scriptscriptfont\ssfam\fivess

\mathchardef\togets="181D

\expandafter\def\csname\string∞\endcsname{\ifmmode\infty\else{\tensy1}\fi}
\expandafter\def\csname\string⇄\endcsname{\togets}

\expandafter\def\csname\string⊠\endcsname{\boxtimes}

\font\msamdisplay=msam10 at 8pt
\font\msamtext=msam10 at 8pt
\font\msamscript=msam10 at 5pt
\font\msamscriptscript=msam10 at 3pt
\def\dochar#1#2{\vcenter{\kern.2ex\hbox{#2#1}\kern.2ex}}
\def\dodochar#1{\mathbin{\mathchoice{\dochar#1\msamdisplay}{\dochar#1\msamtext}{\dochar#1\msamscript}{\dochar#1\msamscriptscript}}}

\chardef\boxtimeschar"02
\def\boxtimes{\dodochar\boxtimeschar}

\chardef\boxemptychar"03
\def\square{\dodochar\boxemptychar}

{\catcode`\'=\active \gdef'{^\bgroup\primes}}
\def\primes{\prime\futurelet\next\primesbis}
\def\primesbis{\ifx'\next\let\nxt\primesbisbis \else\ifx^\next\let\nxt\primet \else\ifx\specialhat\next\let\nxt\primet
  \else\let\nxt\egroup\fi\fi\fi \nxt}
\def\primesbisbis#1{\primes} \def\primet#1#2{#2\egroup}

\begingroup\setbox0=\hbox{$"==colim "==lim "==hom @==op @==Top @==Set "==Fun @==sPSh$}\endgroup

\yearkeytrue
\suppressunusedbibfalse

\inewtoks\title
\inewtoks\author
\inewtoks\address
\inewtoks\email
\fi

\def\<{$$}
\def\>{$$} 
\let\(\relax
\let\)\relax
\def\typesetli{}
\def\typesetpar{\ppar}

\def\mpsfilename{diffeo.mps}
\setetok\preamble{%
outputtemplate := "\mpsfilename";
}

\def\typesettitle{%
\metadata{\the\title}{\the\author}%
{\font\articletitle=cmss17 \articletitle \leftskip0pt plus 20em \rightskip0pt plus 20em \parfillskip0pt \parindent0pt \baselineskip20pt \the\title \bigskip}%
{\font\authortitle=cmr14 \tabskip0pt plus 1fil \halign to\hsize{\hfil##\hfil\cr
\authortitle \the\author\vadjust{\medskip}\cr
Department of Mathematics and Statistics, Texas Tech University\cr
\https://dmitripavlov.org/\vadjust{\medskip}\cr
}}}

\inewif\ifpdf \pdffalse \ifx\pdfoutput\undefined\else\ifx\pdfoutput\relax\else\ifnum\pdfoutput<1 \else\pdftrue\fi\fi\fi

\ifx\documentclass\undefined\else 
\documentclass[nothms]{hha}

\usepackage[fresh-utf8]{inputenc}

\usepackage{newunicodechar}
\usepackage{silence}
\newunicodechar{⇄}{\rightleftarrows}
\newunicodechar{→}{\to}
\newunicodechar{←}{\gets}
\newunicodechar{⊠}{\boxtimes}
\newunicodechar{‖}{\|}
\newunicodechar{∞}{{\ifmmode\infty\else$\infty$\fi}}
\newunicodechar{…}{{\ifmmode\ldots\else\dots\fi}}
\newunicodechar{±}{\pm}

\let\frak\mathfrak

\def\<{$\def\qquad{\ }}
\def\>{$}
\def\(#1\){}

\DeclareMathAlphabet{\matheurm}{U}{eur}{m}{n}
\ifpdf
\usepackage{hyperref} 
\else
\usepackage[dvipdfmx]{hyperref} 
\fi
\usepackage{xcolor}
\definecolor{darkgreen}{rgb}{0,0.45,0}
\hypersetup{colorlinks,urlcolor=blue,citecolor=darkgreen,linkcolor=,linktoc=all,bookmarksnumbered}
\usepackage[capitalize]{cleveref}

\crefformat{equation}{(#2#1#3)}
\crefname{figure}{Figure}{Figures}
\crefformat{section}{\S#2#1#3}
\crefmultiformat{section}{\S\S#2#1#3}{ and~#2#1#3}{, #2#1#3}{, and~#2#1#3}
\usepackage{amsthm}
\counterwithin{equation}{section}
\theoremstyle{definition}

\numberwithin{equation}{section}

\def\catname#1{{\hyperref[cat.#1]{\mathsf{#1}}}}

\def\funname#1{{\hyperref[fun.#1]{\matheurm{#1}}}}

\def\id{{\rm id}}

\def\typesetli{\def\li{}}
\def\typesetpar{}

\def\li{\plainitem{$\bullet$}}
\def\plainitem{\endgraf\plainhang\textindent}
\def\plainhang{\hangindent\parindent}
\def\textindent#1{\indent\llap{#1\enspace}\ignorespaces}

\def\ltoarr#1{\mathop{\count0=#1 \loop\ifnum\count0>0 \smash-\mkern-7mu \advance\count0 -1 \repeat \mathord\rightarrow}\limits} 
\def\lto#1#2{\mathrel{\ltoarr{#1}^{#2}}} 
\def\longto#1^#2_#3{\mathrel{\ltoarr{#1}^{#2}_{#3}}} 
\def\lgetsarr#1{\mathop{\mathord\leftarrow \count0=#1 \loop\ifnum\count0>0 \mkern-7mu\smash-\advance\count0 -1 \repeat}\limits} 
\def\longgets#1^#2_#3{\mathrel{\lgetsarr{#1}\limits^{#2}_{#3}}} 

\def\gmatrix#1#2{\null\,\vcenter{\normalbaselines
        \ialign{#1\crcr
                \mathstrut\crcr\noalign{\kern-\baselineskip}
                #2\crcr\mathstrut\crcr\noalign{\kern-\baselineskip}}}\,}
\def\cdmatrix{\gmatrix{\hfil$##$\hfil&&\enspace\hfil$##$\hfil\enspace&\hfil$##$\hfil}}
\def\sqmatrix{\gmatrix{\hfil$##$&\enspace\hfil$##$\hfil\enspace&$##$\hfil}}
\def\cdbl{\def\normalbaselines{\baselineskip12pt \lineskip1pt \lineskiplimit1pt }}
\def\cd{\cdbl\cdmatrix}
\def\sqcd{\cdbl\let\vagap\;\sqmatrix}

\inewcount\arrowsize \arrowsize3
\def\mapright#1{\smash{\lto\arrowsize{#1}}}
\def\rvagap{\vagap} \def\lvagap{\vagap} \def\rvaskip{\vaskip} \def\lvaskip{\vaskip} \def\vaskip{} \def\vagap{}
\def\mapdown#1{\rvagap\big\downarrow\rlap{$\vcenter{\hbox{$\scriptstyle#1$}}$}\rvaskip}
\def\lmapdown#1{\lvaskip\llap{$\vcenter{\hbox{$\scriptstyle#1$}}$}\big\downarrow\lvagap}

\usepackage{graphicx}
\def\inlinemp#1{\includegraphics[scale=0.3]{diffeo.mps}}
\def\repo#1#2#3{\hskip0pt plus 5in \penalty200 \hskip0pt plus -5in \href{#1#3}{#2#3}}
\def\numdam{\repo{http://www.numdam.org/item/?id=}{numdam:}}
\def\eudml{\repo{https://eudml.org/doc/}{eudml:}}
\def\arXiv{\repo{https://arxiv.org/abs/}{arXiv:}}
\def\doi{\repo{https://doi.org/}{doi:}}
\fi 

\def\lex{\bar\times}
\def\E{{\sf E}}
\def\Ei{\E_{\sf\infty}}
\def\rdf{{\bf R}}
\def\B{{\sf B}} 
\def\tB{{\rm B}} 
\def\csp{{\sf B}_\smallint}
\def\cset{{\frak C}} 

\def\concr{\Upsilon} 
\def\ash{{\bf a}} 

\def\sicat{{\sf Δ}}
\def\sicatinj{\sicat_{\sf inj}}
\def\sSet{{\sf sSet}}
\def\sSetinj{\sSet_{\sf inj}}
\def\Fun{{\sf Fun}}
\def\op{{\sf op}}
\def\Set{{\sf Set}}

\def\inj{{\sf inj}}
\def\closed{{\sf closed}}
\def\EM{{\sf K}}
\def\HH{{\sf H}}
\def\Hom{\mathop{\sf Hom}}

\def\id{{\rm id}}
\def\diag{\mathop{\rm diag}\nolimits}

\mathcode`\:="603A 
\mathchardef\colon="303A 

\def\fbreak{\penalty-200 }
\def\gbreak{\hfil\penalty0 \hfilneg}

\title{Projective model structures on diffeological spaces and smooth sets and the smooth Oka principle}
\author{Dmitri Pavlov}
\ifx\documentclass\undefined\else 
\copyrightyear{2024}
\email{https://dmitripavlov.org/}
\submitted{}
\received{}
\published{}
\shorttitle{Model structures on diffeological spaces}
\address{Department of Mathematics and Statistics\\Texas Tech University} 
\EditFlow{221026-Pavlov}
\received{}
\revised{}
\accepted{}
\published{}
\submitted{J.~Daniel Christensen}
\volumeyear{}
\volumenumber{}
\issuenumber{}

\begin
{document}

{\setbox0\hbox{$\mathfrak{} \matheurm{}$}}

\makeatletter
\def\proof{\par\pushQED{\qed}\normalfont\trivlist\item[\hskip\labelsep\itshape Proof.]\ignorespaces}
\def\HHA@email#1{\href{#1}{\texttt{#1}}}
\makeatother
\fi 

\ifx\documentclass\undefined\typesettitle\fi 



\abstract Abstract.
In the first part of the paper,
we prove that the category of diffeological spaces does not admit a model structure transferred via the smooth singular complex functor from simplicial sets,
resolving in the negative a conjecture of Christensen and Wu, in contrast to Kihara's model structure on diffeological spaces constructed using a different singular complex functor.
Next,
motivated by applications in quantum field theory and topology,
we embed diffeological spaces into sheaves of sets (not necessarily concrete) on the site of smooth manifolds
and study the proper combinatorial model structure on such sheaves transferred via the smooth singular complex functor from simplicial sets.
We show the resulting model category to be Quillen equivalent to the model category of simplicial sets.
We then show that this model structure is cartesian,
all smooth manifolds are cofibrant,
and establish the existence of model structures on categories of algebras over operads.
Finally, we use these results to establish analogous properties for model structures on simplicial presheaves on smooth manifolds,
as well as presheaves valued in left proper combinatorial model categories,
and prove a generalization of the smooth Oka principle established in arXiv:1912.10544.
We apply these results to establish classification theorems for differential-geometric objects
like closed differential forms, principal bundles with connection, and higher bundle gerbes with connection
on arbitrary cofibrant diffeological spaces.

\ifx\documentclass\undefined\else 
\maketitle
\fi 

\tsection Contents

\contents

\def\@{}

\section Introduction

Diffeological spaces were introduced by Souriau\@ [\Souriau], with some closely related preceding work by Chen\@ [\Chen].
Stacey [\Stacey] and Baez–Hoffnung [\BH] give a review and comparison of these and other approaches to categories of smooth spaces.
Diffeological spaces contain many other categories of infinite-dimensional manifolds as full subcategories,
e.g., Fréchet manifolds by a result of Losik [\Losik].

By ^!{Grothendieck quasitopos}, the category $@Diffeo$ of diffeological spaces is a ^{Grothendieck quasitopos},
which is a particularly nice type of a category: it is complete and cocomplete, cartesian closed, and locally cartesian closed category.
Furthermore, $@Diffeo$ contains the category of smooth manifolds as a full subcategory.
This makes $@Diffeo$ a convenient category to work with infinite-dimensional mapping spaces of manifolds
and other smooth spaces more general than manifolds.
A book-length treatment by Iglesias-Zemmour\@ [\IZ] contains many examples illustrating the power of this formalism.

A closely related notion is that of ^{smooth sets}.
A ^{smooth set} (^!{smooth set}) is a sheaf of sets on the site of smooth manifolds and open covers.
A ^{diffeological space} (^!{diffeological space}) is a ^{smooth set}~$F$ that is a ^{concrete sheaf} (^!{concrete sheaf}):
if two sections $s,t∈F(M)$ ($M∈@Man$) coincide on every point $p:`R^0→M$ (meaning $p^*s=p^*t$, where $p^*:F(M)→F(`R^0)$),
then $s=t$.
Morphisms of ^{smooth sets} and ^{diffeological spaces} are simply morphisms of sheaves.
Thus, ^{smooth sets} contain ^{diffeological spaces} as a full subcategory.
The category of ^{smooth sets} is a Grothendieck topos,
so it inherits all the nice properties of ^{diffeological spaces}
and, in addition, it is a balanced category: if a morphism is a monomorphism and epimorphism,
then it is an isomorphism.
This last property is essential for showing that the category of abelian group objects in ^{smooth sets}
is a Grothendieck abelian category, which immediately allows for a development of homological algebra in this setting (to appear in a forthcoming paper).
In contrast, the category of abelian group objects in ^{diffeological spaces} is not an abelian category.

In complete analogy to topological spaces, one can define a (smooth) singular complex functor $"SmSing$ (^!{smooth singular complex}),
which endows the categories of ^{smooth sets} and ^{diffeological spaces} with a relative category structure:
a morphism~$f$ of ^{smooth sets} is a weak equivalence if $"SmSing f$ is a weak equivalence of simplicial sets.
Continuing the analogy to topological spaces, one can then inquire whether the resulting relative categories of ^{smooth sets} and ^{diffeological spaces} can be promoted
to model categories, by creating the class of fibrations using the functor $"SmSing$,
and whether this turns $"SmSing$ into a right Quillen equivalence of model categories.

Model structures in which the classes of weak equivalences and fibrations
are created by a right adjoint functor given by some sort of evaluation procedure
are commonly known as {\it projective model structures}.
For example, the projective model structure on simplicial presheaves is induced by the right adjoint functor
that evaluates a presheaf on all objects of the site and discards the data associated to morphisms of the site.
In our case, the ^{smooth singular complex} functor evaluates on all extended simplices and simplicial maps,
discarding the data associated to the other smooth maps.

The main result of the first part of this paper (^!{start of part 1}–^!{end of part 1})
is that the answer in the case of diffeological spaces is negative (but see ^!{surrogates} as well as Kihara\@ [\Ki],
who constructs a model structure on diffeological spaces using a different variant of the smooth singular complex functor).

\proclaim Theorem.
(See ^!{diffeological model structure} and the second part of ^!{Quillen equivalence for presheaves}.)
The category $@Diffeo$ of ^{diffeological spaces} (^!{diffeological spaces})
does not admit a model structure that is ^{transferred} (^!{transferred})
along the right adjoint functor
\<"SmSing:\fbreak @Diffeo→@sSet\>
(^!{smooth singular complex}), meaning its weak equivalences and fibrations are created by the functor $"SmSing$.
However, the functor $"SmSing$ is a Dwyer–Kan equivalence of relative categories.

The main result of the second part of this paper (^!{start of part 2}–^!{end of part 2})
is that for ^{smooth sets} we do indeed get a Quillen equivalence of model categories, with rather nice properties of involved model structures.
(The mere existence of the projective model structure on ^{smooth sets} is a special case of the Smith recognition theorem combined with Cisinski's results,
and the main difficulty lies in establishing all the additional nice properties listed in the statement below.)

\proclaim Theorem.
^^={main theorem for presheaves}
(See ^!{existence of model structures on presheaves},
^!{model structures on presheaves are cartesian},
^!{manifolds are cofibrant},
^!{smooth Oka principle for presheaves},
^!{properties of model structures on presheaves},
^!{admissibility of operads in presheaves},
^!{quasicategorical algebras in presheaves},
^!{transport of algebras in presheaves}.)
The category $@SmSet$ of ^{smooth sets} (^!{smooth sets})
admits a model structure ^{transferred} (^!{transferred})
along the functor
\<"SmSing:@SmSet→@sSet\>
(^!{smooth singular complex}), meaning its weak equivalences and fibrations are created by the functor $"SmSing$.
Smooth boundary inclusions and smooth horn inclusions form a set of generating cofibrations respectively generating acyclic cofibrations.
This model structure is left and right proper, combinatorial,
^{cartesian} (^!{cartesian}),
h-monoidal, symmetric h-monoidal, and flat (^!{properties of model structures on presheaves}).
All smooth manifolds~$M$ are cofibrant in this model structure
and for every smooth manifold~$M$ the internal hom functor $"Hom(M,-)$ preserves weak equivalences.
The functor $"SmSing$ is a right Quillen equivalence.
Operads in these model categories and algebras over them enjoy a good set of properties, as described
in ^!{admissibility of operads in presheaves}, ^!{quasicategorical algebras in presheaves}, ^!{transport of algebras in presheaves}.
Analogous results hold for the category $@PreSmSet$ of presheaves of sets.

In 1999, Hovey [\Hovey, Problem~2] already inquired whether sheaves on a manifold admit a model structure,
and the model structure studied in this paper can be seen as one possible answer to this question:
for a fixed manifold~$M$ we can take the slice model category of smooth sets over~$M$.

In the third part of the paper (^!{start of part 3}–^!{end of part 3}),
we extend ^!{main theorem for presheaves} to the case of sheaves and presheaves valued in a left proper combinatorial model category~$@V$,
such as simplicial sets or chain complexes.
This is relevant for applications, since many differential-geometric structures of interest
such as the moduli stack of principal $G$-bundles with connection or the moduli stack of higher bundle gerbes with connection
are encoded by such presheaves.
Once again, the mere existence of the projective model structure is a special case of the Smith recognition theorem,
and the main difficulty again lies in procuring the listed properties.

\proclaim Theorem.
^^={main theorem for enriched presheaves}
(See
^!{existence of model structures on enriched presheaves},
^!{operads in enriched presheaves},
^!{Quillen equivalence for enriched presheaves},
^!{smooth Oka principle for varieties}.)
Suppose $@V$ is a left proper combinatorial model category.
The category $@Sm_@V$ of $@V$-valued sheaves
and the category $@PreSm_@V$ of $@V$-valued presheaves
admit a model structure
with weak equivalences created by the ^{shape} functor (^!{shape})
and generating cofibrations analogous to those of ^!{main theorem for presheaves}.
This model structure is left proper, combinatorial,
and inherits from~$@V$ properties like being
monoidal, h-monoidal, symmetric h-monoidal, and flat (^!{existence of model structures on enriched presheaves}).
All smooth manifolds~$M$ are cofibrant in this model structure
and the internal hom functor $"Hom(M,-)$ preserves weak equivalences.
This model structure is Quillen equivalent to~$@V$ via a zigzag of Quillen equivalences.
Operads in these model categories and algebras over them enjoy a good set of properties
analogous to those of ^!{main theorem for presheaves}.

The closest in spirit to our paper is the work of Christensen–Wu\@ [\CW],
who develop the homotopy theory of diffeological spaces using the functor $"SmSing$ (^!{smooth singular simplicial set}).
In particular, we settle several of the conjectures stated in their paper,
including the nonexistence of a transferred model structure on diffeological spaces (^!{diffeological model structure},
which complements the existing work of Kihara\@ [\Ki] that constructs a transferred model structure on diffeological spaces for a different singular functor),
cofibrancy of smooth manifolds (^!{manifolds are cofibrant}),
cartesianness of the model structure on smooth sets (^!{model structures on presheaves are cartesian}),
in addition to the conjecture on
the coincidence of smooth homotopy groups of diffeological spaces with the simplicial homotopy groups of their ^{smooth singular simplicial sets}
(^!{smooth homotopy groups}),
which was already resolved in
Berwick-Evans–Boavida–Pavlov\@ [\BBP, Proposition~2.18].

\ifx\documentclass\undefined 
The reader may find the following chart of logical dependencies between sections useful.

\bigskip\centerline{\vbox{\halign{
&\hfil#\hfil\cr
  &     &   &     & 6   & $→$ & 5   &     & 11                       \cr 
  &     &   &     & $↓$ &     & $↑$ &     & $↓$                      \cr 
2 & $←$ & 3 & $←$ & 4   & $←$ & 7   & $←$ & 8   & $←$ & 9  & $←$ & 10\cr
  &     &   &     &     &     &     &     & $↑$                      \cr
  &     &   &     &     &     &     &     & 12  & $←$ & 13 & $←$ & 14\cr
}}}
\fi 

\subsection Previous work

Kihara [\Ki] constructs a cosimplicial object in diffeological spaces
by introducing a nonstandard diffeology on (nonextended) smooth simplices
that turns smooth horn inclusions into deformation retracts,
proves that the category of diffeological spaces admits a model structure
transferred along the singular complex functor associated to this cosimplicial object,
and shows that the resulting Quillen adjunction between simplicial sets and diffeological spaces
is a Quillen equivalence.
In the resulting model structure all diffeological spaces are fibrant
and by ^!{smooth homotopy groups}
combined with Kihara [\Ki, Theorem~1.4]
its weak equivalences coincide with the weak equivalences of Christensen–Wu [\CW, Definition~4.8],
which we also use in this paper (^!{weak equivalences of smooth sets}).
In particular, Kihara's model structure is connected to the model structure
of ^!{main theorem for presheaves} by a chain of Quillen equivalences.
Clough [\CCGT] continues this line of work,
exploring various Kihara-type model structures on smooth sets and simplicial smooth sets.

Kihara [\KBig, Theorem 1.11] proves that the class
of diffeological spaces that are smoothly homotopy equivalent to cofibrant
diffeological spaces is closed under gluing of D-numerable covers.
In particular, this class contains a large class of infinite-dimensional manifolds
(Kihara [\KBig, Theorem 11.1]).
It would be interesting to see whether cofibrancy in the Kihara model structure
could be established for differential-geometric objects like smooth manifolds.
\typesetpar
We also point out the ongoing work of Haraguchi–Shimakawa [\TS]
on a different model structure on diffeological spaces, which is not cofibrantly generated.

Cisinski [\THT, Théorème~3.9], [\FLA, Exemple~6.1.2, Théorème~6.1.8]
proves a general result that constructs a model structure on ^{smooth sets}
with monomorphisms as cofibrations.
The weak equivalences in Cisinski's model structure are the {\it shape equivalences\/}
(alias {\it Artin–Mazur equivalences}),
which coincide with the class of weak equivalences of ^!{weak equivalences of smooth sets}
by Berwick-Evans–Boavida–Pavlov\@ [\BBP, Proposition~1.3],
which shows that shape equivalences are created by the smooth singular complex functor.
Clough [\CCGT] uses Cisinski's methods to study various model structures on smooth sets and simplicial smooth sets.

Our article and Christensen–Wu [\CW] both use {\it extended\/} smooth simplices (^!{extended smooth simplex}), which are objects of $@Cart$.
This makes it particularly easy to establish the properties of the resulting model structure.
Other definitions of smooth simplices
and the corresponding model structures are explored in the work of Kihara [\Ki], Haraguchi–Shimakawa [\TS], Clough [\CCGT].

In the closely related subject of {\it simplicial smooth sets\/} (i.e., simplicial presheaves on the site of cartesian spaces or the site of smooth manifolds),
Morel–Voevodsky [\AHTS, Proposition~3.3.3] proved that the $`R$-local injective model structure on simplicial sheaves of sets on the site of sufficiently nice topological spaces
is Quillen equivalent to the Kan–Quillen model structure on simplicial sets.
Dugger [\UHT, Proposition~8.3] explicitly states the version for the case of the site of smooth manifolds.
Blander [\LPMS, Theorem~3.1] constructs $`R$-local projective model structures on simplicial presheaves and simplicial sheaves.
Schreiber [\DCCT, Definition~3.4.17] introduces the notion of an ∞-cohesive site
and proves [\DCCT, Proposition~4.4.6] a stronger result that cartesian spaces form an ∞-cohesive site.
Sati–Schreiber [\POC, §3.1.1] give a review of ∞-cohesive toposes.
Bunk [\Rlocal] also reviews and further develops the theory of $`R$-local localizations.
Amabel–Debray–Haine [\DC, §§4–5] develop a quasicategorical version of $`R$-local localizations for presheaves valued in presentable quasicategory.
Ayala–Francis–Rozenblyum [\AFR, §2, Lemma~2.3.16, Theorem~2.4.5]
contains related results that are proved in the more general context of stratified spaces,
although their results are restricted to isotopy sheaves of groupoids,
which excludes many simplicial presheaves, even set-valued ones.
We also point out the work of Sati–Schreiber–Stasheff [\Five, §3] and Fiorenza–Schreiber–Stasheff [\Cech, Appendix~A],
which contain early uses of simplicial presheaves on cartesian spaces in the context of quantum field theory,
as well as an early paper of Kock [\Conv, §5],
who already pointed out that the restriction functor from sheaves on manifolds to sheaves on cartesian spaces is an equivalence of categories.

The smooth Oka principle is due to Berwick-Evans–Boavida–Pavlov\@ [\BBP].
Additional applications of the smooth Oka principle can be found in Sati–Schreiber [\SS].
Another proof of a generalized form of the smooth Oka principle is in Clough [\HTDS, Theorem~B].

\subsection Acknowledgments

I thank
Urs Schreiber for a discussion that led to this paper and for pointing out the results of Cisinski [\THT],
Dan Christensen for pointing out ^!{nonconcrete realizations}, feedback on previous work, including the result of Kihara [\KBig, Theorem~11.1], and additional feedback on the paper.
Kiran Luecke for questions related to ^!{shapes of presheaves},
Adrian Clough for discussions concerning ^!{smooth horn cobase change},
and the anonymous referee of {\it Homology, Homotopy and Applications\/} for a careful reading of the manuscript and additional feedback that improved the paper.

\section Review of diffeological spaces and smooth sets
^^={start of part 1}

\proclaim Definition.
The small category $@=Cart$ of ^={cartesian space[s|]} is the full subcategory of the category of smooth manifolds and smooth maps
on objects~$X$ that are diffeomorphic to $`R^m$ for some $m≥0$ and, furthermore, the underlying set of~$X$ is a subset of $`R^n$ for some $n≥0$.
We turn $@Cart$ into a site by equipping it with the Grothendieck topology generated by the coverage of all open covers
whose finite intersections are empty or diffeomorphic to some $`R^m$ (hence, are objects in $@Cart$).

\proclaim Remark.
The site $@Cart$ (^!{cartesian spaces}) is a ^={concrete site[|s]} (Dubuc [\Dubuc, Definition~1.4])
meaning it has a terminal object $1=`R^0$ such that $"hom(1,-):@Cart→@Set$ is a faithful functor
and for any covering family $\{f_i: U_i→V\}_{i∈I}$ the induced map of sets \<∐_{i∈I}"hom(1,f_i):\fbreak ∐_{i∈I}"hom(1,U_i)→"hom(1,V)\>
is surjective.
On any concrete site one can define a concrete quasitopos (Dubuc [\Dubuc, Definition~1.3]) of ^={concrete shea[ves|f]} (Dubuc [\Dubuc, Definition~1.5]),
where a presheaf \<F:@Cart^@op→@Set\> is ^={concrete} if the canonical map \<F(X)→"hom("hom(1,X),F(1))\>
adjoint to the map \<F(X)⨯"hom(1,X)→F(1)\> induced by the structure maps of the presheaf~$F$
is an injection of sets.
^^={concrete preshea[f|ves]}

\proclaim Remark.
The category of ^{concrete sheaves} on any small ^{concrete site} (^!{concrete site})
is a ^={Grothendieck quasitopos[|es]} (Penon [\PenonI, \PenonII], Dubuc [\Dubuc, Theorem~1.7], Baez–Hoffnung [\BH, Theorem~52 (arXiv); 5.25 (journal)], Johnstone [\SoE, Theorem C2.2.13]).
Any Grothendieck quasitopos is a locally presentable category that is locally cartesian closed.

We now introduce the main categories of this paper.

\proclaim Definition.
The Grothendieck topos \<@=SmSet\>
of ^={smooth set[s|]} is the category of sheaves of sets on the site $@Cart$ (^!{cartesian spaces}).

\proclaim Definition.
The Grothendieck topos \<@=PreSmSet\>
of ^={presmooth set[s|]} is the category of presheaves of sets on the site $@Cart$.

\proclaim Remark.
The inclusion $@SmSet→@PreSmSet$ (^!{smooth sets}, ^!{presmooth sets}) exhibits $@SmSet$ as a reflective subcategory of $@PreSmSet$.
In particular, we have a left adjoint reflection functor $\ash:@PreSmSet→@SmSet$,
known as the ^={associated sheaf functor[|s]}.
^^={associated shea[f|ves]}

A precursor for the following definition can be found in Chen [\Chen],
the modern definition first appeared in Souriau [\Souriau],
and a book-length treatment is given by Iglesias-Zemmour [\IZ].

\proclaim Definition.
The Grothendieck quasitopos \<@=Diffeo\>
of ^={diffeological space[s|]} is the category of ^{concrete sheaves} of sets (^!{concrete sheaf}) on the site $@Cart$ (^!{cartesian spaces}).

\proclaim Definition.
The Grothendieck quasitopos \<@=PreDiffeo\>
of ^={prediffeological space[s|]} is the category of ^{concrete presheaves} of sets on the site $@Cart$.

\proclaim Remark.
The inclusion $@PreDiffeo→@PreSmSet$ (^!{prediffeological spaces}, ^!{presmooth sets}) is a reflective subcategory.
In particular, we have a left adjoint reflection functor,
known as the ^={concretization functor[|s]}
$\concr:@PreSmSet→@PreDiffeo$.
^^={concretization[|s]}
Concretely, the reflection map $F→G$ is the quotient map of presheaves
that identifies two sections $s,t∈F(U)$ ($U∈@Cart$) if for all $u:`Δ^0→U$ we have $su=tu$,
i.e., $s$ and $t$ induce the same maps on the underlying sets of points.
The inclusion $@Diffeo→@SmSet$ (^!{diffeological spaces}, ^!{smooth sets}) is also a reflective subcategory,
with $\ash\concr$ as the reflection functor (^!{concretization functor}, ^!{associated sheaf}).

\proclaim Remark.
^^={colimits of presheaves}
Limits in the categories $@PreSmSet$, $@SmSet$, $@PreDiffeo$, and $@Diffeo$ are computed objectwise,
since the sheaf property and ^{concrete presheaf} property are preserved under limits.
Colimits in these categories are computed as follows.
\typesetli
\li In $@PreSmSet$: objectwise.
\li In $@SmSet$: apply the ^{associated sheaf} functor~$\ash$ (^!{associated sheaf}) to the colimit in $@PreSmSet$.
\li In $@PreDiffeo$: apply the ^{concretization} functor~$\concr$ (^!{concretization}) to the colimit in $@PreSmSet$.
\li In $@Diffeo$: apply $\ash\concr$ to the colimit in $@PreDiffeo$.
Since the latter is always a separated presheaf, the ^{associated sheaf} can be computed using the plus construction.

\proclaim Definition.
The category $@=Man$ of ^={smooth manifold[s|]} has smooth manifolds as objects and smooth maps as morphisms.
To make $@Man$ a small category, we take the full subcategory on smooth manifolds whose underlying set is a subset of some $`R^n$ (ignoring its topology).
We turn $@Man$ into a small site by equipping it with the Grothendieck topology generated by the coverage of all open covers.

\proclaim Remark.
^^={sites for smooth sets}
The restriction functor along the inclusion of sites \<@Cart→@Man\>
(^!{cartesian spaces}, ^!{smooth manifolds})
induces equivalences of categories of sheaves of sets,
as well as ^{concrete sheaves} of sets (^!{concrete sheaf}).

\proclaim Remark.
^^={embedding of manifolds into presheaves}
The (restricted) Yoneda embedding construction induces fully faithful functors (generically denoted by~$`y$)
\<@Man→@Diffeo,\qquad @Man→@PreDiffeo,\qquad @Man→@SmSet, \qquad @Man→@PreSmSet\>
(^!{diffeological spaces}, ^!{prediffeological spaces}, ^!{smooth sets}, ^!{presmooth sets}).
We often omit these functors from our notation when it causes no ambiguity.

\section Smooth singular complex and realization

\proclaim Definition.
The category $@=Δ$ of simplices is the category of finite nonempty totally ordered sets and order-preserving maps.
To make $@Δ$ small, we restrict to the full subcategory of objects given by standard simplices $[m]=\{0<⋯<m\}$ for all $m≥0$.
The category $@=sSet$ of simplicial sets is defined as the category of presheaves of sets on~$@Δ$.

\proclaim Definition.
The functor \<`=Δ:@Δ→@Cart\> (^!{cartesian spaces}) sends a simplex $[m]$ to the ^={extended smooth simpl[ex|ices]} \<`Δ^m=\left\{x∈`R^{[m]}\Bigm| ∑_{i∈[m]}x_i=1\right\}\>
and a map of simplices $f:[m]→[n]$ to the smooth map \<`Δ^f:`Δ^m→`Δ^n,\qquad x↦\left(j↦∑_{i: f(i)=j}x_i\right).\>

\proclaim Definition.
The adjunction \<‖{-}‖:@sSet⇄@PreSmSet:"SmSing\>
is the nerve-\hskip0pt realization adjunction associated to the cosimplicial object \<`Δ:@Δ→@Cart→@PreSmSet\> (^!{extended smooth simplex}, ^!{embedding of manifolds into presheaves}).
The right adjoint is the ^={smooth singular simplicial set[|s]}
(alias ^={smooth singular complex[|es]}) functor \<"=SmSing:@PreSmSet→@sSet,\>
which sends some $F∈@PreSmSet$ to the simplicial set $[n]↦F(`Δ^n)$ and likewise for simplicial structure maps.
The left adjoint is the realization functor associated to~$`Δ$,
which sends a simplicial set~$X$ to \<"colim_{x∈@Δ/X} U(x),\>
where $@Δ/X$ is the category of simplices of~$X$ (objects are pairs $([m]∈@Δ,x∈X_m)$,
morphisms $([m],x)→([n],y)$ are maps of simplices $f:[m]→[n]$ such that $X_f(y)=x$)
and $U:@Δ/X→@Δ→@Cart→@PreSmSet$ denotes the forgetful functor $([m],x)↦[m]$ composed with the functor~$`Δ$ of ^!{extended smooth simplex}
and the Yoneda embedding of ^!{embedding of manifolds into presheaves}.
Analogous adjunctions with $@PreSmSet$ replaced by the categories $@SmSet$, $@PreDiffeo$, $@Diffeo$
have the corresponding restrictions of $"SmSing$ as the right adjoints
and the functors $|{-}|=\ash‖{-}‖$ (the ^={smooth realization[|s]} functor), $\concr‖{-}‖$, $\ash\concr‖{-}‖$
as the left adjoints, respectively.

The following example was inspired by a discussion with J.~Daniel Christensen.

\proclaim Remark.
^^={nonconcrete realizations}
The simplicial sets $A=Δ^2/(d_0 σ \sim d_1 σ)$, where $σ$ is the nondegenerate 2-simplex,
and $B=(Δ^2⊔Δ^2)/(ι_1 ∂Δ^2\sim ι_2 ∂Δ^2)$, where $ι_1$ and $ι_2$ are the embeddings of summands,
have nonconcrete smooth realizations (under $‖{-}‖$ or $|{-}|$), as witnessed by the following example of different sections that have the same underlying map of sets.
For $A$, one section~$s$ is an injective map that traverses the faces $d_0 σ$ and $d_1 σ$ smoothly, with vanishing derivatives at the midpoint,
and the other section~$t$ traverses $d_0 σ$ back and forth.
Once we identify $d_0 σ$ and $d_1 σ$, the two sections have the same underlying map of sets, but are not equal as sections
since their germs at the vertex~2 of~$σ$ are induced by different smooth sections of~$σ$.
For~$B$, take the same section~$s$ together with a section~$r$ of the second copy of~$|Δ^2|$
that traverses the faces $d_0 σ$ and $d_1 σ$ smoothly, with vanishing derivatives at the midpoint.

\proclaim Proposition.
^^={concrete realization[|s]}
Consider the full subcategory $\sSet'⊂\sSet$
comprising simplicial sets~$X$ such that
every nondegenerate simplex in~$X$ yields a monomorphism of simplicial sets $Δ^n→X$
and the intersection (pullback) of any two such simplices is either empty or is another nondegenerate simplex in~$X$.
The restriction of the ^{smooth realization} functor $‖{-}‖:\sSet→@PreSmSet$ (^!{smooth realization})
to the full subcategory $\sSet'⊂\sSet$ factors through $@PreDiffeo$, i.e., lands in ^{concrete presheaves}.

\proof Proof.
Suppose $s,t:U→‖X‖$ are two sections of $‖X‖$ over $U∈@Cart$.
Then we have $s=‖σ‖∘f$ and $t=‖τ‖∘g$, for some $σ:Δ^m→X$, $τ:Δ^n→X$ and $f:U→‖Δ^m‖$, $g:U→‖Δ^n‖$.
By the Eilenberg–Zilber lemma, we can assume $σ$ and $τ$ to be nondegenerate.
We can also assume that $f$ and $g$ do not factor through any proper faces of $‖Δ^m‖$ and $‖Δ^n‖$, respectively.
By assumption, the pullback $Δ^m⨯_X Δ^n$ is a nondegenerate simplex $ρ:Δ^k→X$, through which both $f$ and $g$ must factor.
The maps $f$ and $g$ do not factor through proper faces, so we get $ρ=σ=τ$.
By assumption, the map $σ:Δ^m→X$ is a monomorphism, hence its realization $‖σ‖:‖Δ^m‖→‖X‖$ is also a monomorphism.
The images of $f$ and $g$ under $‖σ‖$ have the same underlying maps of sets, therefore $f$ and $g$ have the same underlying maps of sets.
Since $‖Δ^m‖$ is a concrete sheaf, we obtain $f=g$.

\proclaim Corollary.
^^={examples of concrete realizations}
If $X$ is one of the simplicial sets $Δ^n$ ($n≥0$), $∂Δ^n$ ($n≥0$), or a horn $Λ^n_k$ ($n>0$, $0≤k≤n$),
then $‖X‖$ and $|X|$ are ^{concrete presheaves}.

\proclaim Definition.
^^={weak equivalence[|s] of smooth sets}
The category $@PreSmSet$ (^!{presmooth set}) is turned into a relative category
by postulating that its weak equivalences are precisely
those morphisms whose image under $"SmSing$ (^!{smooth singular simplicial set}) is a weak equivalence of simplicial sets.
The categories $@SmSet$ (^!{smooth set}), $@Diffeo$ (^!{diffeological space}), and $@PreDiffeo$ (^!{prediffeological space}) are turned into relative categories in the same way.

\proclaim Definition.
A ^={smooth homotopy} between morphisms $f,g:A→B$ in $@PreSmSet$
is a morphism of ^{presmooth sets}
\<h:`Δ^1⨯A→B, \qquad h∘ι_0=f, \qquad h∘ι_1=g,\>
where the corresponding inclusion is denoted by \<ι_k:A→\{k\}⨯A→`Δ^1⨯A.\>
A ^={smooth homotopy equivalence} is a map $f:A→B$ in $@PreSmSet$
such that there is a map $g:B→A$ with a smooth homotopy from $\id_A$ to $gf$ and a smooth homotopy from $fg$ to $\id_B$.
A ^={smooth deformation retraction} is a map $f:A→B$ in $@PreSmSet$
that can be made into a smooth homotopy equivalence in such a way that $\id_A=gf$.

\proclaim Proposition.
The functor $"SmSing$ sends
smoothly homotopic maps in the category $@PreSmSet$ to simplicially homotopic maps in $@sSet$,
smooth homotopy equivalences in $@PreSmSet$ to simplicial homotopy equivalences in $@sSet$,
and smooth deformation retractions in $@PreSmSet$ to simplicial deformation retractions in $@sSet$.

\proof Proof.
(See also Christensen–Wu [\CW, Lemma~4.10].)
This follows immediately from the fact that $"SmSing$ is a right adjoint,
in particular, it preserves small limits such as products used in the definition of a smooth homotopy.
The canonical map $Δ^1→"SmSing `Δ^1$ can be used to extract simplicial homotopies from $"SmSing$ evaluated on smooth homotopies.

\section The associated sheaf and concretization

The following result provides a powerful tool to work with colimits of ^{smooth sets},
by allowing us to replace them with colimits of ^{presmooth sets}, which are much easier to work with because colimits of presheaves are computed objectwise.

\proclaim Proposition.
^^={sheafification is a weak equivalence}
Suppose $F∈@PreSmSet$ (^!{presmooth sets}) and $s:F→G=\ash F$ is the canonical morphism from~$F$ to its ^{associated sheaf}~$G=\ash F$ (^!{associated sheaf}).
Then the map $F→G$ is a weak equivalence in the relative category $@PreSmSet$ (^!{weak equivalence of smooth sets}).
More generally, any local isomorphism of presheaves is a weak equivalence in $@PreSmSet$.

\proof Proof.
Consider the model category~$M$ of simplicial presheaves on the site $@Cart$ equipped with its injective model structure
left Bousfield localized at Čech nerves of good open covers.
Consider the functor~$L$ from~$M$ to simplicial sets
that sends a simplicial presheaf~$F$ to the diagonal of the bisimplicial set $n↦F(`Δ^n)$.
The functor~$L$ is a left adjoint functor that preserves monomorphisms and objectwise weak equivalences.
Furthermore, by Borsuk's nerve theorem (for example, combine Weil [\TDR, §5] and Eilenberg [\SHDM, Theorem~II]),
the functor~$L$ sends the Čech nerve of a good open cover to a weak equivalence of simplicial sets.
Thus, $L$ is a left Quillen functor that preserves weak equivalences.
In the model category~$M$, the map $F→G$ is a weak equivalence, hence so is $L(F)→L(G)$.
It remains to observe that for presheaves of sets we have $L="SmSing$, so $"SmSing F→"SmSing G$ is a weak equivalence of simplicial sets
and $F→G$ is a weak equivalence in $@PreSmSet$.

The following special case of ^!{sheafification is a weak equivalence} is important enough to be stated separately.

\proclaim Proposition.
^^={homotopy colimits of presheaves}
Suppose $D:I→@SmSet$ (^!{smooth sets}) is a diagram of ^{smooth sets},
$G$ its colimit,
and $F$ its colimit in the category $@PreSmSet$ (^!{presmooth sets}).
Then the canonical map $F→G$ is a weak equivalence (^!{weak equivalence of smooth sets}).

\proof Proof.
Combine ^!{colimits of presheaves} and ^!{sheafification is a weak equivalence}.

\proclaim Remark.
^^={concretization changes homotopy type}
The analogous result for the concretization functor~$\concr$ (^!{concretization}) is false.
Consider the sheaf~$F$ of closed differential $n$-forms, where $n>0$.
This sheaf is not concrete and its ^{concretization} is~$`Δ^0$ because $F(`Δ^0)$ is a single point.
However, the map $F→`Δ^0$ is not a weak equivalence
because $π_n("SmSing F)≅`R$, with the isomorphism given by integrating a closed differential $n$-form along $n$-dimensional singular simplices;
the Stokes formula then shows that homotopic pointed spheres map to the same real number.

\section Model categories

In this section, we recall some facts about model categories.

\proclaim Proposition.
(The ^={Kan–Quillen model structure on simplicial sets}.)
^^={Kan–Quillen model structure}
The category $@sSet$ admits a ^{cartesian}
combinatorial proper model structure
whose generating cofibrations are ^={boundary inclusion[s|]} \<δ_n:∂Δ^n→Δ^n\qquad (n≥0)\>
and generating acyclic cofibrations are ^={horn inclusion[s|]} \<λ_{n,k}:Λ^n_k→Δ^n\qquad (n>0,\ 0≤k≤n).\>
This model structure is proper, cartesian, its weak equivalences are closed under filtered colimits,
and all objects are cofibrant.

\proclaim Definition.
Suppose $C$ is a model category and $R:D→C$ is a right adjoint functor.
The ^={transferred model structure[|s]} on~$D$ (if it exists)
^^={transferred}
is the unique model structure whose weak equivalences and fibrations are created by the functor~$R$.

\proclaim Proposition.
^^={transfer theorem}
(Crans [\QCMSS, Theorem~3.3], Hirschhorn [\MCL, Theorem~11.3.2].)
Suppose $C$ and $D$ are locally presentable categories, $L⊣R:C⇄D$ is an adjunction,
and $C$ is equipped with a cofibrantly generated (hence combinatorial) model structure.
The ^{transferred model structure} (^!{transferred model structure}) on $D$ exists if and only if
the functor~$R$ sends transfinite compositions of cobase changes of elements of $L(J)$
to weak equivalences in~$C$,
where $J$ denotes any generating set of acyclic cofibrations in~$C$.
Given a set~$I$ of generating (acyclic) cofibrations of~$C$,
the set $L(I)$ is a set of generating (acyclic) cofibrations of~$D$.

\proclaim Proposition.
^^={recognition theorem}
(Barwick [\LR, Proposition~1.7 (arXiv); 2.2 (journal)], Beke [\SHMC, Theorem~1.7],
Lurie [\HTT, Proposition~A.2.6.15 (website); A.2.6.13 (printed)].)
Suppose $C$ is a locally presentable category
and $W$ is a class of morphisms in $C$ that is closed under the 2-out-of-3 property
and is given by the closure under filtered colimits of a set of objects in the category of morphisms and commutative squares in~$C$.
Suppose $I$ is a set of ^{h-cofibrations} (^!{h-cofibrations}) in the relative category $(C,W)$
such that morphisms with the right lifting property with respect to~$I$
necessarily belong to~$W$.
Then $C$ admits a left proper combinatorial model structure
whose class of weak equivalences is given by~$W$
and $I$ is its set of generating cofibrations.

\proclaim Corollary.
^^={h-recognition theorem}
Suppose $C$ is a left proper combinatorial model category and $I$ is a set of ^{h-cofibrations} (^!{h-cofibrations}) in the relative category $(C,W)$
such that $I$ contains some set of generating cofibrations for~$C$.
Then $C$ admits a left proper combinatorial model structure~$M$ whose class of weak equivalences is given by~$W$ and $I$ is its set of generating cofibrations.
The identity functor $C→M$ is a left Quillen equivalence.

\proclaim Definition.
A model category~$C$ is ^={cartesian}
if its underlying category is cartesian closed
(meaning for every $A∈C$ the functor $A⨯-:C→C$ has a right adjoint functor $"=Hom(A,-):C→C$),
the terminal object is cofibrant,
and the pushout product \<A⨯D⊔_{A⨯C}B⨯C→B⨯D\> of a cofibration $A→B$ and an (acyclic) cofibration $C→D$
is an (acyclic) cofibration.

We simplify the unit condition in the following definition since in our case all units are cofibrant.

\proclaim Definition.
A ^={weak monoidal Quillen adjunction[|s]}
^^={weak monoidal Quillen equivalence[|s]}
(Schwede–Shipley [\EMMC, Definition~3.6])
is a Quillen adjunction $L:@C⇄@D:R$ between monoidal model categories
such that the right adjoint functor~$R$ is a lax monoidal functor,
for any cofibrant objects $A,B∈@C$ the {\it comonoidal map\/}
\<L(A⊗B)→LA⊗LB\>
defined as the adjoint of the composition
\<A⊗B\lto7{η_A⊗η_B}RLA⊗RLB\lto7{}R(LA⊗LB)\>
is a weak equivalence,
and the map
\<L1_@C→1_@D\>
adjoint to the map $1_@C→R1_@D$ is a weak equivalence.

\section Projective model structure on diffeological spaces
^^={end of part 1}

In this section we prove that the Kan–Quillen model structure on $@sSet$
does not transfer along the right adjoint functor $"SmSing:@Diffeo→@sSet$ (^!{smooth singular simplicial set}).
This is caused by the pathological behavior of colimits in $@Diffeo$:
colimits in $@PreDiffeo$ are computed as the ^{concretizations} (^!{concretization}) of colimits in $@PreSmSet$
and colimits in $@Diffeo$ are computed as the ^{associated sheaves} of colimits in $@PreDiffeo$.
The concretization functor~$\concr$ (^!{concretization}) can change the homotopy type dramatically, as shown in ^!{concretization changes homotopy type}.
In particular, the concretization functor can interact in a wild way with cobase changes of smooth horn inclusions $|Λ^n_k|→|Δ^n|$,
and this section exploits this behavior to construct a cobase change of the smooth 3-horn that is not a weak equivalence,
which disproves the existence of a transferred model structure.

As shown in the next section, enlarging the category $@Diffeo$ to $@SmSet$ allows us to prove the existence of the transferred model structure.
The content of this section is not used anywhere else in the paper.
Its only purpose is to motivate the enlargement of the category of diffeological spaces to the category of smooth sets.

\proclaim Definition.
^^={bad section}
The injective smooth map \<S:`S^1→|Δ^3|\> is defined as follows.
We parametrize $|Δ^3|=\{(x,y,z)∈`R^3\}$, with the four faces of $|Δ^3|$ being $x+y+z=1$, $x=0$, $y=0$, $z=0$.
Denote by \<\cset=[0,1]∖⋃_{n≥0}⋃_a (z_a+3^{-n-1},z_a+2⋅3^{-n-1})\>
the Cantor set, where $a:\{0,…,n-1\}→\{0,2\}$
and $z_a=∑_{0≤k<n}a_k 3^{-k-1}$.
Denote by $b:`R→`R$ a smooth function that maps $(0,1)$ to itself and vanishes on the complement $(-∞,0]∪[1,∞)$.
Now identify $`S^1=[0,4]/(0∼4)$ and set
$$S(x)=d(x)+∑_{n=2k≥0}∑_a(0,b_{n,a}(x),0)+∑_{n=2k+1≥0}∑_a(0,0,b_{n,a}(x)),$$
where \<b_{n,a}(x)=\allowbreak 3^{-n^2-1} b((x-z_a-3^{-n-1})3^{n+1}),\>
$$\hskip0pt plus 2em d(x)=(c(x),0,0)+(0,c(x-1),0)+(-c(x-2),0,0)+(0,-c(x-3),0),$$
and $c:`R→`R$ is a smooth function such that $c(x)=0$ for all $x≤0$, $c(x)=1$ for all $x≥1$, and $c$ is strictly increasing on~$[0,1]$.
Thus, the image of the smooth map~$d$ looks like a square with vertices $(0,0)$, $(1,0)$, $(1,1)$, $(0,1)$.
The factor of $3^{-n^2-1}$ in $b_{n,a}$ guarantees that the resulting function~$S$ is smooth.

\proclaim Remark.
Taking $d$ and the summands with $n≤2$ in the formula for~$S$ yields a function that can be schematically depicted by the following graph,
where the horizontal axis is~$x$ (depicting only $x∈[0,1]$) and the vertical axis is the normal coordinate with respect to the line $(1,0,0)$;
the part above the horizontal line depicts the coordinate $y≥0$ in the plane $z=0$,
whereas the bottom part depicts the coordinate $z≥0$ in the plane $y=0$.
The remaining values of $x∈[1,4]$ close the loop by a unit square in the half-plane $z=0$, $y≥0$.
$$\inlinemp{
u := 5cm;
draw (0,0)..
(1/27u,0){1,0}..(3/54u,1/27u)..{1,0}(2/27u,0)..
(1/9u,0){1,0}..(3/18u,-1/9u)..{1,0}(2/9u,0)..
(7/27u,0){1,0}..(15/54u,1/27u)..{1,0}(8/27u,0)..
(1/3u,0){1,0}..(3/6u,1/3u)..{1,0}(2/3u,0)..
(19/27u,0){1,0}..(39/54u,1/27u)..{1,0}(20/27u,0)..
(7/9u,0){1,0}..(15/18u,-1/9u)..{1,0}(8/9u,0)..
(25/27u,0){1,0}..(51/54u,1/27u)..{1,0}(26/27u,0)..
(1u,0);
}
$$

The idea behind~$S$ is that it oscillates countably many times between the different faces of $|Δ^3|$,
while not factoring through $|Λ^3_0|→|Δ^3|$ because at points~$p$ in the Cantor set~$\cset$, the restriction of~$f$ to any neighborhood of~$p$
straddles both faces $y=0$ and $z=0$ of~$|Δ^3|$.

\proclaim Definition.
^^={bad cobase change}
Denote by $F$ the subobject of~$|Δ^3|$ given by the intersection of all subobjects that contain $|Λ^3_0|$
together with the plots $`R^n→|Δ^3|$ given by composing the injective smooth map $S:`S^1→|Δ^3|$ (^!{bad section}) with arbitrary smooth maps $`R^n→`S^1$.
We have canonical maps $|Λ^3_0|→F→|Λ^3_0|_{|Δ^3|}$,
whose underlying maps of sets are bijections,
where $|Λ^3_0|_{|Δ^3|}$ comprises all plots of $|Δ^3|$ whose underlying map of sets factors through the underlying set of $|Λ^3_0|$.
Although $|Λ^3_0|$ and $|Λ^3_0|_{|Δ^3|}$ are smoothly contractible, we will see that $F$ is not, thanks to the special properties of the section $S$ explored below.
The (unique) factorization of~$S$ through~$F$ is denoted by $s:`S^1→F$.

\proclaim Remark.
^^={properties of the bad section}
By construction, the maps $s:`S^1→F$, $|Λ^3_0|→F$, and $F→|Δ^3|$ are monomorphisms.
The restriction of $S:`S^1→|Δ^3|$ to $`S^1∖\cset$ factors through the inclusion $|Λ^3_0|→|Δ^3|$.
If $U⊂`S^1$ is an open subset such that $U∩\cset≠∅$, then the restriction of $S$ to~$U$ does not factor through the inclusion $|Λ^3_0|→|Δ^3|$.

\proclaim Proposition.
^^={smooth horn cobase change}
The cobase change of the smooth horn inclusion $|Λ^3_0|→|Δ^3|$ along the inclusion $|Λ^3_0|→F$ (^!{bad cobase change}) in the category $@Diffeo$ is not a weak equivalence.

\proof Proof.
The objects in the statement are ^{concrete presheaves} by ^!{examples of concrete realizations}.
The cobase change in the category $@Diffeo$ is the ^{associated sheaf} (^!{associated sheaf}) of the ^{concretization} (^!{concretization}) 
of the cobase change in the category $@PreSmSet$.
The ^{concretization functor}~$\concr$ identifies the section~$s:`S^1→F$ that was manually added to $|Λ^3_0|$ to form~$F$ with the section $S:`S^1→|Δ^3|$.
Therefore, the pushout in the category $@PreDiffeo$ is the inclusion $F→|Δ^3|$.
Since $"SmSing |Δ^3|$ is contractible, we have to show that $"SmSing F$ is not contractible.
By Berwick-Evans–Boavida–Pavlov\@ [\BBP, Proposition~2.18], it suffices to show that the morphism $s:`S^1→F$ is not smoothly homotopic to a constant map $t:`S^1→F$
via a smooth homotopy $H:`R⨯`S^1→F$ such that $H|_{0⨯`S^1}=s$, $H|_{1⨯`S^1}=t$.
We continue to use the conventions of ^!{bad section},
identifying $`S^1=[0,4]/(0∼4)$ and equipping it with the positive orientation induced from~$[0,4]$.
\ppar
Fix some $r∈`R$.
Denote by $c_r:`S^1≅\{r\}⨯`S^1→`R⨯`S^1$ the canonical inclusion.
Every section of~$F$ factors locally through the map $|Λ^3_0|→F$ or the map $s:`S^1→F$.
The image of~$c_r$ in $`R⨯`S^1$ is compact and therefore can be covered by open subsets of $`R⨯`S^1$
such that the restriction of $H:`R⨯`S^1→F$ to every subset factors through $|Λ^3_0|→F$ or $s:`S^1→F$.
By definition of the product topology on $`R⨯`S^1$,
we may assume these open sets to be products of an open interval~$R$ in~$`R$ and an open interval~$W$ in~$`S^1$.
Since $`S^1$ is compact, we can assume there are only finitely many such sets $R_i⨯W_i$.
By replacing every $R_i$ with the intersection $R=⋂_i R_i$ (which contains $r∈`R$), we may further assume that the interval~$R$ is the same for all subsets.
By shrinking and refining the intervals~$W_i$ in~$`S^1$ as necessary, we get a cyclically ordered set of open intervals~$W_i⊂`S^1$
such that nonconsecutive intervals have disjoint closures,
the restriction of $H$ to~$R⨯W_{2i+1}$ factors through the inclusion $|Λ^3_0|→F$ as a (unique) map \<h_{2i+1}:R⨯W_{2i+1}→|Λ^3_0|,\>
and the restriction of $H$ to~$R⨯W_{2i}$ factors through $s:`S^1→F$ as a (unique) map \<h_{2i}:R⨯W_{2i}→`S^1.\>
The maps $h_i$ are uniquely defined because the maps $|Λ^3_0|→F$ and $`S^1→F$ are monomorphisms.
\ppar
Having made a choice of $R$ and $\{W_i\}_i$ (and thus also $\{h_i\}_i$) for every $r∈`R$,
since the interval $[0,1]$ is compact and the intervals $R$ cover $[0,1]$,
we pick finitely many $r∈`R$ so that the corresponding intervals~$R$ cover $[0,1]$
and therefore the finite family of open sets $R⨯W_i⊂`R⨯`S^1$ constructed above covers $[0,1]⨯`S^1$.
\ppar
The remainder of the proof analyzes the maps $h_{2i}:R⨯W_{2i}→`S^1$.
For a generic point $c∈\cset$ we will define an appropriate version of a {\it local degree\/} of the collection of maps~$h_{2i}$ at~$c$.
We will then show that the local degree is independent of the parameter $r∈`R$.
For $r=0$ the degree is~1, whereas for $r=1$ the degree is~0, which contradicts the existence of~$H$.
To define generic points, we need to exclude finitely many {\it special points\/} $c∈\cset$.
This is done in two stages: first, for every interval~$R$ we exclude a certain pair of points for every consecutive intervals $W_i$, $W_{i+1}$,
ensuring the local degree is well-defined for a fixed~$R$.
Second, for every pair of intervals $R$, $\bar R$, we exclude a pair of points for every intersection $W_i∩\bar W_j$,
ensuring the local degree does not change when switching from~$R$ to~$R'$.
\ppar
Recall (^!{properties of the bad section}) that the restriction of~$s$ to any open neighborhood of any $c∈\cset$ does not factor through $|Λ^3_0|→F$.
Therefore, if for a map $f:U→`S^1$ the composition $sf:U→F$ factors through $|Λ^3_0|→F$,
then the image of~$f$ may not be an open neighborhood of any $c∈\cset$.
Thus, locally on~$U$ the map~$f$ must factor through a closed interval $[u,v]⊂`S^1$ such that $(u,v)⊂`S^1∖\cset$.
(In particular, $f$ can be constant.)
If $U$ is connected, this description is valid globally on~$U$.
\ppar
We work with a fixed interval~$R$ and $(W_i,h_i)$ as defined above.
On the (connected) intersection $R⨯(W_{i}∩W_{i+1})$, the map~$H$ factors through both $|Λ^3_0|→F$ and $s:`S^1→F$.
Therefore, for every~$i$ the maps \<h_{2i}|_{R⨯(W_{2i}∩W_{2i-1})},\qquad h_{2i}|_{R⨯(W_{2i}∩W_{2i+1})}\> factor through some intervals $[u^-_{2i},v^-_{2i}],[u^+_{2i},v^+_{2i}]⊂`S^1$, where $(u^±_{2i},v^±_{2i})⊂`S^1∖\cset$.
We refer to $u^±_{2i}$, $v^±_{2i}$ as {\it special points}.
There are only finitely many special points in~$`S^1$ since there are only finitely many choices for $R$ and~$i$.
The set of special points will be further enlarged below, when we discuss the independence of local degree from the choice of~$R$.
\ppar
Given $r∈R$, denote by $h_{i,r}$ the restriction of~$h_i$ to $W_i≅\{r\}⨯W_i⊂R⨯W_i$. 
A generic point $p∈`S^1$
is a regular value of the maps~$h_{2i,r}:W_{2i}→`S^1$.
In particular, the local degree of the map~$h_{2i,r}$ at~$p$
is well-defined and can be computed as the difference between
the number of points $a∈W_{2i}$ such that $h_{2i,r}(a)=p$ and $h_{2i,r}'(a)>0$
and the number of points $b∈W_{2i}$ such that $h_{2i,r}(b)=p$ and $h_{2i,r}'(b)<0$.
\ppar
Given a nonspecial point $c∈\cset$,
we can choose an open interval $U⊂`S^1$ that contains~$c$ and is disjoint from all intervals $[u^±_{2i},v^±_{2i}]$ associated to the given interval~$R$.
For a generic point $p∈U$, the local degree of $h_{2i,r}$ at~$p$ is independent of~$p$
because the restrictions $h_{2i,r}|_{W_{2i±1}∩W_{2i}}$ factor through the intervals $[u^±_{2i},v^±_{2i}]$ as described above, and the latter intervals are disjoint from~$U$.
We refer to the resulting common local degree as the {\it local degree of $h_{2i,r}$ at~$c$}.
For the same reason, the local degree of $h_{2i,r}$ at $c∈\cset$ is independent of the choice of $r∈R$,
so we refer to it as the {\it local degree of $h_{2i}$ at~$c$}, where the interval~$R$ is implied in the notation.
Finally, taking the sum over all~$i$, we talk about the {\it local degree of~$H$ at~$c$}, where $R$ is implicit again.
\ppar
Next, we analyze the dependence of the local degree of~$H$ at a nonspecial point $c∈\cset$ on the interval~$R$.
Suppose $r∈`R$ satisfies $r∈R$ and $r∈\bar R$ for some previously constructed intervals $R$ and~$\bar R$ together with open intervals $\{W_i\}_i$, $\{\bar W_j\}_j$.
Set \<W_{[0]}=⋃_i W_{2i},\qquad W_{[1]}=⋃_i W_{2i+1},\qquad \bar W_{[0]}=⋃_j \bar W_{2j},\qquad \bar W_{[1]}=⋃_j \bar W_{2j+1}.\>
Consider the open subset
\<M=(W_{[1]}∩\bar W_{[0]}) ∪ (W_{[0]}∩\bar W_{[1]}),\>
which is a disjoint union of finitely many open intervals $I_k⊂`S^1$.
By construction, the restriction of~$H$ to the product $(R∩\bar R)⨯M$ factors through the maps $|Λ^3_0|→F$ and $s:`S^1→F$.
Thus, the restriction of $H$ to every $(R∩\bar R)⨯I_k$ factors through some interval $[u,v]⊂`S^1$ such that $(u,v)⊂`S^1∖\cset$.
Since there are only finitely many intervals~$R$ and, therefore, finitely many intervals~$I_k$ for all pairs $R$ and~$\bar R$,
we can retroactively add the endpoints $u$ and~$v$ constructed above to the finite list of special points.
From now on we use the resulting more restrictive notion of a nonspecial point, assuming $c$ to be such a nonspecial point.
Furthermore, we choose the open interval~$U$ around~$c$ to be disjoint also from all the newly constructed intervals~$[u,v]$.
\ppar
The adjustments made to the list of special points and to the open interval~$U$ guarantee that
the local degree of~$H|_{(R∩\bar R)⨯M}$ at a generic point $p∈U$ vanishes.
Since \<W_{[0]}∪M=(W_{[0]}∩\bar W_{[0]})∪M=\bar W_{[0]}∪M,\>
the local degree of~$H|_{(R∩\bar R)⨯W_{[0]}}$ at~$c$ coincides with the local degree of~$H|_{(R∩\bar R)⨯\bar W_{[0]}}$ at~$c$,
which shows that the local degrees of~$H$ at~$c$ computed for the intervals $R$ and $\bar R$ are equal.
\ppar
Thus, given some $r∈`S^1$, every nonspecial point $c∈\cset$ has a well-defined local degree that does not depend on the interval $R$ that contains the point $r∈`S^1$.
Previously, we also proved that for any fixed interval~$R$ the local degree of~$c$ does not depend on~$r∈R$.
Since the intervals~$R$ cover the interval $[0,1]⊂`R$, the local degree of a nonspecial point $c∈\cset$ is independent of the choice of an interval~$R$ as well as a point $r∈R$.
\ppar
If the interval~$R$ contains~0, the local degree of all nonspecial points is~1 since $H|_{\{0\}⨯`S^1}=s$.
On the other hand, if the interval~$R$ contains~1, the local degree of all nonspecial points is~0 because $H|_{\{1\}⨯`S^1}$ factors through a constant map.
The resulting contradiction shows that the map~$H$ does not exist.

\proclaim Theorem.
^^={diffeological model structure}
The category $@Diffeo$ does not admit a model structure ^{transferred} (^!{transferred})
along the right adjoint functor $"SmSing: @Diffeo→@sSet$ (^!{smooth singular simplicial set})
from the ^{Kan–Quillen model structure} on the category $@sSet$.

\proof Proof.
The simplicial horn inclusion $Λ^3_0→Δ^3$ is an acyclic cofibration in $@sSet$.
By ^!{examples of concrete realizations}, $|Λ^3_0|≅\ash\concr‖Λ^3_0‖$ and $|Δ^3|≅\ash\concr‖Δ^3‖$.
Therefore, the smooth horn inclusion $|Λ^3_0|→|Δ^3|$
is an acyclic cofibration in the transferred model structure on $@Diffeo$, if it exists.
Thus, any cobase change of $|Λ^3_0|→|Δ^3|$ must be a weak equivalence in $@Diffeo$.
By ^!{smooth horn cobase change}, the cobase change of $|Λ^3_0|→|Δ^3|$ along the map $|Λ^3_0|→F$ constructed in ^!{bad cobase change}
is not a weak equivalence in $@Diffeo$, contradicting the existence of the transferred model structure on $@Diffeo$.

\proclaim Example.
The pushout of $|Δ^3|←|Λ^3_0|→F$ considered in ^!{smooth horn cobase change}
is a span of monomorphisms of ^{diffeological spaces} whose pushout in the category of ^{smooth sets}
is not a ^{diffeological space}, resolving in the negative a conjecture of Clough [\CCGT, Proposition 5.2.15].
Indeed, in the pushout in $@SmSet$ the section~$s:`S^1→F$ that was manually added to $|Λ^3_0|$ to form~$F$ is different from the section $S:`S^1→|Δ^3|$.
Since $s$ and $S$ have the same map of underlying sets, this proves that the pushout is not a ^{concrete sheaf}.

\section The projective model structure on smooth sets
^^={start of part 2}

In this section we prove that the categories $@PreSmSet$ and $@SmSet$
admit model structures transferred along the right adjoint functor $"SmSing$
and prove that $"SmSing$ is a right Quillen equivalence in both cases.

\proclaim Definition.
(Grothendieck; Batanin–Berger [\HTAPM, Definition~1.1].)
A morphism $f:X→Y$ in a relative category~$C$
is an ^={h-cofibration[|s]} if the cobase change functor
\<f_!: X/C→Y/C\>
preserves weak equivalences.

A model category is left proper if and only if all cofibrations are h-cofibrations
and in a left proper model category, cobase changes along h-cofibrations are homotopy cobase changes.
See, for example, Pavlov–Scholbach\@ [\PS, Definition~2.3] and references therein for more information.

\proclaim Proposition.
In the relative categories $@PreSmSet$ and $@SmSet$ (^!{weak equivalences of smooth sets}), all ^=:{monomorphisms are h-cofibrations}.
Furthermore, the functor $"SmSing$ reflects h-cofibrations.

\proof Proof.
In the relative category $@PreSmSet$, monomorphisms are h-cofibrations because the functor $"SmSing$ preserves colimits, monomorphisms, and weak equivalences,
so the image under $"SmSing$ of the diagram of pushout squares
$$\cd{
X&\mapright{}&A&\mapright{w}&B\cr
\mapdown{f}&&\mapdown{}&&\mapdown{}\cr
Y&\mapright{}&A'&\mapright{w'}&B',\cr
}$$ where $f$ is a monomorphism and $w$ is a weak equivalence,
is a diagram of pushout squares in $@sSet$, where the image of~$f$ is a monomorphism and the image of~$w$ is a weak equivalence.
Thus, the image of $w'$ is a weak equivalence of simplicial sets, hence the map~$w'$ is a weak equivalence.
Since the functor $"SmSing$ preserves and reflects weak equivalences, it reflects h-cofibrations.
\ppar
Applying ^!{sheafification is a weak equivalence}, we deduce that in $@SmSet$ all monomorphisms are h-cofibrations and $"SmSing$ reflects h-cofibrations.

\proclaim Proposition.
^^={weak equivalences of presheaves are closed under filtered colimits}
Weak equivalences (^!{weak equivalence of smooth sets}) in $@PreSmSet$ and $@SmSet$
are closed under filtered colimits, hence also transfinite compositions.

\proof Proof.
For $@PreSmSet$ this holds because $"SmSing$ preserves colimits
and weak equivalences of simplicial sets are closed under filtered colimits (^!{Kan–Quillen model structure on simplicial sets}).
For $@SmSet$ we use ^!{homotopy colimits of presheaves} to reduce to the previous case.

The following theorem establishes the transferred model structures on $@PreSmSet$ and $@SmSet$.
Model structures on $@SmSet$ were constructed by Cisinski [\THT, Thé\-orème 3.9] and Clough [\CCGT, Proposition~6.1.4], see the proof for details.

\proclaim Theorem.
^^={existence of model structures on presheaves}
The categories $@PreSmSet$ (^!{presmooth set}) and $@SmSet$ (^!{smooth set})
admit left proper combinatorial model structures
transferred (^!{transferred model structure}) via the ^{smooth singular simplicial set} functor $"SmSing$ (^!{smooth singular simplicial set})
from the ^{Kan–Quillen model structure on simplicial sets} (^!{Kan–Quillen model structure on simplicial sets}).
The ^{associated sheaf} functor $\ash:@PreSmSet→@SmSet$ is a left Quillen equivalence.

\proof Proof.
The mere existence of transferred model structure on $@SmSet$ is a special case of the Smith recognition theorem (Barwick [\LR, Proposition~1.7 (arXiv); 2.2 (journal)], Beke [\SHMC, Theorem~1.7]).
The existence of the model structure of Cisinski [\THT, Théorème~3.9] proves all conditions in the Smith theorem except $\inj(|I|)⊂W$,
where $I$ is the set of simplicial ^{boundary inclusions} (^!{boundary inclusions}).
By adjunction $|{-}|⊣"SmSing$, the condition $\inj(|I|)⊂W$ is equivalent to $\inj(I)⊂"SmSing(W)$, which holds by
Berwick-Evans–Boavida–Pavlov\@ [\BBP, Proposition~1.3].
Clough [\CCGT, Proposition~6.1.4] shows that $@SmSet$ admits a model structure with the same weak equivalences as Cisinski [\THT, Théorème~3.9]
and $|I|$ as generating cofibrations (replace the reference to Proposition~3.4.3 there with a reference to Crans [\QCMSS, Theorem~3.3] or Hirschhorn [\MCL, Theorem~11.3.2]).
Combined with Berwick-Evans–Boavida–Pavlov\@ [\BBP, Proposition~1.3], which shows that Cisinski's weak equivalences coincide with weak equivalences transferred along $"SmSing$,
this yields another proof of the existence of the transferred model structure on $@SmSet$.
More recently, a revised version of this argument has appeared in Clough [\HTDS, Proposition~7.1.5].
\typesetpar
Below, we give self-contained proofs of the existence of the model structures on $@SmSet$ and $@PreSmSet$ that do not rely on Cisinski's result and obviate the need to compare
our definition of weak equivalences in $@SmSet$ to Cisinski's (whose equivalence is established by Berwick-Evans–Boavida–Pavlov\@ [\BBP, Proposition~1.3]).
\ppar
By ^!{transfer theorem}, the ^{transferred model structure} on $@PreSmSet$ exists if and only if
the functor~$"SmSing$ sends transfinite compositions of cobase changes of elements of $‖J‖$ (^!{smooth realization})
to weak equivalences in~$@sSet$, where $J$ denotes the set of simplicial ^{horn inclusions} (^!{horn inclusions}).
Cobase changes of elements of $‖J‖$ in $@PreSmSet$ are weak equivalences because $"SmSing$ preserves colimits and weak equivalences,
and the simplicial map $"SmSing(‖λ_{n,k}‖)$ is a simplicial homotopy equivalence.
Cobase changes of elements of $|J|$ in $@SmSet$ by ^!{homotopy colimits of presheaves}, which reduces the problem
to the case of $@PreSmSet$, since the associated sheaf functor sends $‖λ_{n,k}‖$ to $|λ_{n,k}|$.
By ^!{weak equivalences of presheaves are closed under filtered colimits}, weak equivalences in $@PreSmSet$ are closed under transfinite compositions, completing the proof in the case of $@PreSmSet$.
The same argument establishes the case of $@SmSet$ using $|J|$ instead of $‖J‖$ and invoking ^!{sheafification is a weak equivalence}.
A model category is left proper if and only if all cofibrations are h-cofibrations.
All cofibrations are monomorphisms by construction and all monomorphisms are h-cofibrations by ^!{monomorphisms are h-cofibrations}.
\ppar
An alternative proof could be given using ^!{recognition theorem}.
The class of weak equivalences satisfies the desired properties by ^!{weak equivalences of presheaves are closed under filtered colimits}
and Makkai–Paré [\AccCat, Theorem~5.1.6] combined with the combinatoriality of the Kan–Quillen model structure (^!{Kan–Quillen model structure on simplicial sets}).
Morphisms with the right lifting property with respect to $‖I‖$ (respectively $|I|$)
are weak equivalences by adjunction $‖{-}‖⊣"SmSing$ (respectively $|{-}|⊣"SmSing$).
Finally, elements of $‖I‖$ (respectively $|I|$)
are ^{h-cofibrations} by ^!{monomorphisms are h-cofibrations}.
\typesetpar
Christensen–Wu [\CW, Proposition~4.24] observed that the relative category $@Diffeo$ is right proper
for trivial reasons: the functor $"SmSing$ preserves pullback squares, and the Kan–Quillen model structure
on simplicial sets is right proper.
The same argument shows the right properness of relative categories $@SmSet$, $@PreSmSet$, and $@PreDiffeo$.
\typesetpar
The ^{associated sheaf} functor~$\ash$ sends the generating (acyclic) cofibrations of $@PreSmSet$ to those of $@SmSet$, therefore is a left Quillen functor.
The functor~$\ash$ and its right adjoint functor (the inclusion $@SmSet→@PreSmSet$) both
preserve and reflect weak equivalences.
Furthermore, the unit map is a weak equivalence by ^!{sheafification is a weak equivalence} and the counit map is an isomorphism.
Thus, the ^{associated sheaf} functor is a left Quillen equivalence.

\proclaim Corollary.
^^={existence of intermediate model structures on presheaves}
The categories $@PreSmSet$ (^!{presmooth set}) and $@SmSet$ (^!{smooth set})
admit left proper combinatorial model structures
whose weak equivalences coincide with that of ^!{existence of model structures on presheaves}
and the set of generating cofibrations is given by an arbitrary set of monomorphisms that contains the generating cofibrations of ^!{existence of model structures on presheaves}.
The resulting model structures are Quillen equivalent to those of ^!{existence of model structures on presheaves}.

\proof Proof.
Combine ^!{existence of model structures on presheaves} with ^!{monomorphisms are h-cofibrations} and ^!{h-recognition theorem}.

\proclaim Remark.
Using the class of monomorphisms as generating cofibrations (which is generated by a set),
^!{existence of intermediate model structures on presheaves}
implies the existence of a model structure on $@SmSet$ with the same weak equivalences as in ^!{existence of model structures on presheaves}
and monomorphisms as cofibrations.
This recovers the model structure on $@SmSet$ constructed by Cisinski [\THT, Théorème~3.9],
taking the class of weak equivalences of ^!{weak equivalences of smooth sets}.

\proclaim Remark.
The proof of ^!{existence of intermediate model structures on presheaves}
also gives a new proof of the existence of
Kihara's model structure on diffeological spaces
(Kihara [\Ki, Theorem~1.3]).
Indeed, set the set~$I$ of generating cofibrations
to the set $\{|δ_n|_K\mid n≥0\}$ of realizations of simplicial boundary inclusions
with respect to Kihara's cosimplicial object (Kihara [\Ki, Definition~1.2]).
Since elements of $I$ are monomorphisms,
to show that the transferred model structure exists,
it suffices to prove that morphisms with the right lifting property with respect to $|I|_K$
are weak equivalences,
which is shown in Kihara [\Ki, Lemma~9.6.(2)].

\proclaim Theorem.
^^={Quillen equivalence for presheaves}
The Quillen adjunctions of ^!{existence of model structures on presheaves} between the model category $@sSet$
(^!{Kan–Quillen model structure on simplicial sets})
and the model categories
$@PreSmSet$ (^!{presmooth sets})
or
$@SmSet$ (^!{smooth sets})
are Quillen equivalences, in fact, weak monoidal Quillen equivalences in the sense of Schwede–Shipley [\EMMC, Definition~3.6].
The relative adjunctions between the relative category $@sSet$
and the relative categories $@PreDiffeo$ (^!{prediffeological spaces})
and $@Diffeo$ (^!{diffeological spaces}) are Dwyer–Kan equivalences of relative categories.

\proof Proof.
We give a proof for all four adjunctions simultaneously.
It suffices to show that the unit maps are weak equivalences.
Indeed, the functor $"SmSing$ reflects weak equivalences,
which implies that the left adjoint preserves weak equivalences and the triangle identity shows that the counit maps are weak equivalences.
Thus, both adjoints preserve weak equivalences and the unit and counit maps are weak equivalences, completing the proof.
The functor $"SmSing$ preserves colimits in the category $@PreSmSet$, so the unit map of $X∈@sSet$ is cocontinuous in~$X$.
Since weak equivalences in $@sSet$ are closed under filtered colimits,
we can present $X$ as a transfinite composition of cobase changes of ^{boundary inclusions} (^!{boundary inclusions})
and reduce the problem to the following elementary step:
if $X→Y$ is a cobase change of a ^{boundary inclusion}
and the unit map of~$X$ is a weak equivalence, then so is the unit map of~$Y$.
Specializing to the adjunction for $@PreSmSet$, we have a natural transformation
$$\sqcd{
∂Δ^n&\mapright{}&X\cr
\mapdown{}\;\;&&\mapdown{}\cr
Δ^n&\mapright{}&Y\cr
}
\quad⟹\quad
\sqcd{
"SmSing‖∂Δ^n‖&\mapright{}&"SmSing‖X‖\cr
\mapdown{}\hskip.4in&&\hskip.3in\mapdown{}\cr
"SmSing‖Δ^n‖&\mapright{}&"SmSing‖Y‖\cr
}
$$
of corresponding pushout squares.
The component \<X→"SmSing‖X‖\> is a weak equivalence by assumption.
The component \<Δ^n→"SmSing‖Δ^n‖\> is a weak equivalence because its source and target are contractible.
The component \<∂Δ^n→"SmSing‖∂Δ^n‖\> is a weak equivalence by inductive assumption (prove the claim by induction on the dimension of~$X$).
The left maps are monomorphisms, hence both squares are homotopy pushout squares in $@sSet$
and the component \<Y→"SmSing‖Y‖\> is a weak equivalence.
The argument for $@SmSet$, $@PreDiffeo$, and $@Diffeo$ is analogous, replacing $‖{-}‖$ with $|{-}|$, $\concr‖{-}‖$, and $\ash\concr‖{-}‖$, respectively.
\typesetpar
Finally, to show that the established Quillen equivalences
are ^{weak monoidal Quillen equivalences}
(^!{weak monoidal Quillen equivalences}),
observe that passing to adjoint maps preserves weak equivalences
because the unit and counit maps are weak equivalences.
The comonoidal map is adjoint to the product of two unit maps,
which are weak equivalences.

\proclaim Remark.
^^={surrogates}
Kihara [\KiQ, Theorem~1.1.(1)] establishes a Quillen equivalence
between simplicial sets and diffeological spaces equipped with the model structure
constructed in Kihara [\Ki, Theorem~1.3], which shows that an analogue
of the second part of ^!{Quillen equivalence for presheaves}
holds for the singular complex functor
associated to Kihara's cosimplicial diffeological space (Kihara [\Ki, Definition~1.2]).
Kihara's cosimplicial diffeological space embeds into the standard cosimplicial diffeological space
(Kihara [\Ki, Lemma~3.1]),
and this embedding induces a natural transformation between the corresponding smooth realization functors.
The cube lemma (Hovey [\MC, Lemma~5.2.6]) then shows this natural transformation to be a weak equivalence.
This provides an alternative proof of Kihara [\KiQ, Theorem~1.1.(1)].

\section The projective model structure is cartesian

We start by recalling the notion of a {\it semisimplicial set}.

\proclaim Definition.
Denote by $\sicatinj$ the subcategory of~$\sicat$ given by the same objects and injective maps of finite nonempty ordered sets.
Denote by $\sSetinj$ the subcategory of~$\sSet$ given by the essential image of the left adjoint of the restriction functor
\<\sSet=\Fun(\sicat^\op,\Set)→\Fun(\sicatinj^\op,\Set).\>

\proclaim Remark.
The left adjoint functor \<\Fun(\sicatinj^\op,\Set)→\sSet\> is faithful,
so we have an equivalence of categories \<\Fun(\sicatinj^\op,\Set)→\sSetinj.\>
Objects and morphisms in the category $\Fun(\sicatinj^\op,\Set)$ are known as ^={semisimplicial set[s|]} and ^={semisimplicial map[s|]} respectively.
Objects in $\sSetinj$ are precisely those simplicial sets for which face maps preserve nondegenerate simplices.
Morphisms in $\sSetinj$ are precisely those simplicial maps that preserve nondegenerate simplices.

\proclaim Remark.
^^={cocontinuous extension}
If $D$ is a cocomplete category,
the restriction functor along the Yoneda embedding
\<\Fun(\sSetinj^\op,D)→\Fun(\sicatinj^\op,D)\>
becomes an equivalence of categories
if we take the full subcategory of cocontinuous functors on the left side.
Likewise, the restriction functor
\<\Fun(\sSetinj^\op⨯\sSetinj^\op,D)→\Fun(\sicatinj^\op⨯\sicatinj^\op,D)\>
becomes an equivalence of categories
if on the left side we take the full subcategory of functors that are separately cocontinuous in each variable.
We use this observation to construct functors of the form $\sSetinj^\op⨯\sSetinj^\op→D$ and natural transformations between them.

\proclaim Definition.
^^={weird product}
The functor
\<⊙:\sSetinj⨯\sSetinj→\sSetinj,\qquad (K,L)↦K⊙L\>
is defined as the separately cocontinuous extension (^!{cocontinuous extension}) of the product functor
\<⊙:\sicatinj⨯\sicatinj→\sSetinj, \qquad ([m],[n])↦Δ^{[m]\lex[n]}.\>
Here $\lex$ denotes the ordinary product of finite sets with the lexicographic order.
This construction is manifestly functorial with respect to injective maps of simplices.

\proclaim Remark.
To better understand the natural of the functor $⊙$,
observe that there is a natural weak equivalence
\<⨯→⊙:\sSetinj⨯\sSetinj→\sSetinj,\> given by sending the pair \<(K,L)↦(K⨯L→K⊙L).\>
We do not need this claim later, but details of the proof can be found in Version~1 on arXiv.

\proclaim Proposition.
^^={products are retracts}
Denote by $@C$ the category $@SmSet$.
Recall the functors $‖{-}‖$ and $|{-}|$ (^!{smooth realization}).
The functor
\<|{-}|⨯|{-}|:\sSetinj⨯\sSetinj→@C,\qquad (K,L)↦|K|⨯|L|\>
is a retract of the functor
\<|{-}⊙{-}|:\sSetinj⨯\sSetinj→@C,\qquad (K,L)↦|K⊙L|.\>
The same is true for the category $@C=@PreSmSet$,
with the functor $|{-}|$ replaced by $‖{-}‖$.

\proof Proof.
By ^!{cocontinuous extension}, it suffices to exhibit the functor
\<|{-}|⨯|{-}|:\sicatinj⨯\sicatinj→@C,\qquad\allowbreak (K,L)↦|K|⨯|L|\>
as a retract of the functor
\<|{-}⊙{-}|:\sicatinj⨯\sicatinj→@C,\> which sends \<(K,L)↦|K⊙L|.\>
The natural inclusion \<ι:`Δ^m⨯`Δ^n→`Δ^{[m]\lex[n]}\>
sends $$(x_0,…,x_m,\fbreak y_0,…,y_n)↦(x_0y_0,x_0y_1,…,x_0y_n,x_1y_0,…,x_1y_n,…,x_my_0,…,x_my_m).$$
The natural retraction \<ρ:`Δ^{[m]\lex[n]}→`Δ^m⨯`Δ^n\>
sends
$(z_{0,0},…,z_{m,n})$ to the point $$(z_{0,0}+⋯+z_{0,n},…,z_{m,0}+⋯+z_{m,n},\fbreak z_{0,0}+⋯+z_{m,0},…,z_{0,n}+⋯+z_{m,n}).$$
The composition $ρι$ is the identity map by construction.

\proclaim Proposition.
^^={pushout product of cofibrations}
Given $m≥0$, $n≥0$,
the pushout product \<p:P→‖Δ^m‖⨯‖Δ^n‖\> of the maps \<‖δ_m‖:‖∂Δ^m‖→‖Δ^m‖, \qquad ‖δ_n‖:‖∂Δ^n‖→‖Δ^n‖\>
(^!{boundary inclusion}, ^!{smooth realization}) is a cofibration in $@PreSmSet$
(^!{existence of model structures on presheaves}).
Likewise,
the pushout product \<p:P→|Δ^m|⨯|Δ^n|\> of the maps \<|δ_m|:|∂Δ^m|→|Δ^m|, \qquad |δ_n|:|∂Δ^n|→|Δ^n|\>
(^!{boundary inclusion}, ^!{smooth realization}) is a cofibration in $@SmSet$
(^!{existence of model structures on presheaves}).

\proof Proof.
We can apply ^!{products are retracts}, since the involved maps are morphisms in $\sSetinj$.
Consider the simplicial map~$q$ given by the pushout product of $∂Δ^m→Δ^m$ and $∂Δ^n→Δ^n$ with respect to the operation~$⊙$ of ^!{weird product}.
The operation~$⊙$ preserves colimits in each argument, so every simplex $σ:Δ^k→A⊙B$ ($A,B∈@sSet$) factors through the map $a⊙b:Δ^m⊙Δ^n→A⊙B$
for some $a:Δ^m→A$, $b:Δ^n→B$.
If we require that $σ$ does not factor through the maps $a'⊙b$ or $a⊙b'$ induced by a proper face~$a'$ of~$a$ or $b'$ of~$b$,
then the pair $(a,b)$ is uniquely determined by~$σ$.
Thus, if the map $q$ sends two simplices in its domain to the same simplex in its codomain, both simplices must have the same pair $(a,b)$.
In particular, they must come from the same summand in the pushout and therefore must be equal as simplices of that summand.
Therefore, the map~$q$ is a monomorphism, i.e., a cofibration of simplicial sets.
\ppar
The natural retraction defined in ^!{products are retracts} exhibits $p$ as a retract of~$‖q‖$ respectively~$|q|$.
Since $‖q‖$ respectively $|q|$ is a cofibration, so is $p$.

\proclaim Proposition.
^^={pushout product of a cofibration and acyclic cofibration}
Given $m>0$, $0≤k≤m$, $n≥0$,
the pushout product of the maps \<|λ_{m,k}|:\fbreak |Λ^m_k|→|Δ^m|, \qquad |δ_n|:|∂Δ^n|→|Δ^n|\>
(^!{horn inclusion}, ^!{smooth realization}) is a weak equivalence in $@SmSet$
(^!{weak equivalence of smooth sets}).
The same is true for the category $@PreSmSet$,
with the functor $|{-}|$ replaced by $‖{-}‖$.

\proof Proof.
The inclusion of the apex $|Δ^0|→|Λ^n_k|$ is a ^{smooth homotopy equivalence} (^!{smooth homotopy equivalence}).
Therefore, its pushout product with $|δ_n|:|∂Δ^n|→|Δ^n|$ is also a smooth homotopy equivalence.
Smooth homotopy equivalences are weak equivalences, completing the proof.
The case of $‖{-}‖$ is treated in the same way.

The following result implies (as a special case)
an affirmative answer to a conjecture of Christensen–Wu [\CW, Proposition~4.38]:
the internal hom from a cofibrant diffeological space to a fibrant diffeological space
is a fibrant diffeological space.

\proclaim Proposition.
^^={model structures on presheaves are cartesian}
The model structures of ^!{existence of model structures on presheaves} are ^{cartesian} model structures (^!{cartesian}).

\proof Proof.
The same proof works for
both
model categories.
By ^!{pushout product of cofibrations}, the pushout product of generating cofibrations is a cofibration.
Thus, the pushout product of cofibrations is a cofibration.
By ^!{pushout product of a cofibration and acyclic cofibration},
the pushout product of a generating cofibration and a generating acyclic cofibration
is a weak equivalence.
Since it is also a cofibration, it must be an acyclic cofibration.
Therefore, the pushout product of a cofibration and an acyclic cofibration is an acyclic cofibration.
Finally, the terminal object (given by a point) is cofibrant.

\proclaim Proposition.
The categories $@PreSmSet$ (^!{presmooth set}) and $@SmSet$ (^!{smooth set})
admit cartesian left proper combinatorial model structures
whose weak equivalences coincide with that of ^!{existence of model structures on presheaves}
and the set of generating cofibrations is given by an arbitrary set of monomorphisms that is closed under pushout products
and contains the generating cofibrations of ^!{existence of model structures on presheaves}.

\proof Proof.
Combine ^!{existence of intermediate model structures on presheaves} with the fact that the pushout product axiom can be checked on generating cofibrations.

\section Cofibrancy of manifolds

By Christensen–Wu [\CW, Corollary~4.36], every manifold is fibrant in $@PreSmSet$ and $@SmSet$.
In this section, we show that every manifold is cofibrant in $@SmSet$, resolving in the affirmative (^!{manifolds are cofibrant}) a conjecture of Christensen–Wu [\CW, §4.2].

\proclaim Proposition.
^^={blur}
(Berwick-Evans–Boavida–Pavlov\@ [\BBP, Proposition~4.17].)
For any simplicial set~$K$ and a rectilinear (hence smooth) triangulation $j:|K|→U$
of an open subset $U⊂`R^n$ ($n≥0$),
we can find a morphism $r:U→|K|$ of ^{smooth sets} with the following properties.
\typesetli
\li The map $r$ collapses an open neighborhood~$U_σ$ of every closed simplex~$σ$
(given by taking $x_i≥0$ in ^!{extended smooth simplex})
in the triangulation~$j$ to~$σ$.
\li There is a smooth homotopy $h:`Δ^1⨯U→U$ from the identity map on~$U$ to~$jr$.
This homotopy preserves the image of every closed simplex in~$|K|$.
\li The smooth homotopy $h$ restricts to a smooth homotopy $`Δ^1⨯|K|→|K|$
from the identity map on~$|K|$ to $rj$.
This homotopy preserves every closed simplex in~$|K|$.

\proclaim Proposition.
^^={manifolds are cofibrant}
Any (paracompact Hausdorff) smooth manifold is cofibrant in the model category $@SmSet$
(^!{existence of model structures on presheaves}).

\proof Proof.
Coproducts of cofibrant objects are cofibrant, so we can assume the manifold to be connected, hence second countable.
Any second countable Hausdorff manifold is a retract of a tubular neighborhood of the image of its embedding into some $`R^n$.
Thus, it remains to treat the case when $M$ is an open subset of $`R^n$.
\ppar
Pick smooth functions $f_1,…,f_n:M→(0,∞)$
such that for every~$i$ the vector field $e_if_i$ has an everywhere defined flow $a_i:`R⨯M→M$,
where $e_i$ are elements of the standard basis of~$`R^n$.
The various $a_i$ combine into a smooth map $a:`R^n⨯M→M$ that sends a point $(t,x)$ to \<a_n(t_n,a_{n-1}(t_{n-1},… a_1(t_1,x)…)).\>
For any $m∈M$ the map $b_m=a(-,m):`R^n→M$ is an open embedding that sends $0$ to~$m$.
In particular, the map \<b_m^{-1}:D_m→`R^n\> is well defined, with its domain~$D_m$ being the open subset of~$M$ given by the image of~$b_m$,
so that \<a(b_m^{-1}(x),m)=x\> for all $x∈D_m$.
The maps $b_m^{-1}$ combine into the smooth map \<c:D→`R^n,\qquad (x,m)↦b_m^{-1}(x),\>
whose source \<D=\{(x,m)∈M⨯M\mid x∈D_m\}\> is an open subset of $M⨯M$.
We have $a(c(x,m),m)=x$ for all $(x,m)∈D$.
\ppar
Pick a rectilinear triangulation~$K$ of~$M$,
with the induced map $ι:|K|→M$.
(Since $M$ is an open subset of~$`R^n$,
such a triangulation can be constructed in an elementary fashion without using the full strength of the triangulation theorem for smooth manifolds.)
We now exhibit~$M$ as a retract of $`Δ^n⨯|K|$.
The latter object is cofibrant by ^!{model structures on presheaves are cartesian}, which implies that $M$ is also cofibrant.
\ppar
\begingroup
Using ^!{blur}, pick a map $α:M→|K|$ with the following properties.
\typesetli
\li Given a simplex~$σ$ in~$K$, consider its associated map $ι:`Δ^k→M$.
Denote by $V_σ⊂M$ the $ι$-image of the closed simplex~$`Δ_c^k⊂`Δ^k$,
given by the subpresheaf of $`Δ^k$ (^!{extended smooth simplex}) with coordinates $x_i≥0$ for all~$i$.
We require that $α$ maps some open neighborhood~$U_σ$ of~$V_σ$ to the image of~$`Δ_c^k→`Δ^k→|K|$, where the map $`Δ^k→|K|$ is induced by~$σ$.
\li Additionally, we require that for any $m∈U_σ$ we have $m∈D_{ι(α(m))}$.
We can always shrink $U_σ$ to a smaller open neighborhood of $V_σ$ so that it satisfies this condition,
since $ι(α(m))∈V_σ$ and $V_σ$ is compact, so there is $ε>0$ such that for any $m∈V_σ$ and any $x∈M$ with $‖x-m‖<ε$ we have $x∈D_m$,
and for any $ε>0$ we can choose $α$ so that for all $m∈V_σ$ we have $‖m-ι(α(m))‖<ε$.
\ppar
\endgroup
The retraction~$r$ is given by the composition \<r:`Δ^n⨯|K|\lto6{`Δ^n⨯ι}`Δ^n⨯M\lto1{a}M.\>
Consider the inclusion \<i:M→`Δ^n⨯|K|, \qquad m↦(c(m,α(m)),ι^{-1}(α(m))).\>
By definition of $α$ we have $m∈D_{α(m)}$, so $(m,α(m))∈D$ and the first component is well defined and smooth.
The point $α(m)$ belongs to the $ι$-image of a unique interior simplex $`Δ_i^k⊂|K|$,
where $`Δ_i^k$ is the subpresheaf of $`Δ^k$ (^!{extended smooth simplex}) with coordinates $x_i>0$ for all~$i$.
Thus, the second map is well defined on individual points.
To show that it is induced by a (necessarily unique) morphism of sheaves, it suffices to observe that
for any $k$-simplex $σ∈K$ the restriction of~$i$ to~$U_σ⊂M$ is given by the composition of morphisms of sheaves
$$f_σ:U_σ\lto5{\diag}U_σ⨯U_σ\lto5{\id⨯α}U_σ⨯U_σ\lto5{(c,π_2)}`Δ^n⨯`Δ^k\lto5{\id⨯σ}`Δ^n⨯|K|.$$
The collection $\{U_σ\}_{σ∈K}$ is an open cover of~$M$
and the family $\{f_σ\}_{σ∈K}$ is compatible because it is compatible on underlying sets by construction and the sheaf $`Δ^n⨯|K|$ is concrete
because $K$ satisfies the assumptions of ^!{concrete realization}.
Thus, the compatible family $\{f_σ\}_{σ∈K}$ can be glued to a morphism of sheaves~$i$.
\typesetpar
The composition $ri:M→M$ sends $m∈M$ to \<a(c(m,α(m)),α(m))=m,\> so $ri=\id_M$ by concreteness of~$M$.

\section The smooth Oka principle for smooth sets

The following result improves on the usual way of computing derived internal homs in cartesian model categories
by eliminating the fibrant replacement functor.
The proof of a more general result (discussed in ^!{smooth Oka principle for simplicial presheaves} below)
can be found in Berwick-Evans–Boavida–Pavlov\@ [\BBP, Theorem~1.1].
The name “smooth Oka principle” was suggested by Urs Schreiber (Sati–Schreiber [\SS, Theorem~3.3.53]).

\proclaim Proposition.
^^={smooth Oka principle for presheaves}
(The smooth Oka principle for smooth sets and diffeological spaces.)
If $X$ is a smooth manifold, the functor
\<"Hom(X,-):@SmSet→@SmSet\>
preserves weak equivalences (^!{weak equivalence of smooth sets})
and therefore computes the derived internal hom in the model structure of ^!{existence of model structures on presheaves}.

The following result was already established in Berwick-Evans–Boavida–Pavlov\@ [\BBP, Proposition 2.18].
It resolves in the affirmative a conjecture of Christensen–Wu [\CW, §1].
We reproduce the proof here for the sake of completeness, adding a few more details.

\proclaim Corollary.
^^={smooth homotopy groups}
For every $X∈@SmSet$, the canonical map from the $n$th smooth homotopy group of~$X$ at point $x_0∈X$
to the $n$th simplicial homotopy group of $"SmSing X$ at point~$x_0$
is an isomorphism.
Here the $n$th smooth homotopy group of~$X$ at point $x_0∈X$
is defined as the quotient of the set of morphisms $s:`S^n→X$ that send $*∈`S^n$ to~$x_0$
modulo the equivalence relation that identifies $s∼s'$ if there is a morphism $h:`Δ^1⨯`S^n→X$
whose restriction to $`Δ^1⨯\{x_0\}$ is the constant map given by the composition $`Δ^1→`Δ^0\lto3{x_0}X$.

\proof Proof.
Recall that the simplicial homotopy group $π_n("SmSing X,x_0)$
can be computed as the set of connected components of the homotopy fiber of the map
of derived mapping simplicial sets
$$\rdf "Hom("SmSing `S^n,"SmSing X)→\rdf "Hom("SmSing `Δ^0,"SmSing X).$$
By ^!{smooth Oka principle for presheaves}, the latter map is weakly equivalent to $"SmSing$ applied to the map
\<"Hom(`S^n,X)→"Hom(`Δ^0,X).\>
\ppar
The set of connected components of the homotopy fiber of the latter map can be computed as the following quotient.
Elements are morphisms $S:`S^n→X$ together with a map $P:`Δ^1→X$ that sends $1↦s(*)$ and $0↦x_0$.
The pair $(S,P)$ can be encoded as a single map $`S^n⊔_{`Δ^0}`Δ^1→X$.
We identify $(S,P)∼(S',P')$ if there is a smooth homotopy \<`Δ^1⨯(`S^n⊔_{`Δ^0}`Δ^1)→X\> between them.
The canonical map $`S^n⊔_{`Δ^0}`Δ^1→`S^n$ that projects $`Δ^1$ to $*∈`S^n$ is a smooth homotopy equivalence,
the inverse map $`S^n→`S^n⊔_{`Δ^0}`Δ^1$ is constructed by projecting a disk of small radius $ε>0$
around~$*$ to the interval $[0,1]⊂`Δ^1$ using the appropriately smoothened distance function from~$*$.
Since this smooth homotopy equivalence preserves the basepoint,
this proves that the set of connected components of the homotopy fiber is isomorphic to the $n$th smooth homotopy group of~$X$.

The following result answers a question by Sati–Schreiber [\SS, Remark~2.2.9].
We remark that the extended simplex $`Δ^1$ can be replaced with the interval $[0,1]$
in the statement below, since both simplices give rise to the same notion of concordance.
The result is applicable when $X$ is a manifold, since these are cofibrant by ^!{manifolds are cofibrant}.

\proclaim Proposition.
Suppose $P_0→X$ and $P_1→X$ are diffeological principal bundles over a cofibrant diffeological space~$X$, e.g., a smooth manifold.
Suppose $P_0→X$ and $P_1→X$ are concordant, meaning there is a diffeological principal bundle over $`Δ^1⨯X$
whose pullback to $\{i\}⨯X$ is isomorphic to $P_i→X$.
Then $P_0→X$ and $P_1→X$ are isomorphic.

\proof Proof.
As pointed out in Sati–Schreiber [\SS, Theorem~2.2.8 and Remark~2.2.9],
it suffices to show that $X→`Δ^1⨯X$ is an acyclic cofibration
and every diffeological fiber bundle is a fibration.
The former holds by ^!{model structures on presheaves are cartesian} and the latter holds by Christensen–Wu [\CW, Propositions 4.28 and 4.30].

\section Algebras over operads in smooth sets
^^={end of part 2}

In this section, we establish model structures on operads and algebras over operads in (pre)smooth sets
and compare them to the existing constructions in the simplicial and quasicategorical settings.

\proclaim Proposition.
^^={properties of model structures on presheaves}
The model categories $@PreSmSet$ and $@SmSet$ of ^!{existence of model structures on presheaves} are
h-mon\-oidal, symmetric h-monoidal, and flat (Pavlov–Scholbach [\PS, Definitions 3.2.2, 4.2.4, 3.2.4]).

\proof Proof.
For h-monoidality,
since these model structures are ^{cartesian} by ^!{model structures on presheaves are cartesian},
it suffices to show that the product of any object and an (acyclic) cofibration is an (acyclic) h-cofibration.
The nonacyclic part holds because cofibrations are monomorphisms, the product of an object and a monomorphism is a monomorphism,
and monomorphisms are h-cofibrations by ^!{monomorphisms are h-cofibrations}.
The acyclic part holds because $"SmSing$ preserves and reflects weak equivalences.
\typesetpar
For symmetric h-monoidality, the argument is the same, using the fact that $"SmSing$ preserves colimits in $@PreSmSet$.
For $@SmSet$ we need to further observe that the ^{associated sheaf} functor preserves monomorphisms and weak equivalences by ^!{sheafification is a weak equivalence}.
\ppar
Flatness in $@PreSmSet$ follows from the fact that $"SmSing$ preserves products and pushouts,
and the model category $@sSet$ is flat (Pavlov–Scholbach [\PS, \S7.1]).
Flatness in $@SmSet$ then follows from ^!{sheafification is a weak equivalence}.

Recall (Pavlov–Scholbach [\PSa, Definition~2.1] that a map $f:A→B$ in a symmetric monoidal model category is {\it flat\/}
if $f$ is a weak equivalence and the pushout product $f◻s$ is a weak equivalence for any cofibration~$s$.
In $@SmSet$ and $@PreSmSet$, flat maps coincide with weak equivalences.
Likewise, a $Σ_n$-equivariant map~$f$ is {\it symmetric flat\/} if $f◻_{Σ_n}s^{◻n}$ is a weak equivalence for any multi-index~$n$ and finite family of cofibrations~$s$.
A sufficient condition is given in Pavlov–Scholbach [\PSa, Lemma~7.6], essentially requiring the $Σ_n$-action to be projectively cofibrant.

\proclaim Proposition.
^^={admissibility of operads in presheaves}
^^={rectification of algebras in presheaves}
Suppose $O$ is a colored (symmetric) operad in $@PreSmSet$ or $@SmSet$ (^!{existence of model structures on presheaves}).
The category of algebras over~$O$ admits a model structure transferred along the forgetful functor that extracts underlying objects.
If $f:O→O'$ is a weak equivalence of colored (symmetric) operads, then it induces a Quillen equivalence of model categories of algebras over $O$ and $O'$
if and only if $f$ is a (symmetric) flat map.
(In the nonsymmetric case, flat maps coincide with weak equivalences.)

\proof Proof.
Combine ^!{properties of model structures on presheaves}
together with Pavlov–Scholbach [\PSa, Theorems 5.11, 7.5, 7.11].

\proclaim Proposition.
^^={quasicategorical algebras in presheaves}
Suppose $O$ is a $Σ$-cofibrant colored symmetric operad in the category $@PreSmSet$ or $@SmSet$ (^!{existence of model structures on presheaves}),
where an operad~$O$ is {\it $Σ$-cofibrant\/} if the unit map $1→O(a,a)$ is a cofibration for every color~$a$
and every component of~$O$ is projectively cofibrant as an object in $@PreSmSet$ or $@SmSet$ with respect to the action of the appropriate symmetric group.
Then the functor of quasicategories
$$@Alg_O(@SmSet)^c[W_O^{-1}]→@Alg_O(@SmSet[W^{-1}])$$
is an equivalence of quasicategories.
Here $@=Alg_O$ on the left denotes the category of algebras over the operad~$O$, $@Alg_O$ on the right denotes the quasicategory of quasicategorical algebras over
the operad~$O$,
the brackets $[-]$ denote quasicategorical localizations,
superscript~$c$ denotes the full subcategory of cofibrant objects,
and $W_O$ and $W$ denotes the weak equivalences with respect to the corresponding model structures.
In particular, since the quasicategory $@SmSet[W^{-1}]$ is equivalent to the underlying quasicategory of $@sSet$
by ^!{Quillen equivalence for presheaves},
the right side is equivalent to the quasicategory of algebras over the operad $"SmSing(O)$ in spaces.
All statements also hold if $@SmSet$ is replaced by $@PreSmSet$.

\proof Proof.
Combine ^!{properties of model structures on presheaves} and Haugseng [\Haug, Theorem~4.10].

\proclaim Proposition.
^^={transport of algebras in presheaves}
There is a Quillen equivalence
\<L⊣R:@=Oper_@sSet ⇄ @Oper_@PreSmSet\>
of model categories of colored symmetric operads in $@sSet$ and $@PreSmSet$ (constructed using ^!{admissibility of operads in presheaves}).
Here the right adjoint functor~$R$ applies the functor $"SmSing$ componentwise to a given operad in $@PreSmSet$.
For any cofibrant (in the model category $@Oper_@sSet$) colored symmetric simplicial operad~$O$, there is a Quillen equivalence
\<L_O⊣R_O:@Alg_O(@sSet) ⇄ @Alg_{LO}(@PreSmSet),\>
where the right adjoint functor~$R_O$ applies $"SmSing$ to components of a given algebra over~$LO$,
and equips the result with an action of~$O$ using the unit map $O→RLO$.
For any fibrant (in $@PreSmSet$) operad~$P$, there is a Quillen equivalence
\<L_P⊣R_P:@Alg_{RP}(@sSet) ⇄ @Alg_P(@PreSmSet),\>
where the right adjoint functor~$R_P$ applies $"SmSing$ to components of a given algebra over~$P$.
All statements also hold if $@PreSmSet$ is replaced by $@SmSet$.
Also, without (co)fibrancy conditions on $O$ and $P$ we still get Quillen adjunctions.

\proof Proof.
Combine Pavlov–Scholbach [\PSa, Theorem~8.10], ^!{properties of model structures on presheaves}, and ^!{Quillen equivalence for presheaves}.

\proclaim Example.
As a special case, we see that (strict) monoids in ^{smooth sets} are Quillen equivalent to simplicial monoids.
Likewise, $\Ei$-monoids in ^{smooth sets} (where $\Ei$ denotes a $Σ$-cofibrant operad in smooth sets weakly equivalent to the terminal operad)
are Quillen equivalent to $Γ$-spaces and $\Ei$-monoids in simplicial sets,
which can be seen by combining the second part of ^!{transport of algebras in presheaves}
with ^!{rectification of algebras in presheaves}.

\section Model structures on enriched presheaves
^^={start of part 3}

In this section, we extend the results obtained so far to the case of simplicial presheaves on the site $@Cart$,
i.e., simplicial objects in the categories $@PreSmSet$ and $@SmSet$.
This is of crucial importance to applications, many of which involve objects that have higher homotopy groups,
such as the stack of vector bundles with connections or the stack of bundle gerbes.

More generally, we construct a model structure on presheaves and sheaves on $@Cart$ valued in a left proper combinatorial model category~$@V$.
Its weak equivalences are precisely those morphisms $F→G$ of $@V$-valued presheaves (or sheaves) on manifolds
such that the induced map on ^{shapes} (^!{shape}) is a weak equivalences in~$@V$.

\proclaim Examples.
We have the following principal examples of left proper combinatorial model categories $@V$:
\typesetli
\li $@V=@sSet$: suitable for encoding structures such as principal $G$-bundles and higher nonabelian bundles;
\li $@V=@=Ch_{≥0}$: suitable for encoding abelian sheaf cohomology, e.g., bundle $n$-gerbes with connection;
\li $@V=@=Sp_{≥0}$: suitable for encoding extraordinary differential cohomology, e.g., differential K-theory;
\li $@V=@Ch$ and $@V=@Sp$ are also examples, although they do not satisfy the conditions of ^!{smooth Oka principle for varieties}.

\proclaim Definition.
^^={enriched preshea[f|ves]}
Given a cocomplete and complete category~$@V$,
denote by $@=PreSm_@V$ respectively $@=Sm_@V$ the category of presheaves respectively sheaves on the site $@Cart$ valued in~$@V$.
In particular, for $@V=@sSet$ objects in $@Sm_@sSet$ are simplicial objects in ^{smooth sets}, i.e., ^={simplicial smooth set[s|]}.
Denote by
\<⊗:@V⨯@Set→@V,\qquad (V,S) ↦ ∐_S V\>
the tensoring of~$@V$ over sets.
Denote by
\<⊗:@V⨯@PreSm_@V→@PreSm_@V\>
the functor sending
\<(X,F)↦(W↦X⊗F(W))\>
and by
\<⊗:@V⨯@Sm_@V→@Sm_@V\>
the functor that takes the ^{associated sheaf} (^!{associated sheaf}) of the tensoring in $@PreSm_@V$.
Then denote by
\<⊠:@V^→⨯@PreSm_V^→→@PreSm_@V^→, \qquad ⊠:@V^→⨯@Sm_V^→→@Sm_@V^→\>
the associated pushout product functors.
Objects in $@PreSmSet$ (^!{presmooth set}) and $@SmSet$ (^!{smooth set}) can be (silently) promoted to objects in $@PreSm_V$ respectively $@Sm_V$
using the cocontinuous functor $1⊗{-}:@Set→@V$, where $1$ is the terminal object in~$@V$.

\proclaim Definition.
^^={model structure on simplicial objects}
Suppose $@V$ is a left proper combinatorial model category.
Denote by $@V_@Δ$ the category of simplicial objects in~$@V$.
Turn $@V_@Δ$ into a relative category by creating its weak equivalences
using the homotopy colimit functor $@V_@Δ→@V$.
Turn $@V_@Δ$ into a model category by equipping it with the left Bousfield localization of the projective model structure
at maps of representable presheaves $Δ^n→Δ^0$ tensored with an arbitrary object of~$@V$.
(It suffices to take the set of $λ$-small objects in~$@V$ for a sufficiently large regular cardinal~$λ$.)
We also have a left Quillen equivalence $"colim:@V_@Δ→@V$, which takes the colimit of a simplicial object.
It is a ^{weak monoidal Quillen equivalence} (^!{weak monoidal Quillen equivalence}).

\proclaim Definition.
Suppose $@V$ is a left proper combinatorial model category.
Denote by
\<\strut‖{-}‖:@V_@Δ⇄@PreSm_@V:"=Sing_@V,\qquad|{-}|:@V_@Δ⇄@Sm_@V:"Sing_@V\>
the adjunctions constructed as follows.
The right adjoint $"Sing_@V$ evaluates the given presheaf on smooth simplices~$`Δ^n$.
The left adjoints send $V⊗Δ^n$ to $V⊗‖Δ^n‖$ respectively $V⊗|Δ^n|$.
^^={smooth singular bisimplicial set[|s]}
^^={smooth simplicial set[|s]}
^^={presmooth simplicial set[|s]}
Equip $@PreSm_@V$ and $@Sm_@V$ with weak equivalences created by the functor $"Sing_@V$, which turns them into relative categories.
^^={weak equivalence[|s] of smooth simplicial sets}

\proclaim Proposition.
^^={infinity-sheafification is a weak equivalence}
Given a left proper combinatorial model category~$@V$,
any Čech-local (equivalently, stalkwise) weak equivalence in $@PreSm_@V$ is a weak equivalence in $@PreSm_@V$.
As a special case, the map $F→LF$ that takes the associated sheaf of a presheaf~$F$ is a weak equivalence in $@PreSm_@V$.

\proof Proof.
(Compare ^!{sheafification is a weak equivalence}.)
Consider the model structure on the category $@PreSm_@V$ given by the injective model structure left Bousfield localized at Čech nerves of good open covers.
Consider the model structure on the category $@V_@Δ$ given by the injective model structure left Bousfield localized at maps of representable presheaves $Δ^n→Δ^0$.
(In both cases we tensor the representable presheaves with an arbitrary $λ$-small object of~$@V$, for a sufficiently large cardinal~$λ$.)
Consider the functor \<"Sing_@V:@PreSm_@V→@V_@Δ.\>
The functor~$"Sing_@V$ is a left adjoint functor that preserves injective cofibrations and injective weak equivalences.
Furthermore, by Borsuk's nerve theorem (for example, combine Weil [\TDR, §5] and Eilenberg [\SHDM, Theorem~II]),
the functor~$"Sing_@V$ sends the Čech nerve of a good open cover to a weak equivalence in~$@V_@Δ$.
Thus, $"Sing_@V$ is a left Quillen functor that preserves weak equivalences.
The map $F→G$ is a Čech-local weak equivalence by assumption.
Thus, the map $"Sing_@V F→"Sing_@V G$ is a weak equivalence in~$@V$,
therefore $F→G$ is a weak equivalence in $@PreSm_@V$.

\proclaim Remark.
^^={existence of local model structures}
Weak equivalences in $@PreSm_@V$ (and $@Sm_@V$)
are precisely the weak equivalences in the $`R$-local
projective or injective model structure
on $@V$-valued presheaves on~$@Cart$,
defined as the left Bousfield localization of the projective or injective model structure at maps $`R^n→`R^0$,
which exists by the Smith theorem (Barwick [\LR, Theorem~2.1 (arXiv); 4.7 (journal)]).

\proclaim Theorem.
^^={existence of model structures on enriched presheaves}
Given a left proper combinatorial model category~$@V$,
the categories $@PreSm_@V$ and $@Sm_@V$ (^!{enriched presheaves}) admit
left proper combinatorial model structures
whose weak equivalences are as in ^!{weak equivalences of smooth simplicial sets}
and generating cofibrations are given by the maps
$i⊠‖δ_n‖$ (respectively $i⊠|δ_n|$),
where $i$ belongs to a fixed set of generating cofibrations in~$@V$,
the map $δ_n:∂Δ^n→Δ^n$ is a simplicial boundary inclusion ($n≥0$),
and $⊠$ is defined in ^!{enriched presheaves}.
Both model structures have the following properties.
\typesetli
\li If weak equivalences in~$@V$ are closed under filtered colimits, then so are weak equivalences in $@PreSm_@V$ and $@Sm_@V$.
\li Objectwise ^{h-cofibrations} are h-cofibrations in $@PreSm_@V$ and $@Sm_@V$, and the functor $"Sing_@V$ reflects h-cofibrations.
\li (Compare ^!{properties of model structures on presheaves}.)
The model categories $@PreSm_@V$ and $@Sm_@V$ inherit from~$@V$
properties such as being
monoidal (with respect to the objectwise monoidal product), tractable, h-monoidal, and flat (Pavlov–Scholbach [\PS, Definitions 2.1, 3.2.2, and 3.2.4]),
symmetric h-monoidal (Pavlov–Scholbach [\PS, Definition 4.2.4]).

\proof Proof.
The mere existence of the model structures is a special case of the Smith recognition theorem (Barwick [\LR, Proposition~1.7 (arXiv); 2.2 (journal)], Beke [\SHMC, Theorem~1.7])
and Smith's theorem on the existence of left Bousfield localizations (Barwick [\LR, Theorem~2.1 (arXiv); 4.7 (journal)]),
used to construct the model structures of ^!{existence of local model structures}.
The existence of the local injective model structure proves all conditions in the Smith theorem except $\inj(I)⊂W$,
where $I$ is the set of generating cofibrations.
The latter condition then follows from the existence of the local projective model structure,
since all projective cofibrations are also cofibrations in the model structure under consideration.
\typesetpar
Below, we give a self-contained proof of the existence of the model structure using ^!{recognition theorem},
which does not rely on the local projective or local injective model structures.
\ppar
If weak equivalences in $@V$ are closed under filtered colimits,
then filtered colimits in $@V$ are also homotopy colimits.
Therefore, weak equivalences in $@V_@Δ$ are closed under filtered colimits
because filtered homotopy colimits commute with homotopy colimits of simplicial objects.
Since the functor $"Sing_@V$ preserves colimits,
weak equivalences in $@PreSm_@V$ are closed under filtered colimits.
For $@Sm_@V$ we use ^!{infinity-sheafification is a weak equivalence} to reduce to the case of $@PreSm_@V$.
\ppar
In the relative category $@PreSm_@V$, objectwise h-cofibrations are h-cofibrations because the functor $"Sing_@V$ preserves colimits, objectwise h-cofibrations,
and weak equivalences,
so the image under $"Sing_@V$ of the diagram of pushout squares
$$\cd{
X&\mapright{}&A&\mapright{w}&B\cr
\mapdown{f}&&\mapdown{}&&\mapdown{}\cr
Y&\mapright{}&A'&\mapright{w'}&B',\cr
}$$ where $f$ is an objectwise h-cofibration and $w$ is a weak equivalence,
is a diagram of pushout squares in $@V_@Δ$, where the image of~$f$ is an objectwise h-cofibration and the image of~$w$ is a weak equivalence.
Interpreting the resulting pushout squares in $@V_@Δ$ as a simplicial object in the category of diagrams of homotopy pushout squares in~$@V$,
its homotopy colimit is also a diagram of homotopy pushout squares in~$@V$.
Hence the map~$w'$ is a weak equivalence.
Since the functor $"Sing_@V:@PreSm_@V→@V_@Δ$ preserves and reflects weak equivalences, it reflects h-cofibrations.
Applying ^!{infinity-sheafification is a weak equivalence}, we deduce that in $@Sm_@V$ all objectwise h-cofibrations are h-cofibrations and $"Sing_@V$ reflects h-cofibrations.
\typesetpar
All generating cofibrations $i⊠‖δ_n‖$ (respectively $i⊠|δ_n|$)
are objectwise (coproducts of) cofibrations, hence also objectwise h-cofibrations by left properness of~$@V$,
therefore they are h-cofibrations.
\ppar
To show the existence of the model structure on $@PreSm_@V$ (respectively $@Sm_@V$),
we apply ^!{recognition theorem} to the set of generating cofibrations $i⊠‖δ_n‖$ (respectively $i⊠|δ_n|$).
The class of weak equivalences satisfies the desired properties because the functor $"Sing_@V$ preserves filtered colimits
and weak equivalences in $@V_@Δ$ satisfy the desired properties.
Morphisms~$f$ with the right lifting property with respect to generating cofibrations $i⊠‖δ_n‖$ (respectively $i⊠|δ_n|$)
are weak equivalences by adjunction of ^!{presmooth simplicial set},
which forces the Reedy matching maps of $"Sing_@V f$ to have the right lifting property with respect to generating cofibrations~$i$,
making them into acyclic fibrations in~$@V$.
This implies that $"Sing_@V f$ is a Reedy acyclic fibration, hence an objectwise weak equivalence, hence $f$ is a weak equivalence.
Since the generating cofibrations are h-cofibrations, this proves the existence of the model structure.
\ppar
To show that the model structures on $@PreSm_@V$ and $@Sm_@V$ are monoidal (with respect to objectwise monoidal products of presheaves)
whenever $@V$ is a monoidal model category, observe that the pushout product of generating cofibrations can be rewritten as follows:
\<(i⊠‖δ_m‖)◻(j⊠‖δ_n‖)=(i◻j)⊠(‖δ_m‖◻‖δ_n‖).\>
The pushout product $‖δ_m‖◻‖δ_n‖$ is a cofibration in $@PreSmSet$ by ^!{model structures on presheaves are cartesian}
and the pushout product $i◻j$ is a cofibration in the model category~$@V$ because the model structure on~$@V$ is monoidal.
Thus, the pushout product of cofibrations in $@PreSm_@V$ is a cofibration, and likewise for $@Sm_@V$.
On $@PreSm_@V$, the functor $"Sing_@V$ preserves pushouts, monoidal products, and tensorings.
The functor $"hocolim:@V_@Δ→@V$ preserves homotopy pushout squares, and also preserves and reflects weak equivalences.
The cocartesian square for the pushout product of a cofibration and acyclic cofibration in $@PreSm_@V$ is a homotopy pushout square.
Therefore, its image under $"Sing_@V$ followed by $"hocolim$ is a homotopy pushout square.
Therefore, the pushout product is a weak equivalence by the 2-out-of-3 property.
Thus, the pushout product of a cofibration and acyclic cofibration in $@PreSm_@V$ is a weak equivalence.
By ^!{infinity-sheafification is a weak equivalence}, the same holds for $@Sm_@V$.
\ppar
Assuming $@V$ is tractable, h-monoidal, and flat,
the model category $@PreSm_@V$ is tractable because $i⊠‖δ_n‖$ has a cofibrant domain
since $i$ has cofibrant domain, and likewise for $@Sm_@V$.
The nonacyclic part of h-monoidality holds because cofibrations in $@PreSm_@V$ are objectwise h-cofibrations,
the monoidal product of an object and an objectwise h-cofibration is an objectwise h-cofibration by h-monoidality of~$@V$,
and objectwise h-cofibrations are ^{h-cofibrations}.
Flatness in $@PreSm_@V$ follows from the same argument as the acyclic part of the pushout product axiom,
using the fact that the cocartesian square for the pushout product of a cofibration
and a weak equivalence
is a homotopy pushout product square
by the nonacyclic part of h-monoidality.
Flatness in $@Sm_@V$ then follows from ^!{infinity-sheafification is a weak equivalence}.
The acyclic part of h-monoidality holds by Pavlov–Scholbach [\PS, Theorem~3.2.6, Corollary~3.2.8].
(Pretty smallness in the cited results is only used to show that weak equivalences are closed under filtered colimits, which indeed holds in our case.)
\ppar
For symmetric h-monoidality, the argument is the same, using
Pavlov–Scholbach [\PS, Theorem~3.2.7] and
the fact that $"Sing_@V$ preserves colimits in $@PreSm_@V$.
For $@Sm_@V$ we need to further observe that the ^{associated sheaf} functor preserves objectwise h-cofibrations and weak equivalences
by ^!{infinity-sheafification is a weak equivalence}.

\proclaim Example.
(Pavlov–Scholbach [\PS, §7].)
The following model categories satisfy the properties that occur in the statement of ^!{existence of model structures on enriched presheaves}
and ^!{operads in enriched presheaves}.
\typesetli
\li Simplicial sets with simplicial weak equivalences: all properties.
\li Chain complexes (unbounded or nonnegatively graded): all properties except symmetric h-monoidality.
\li Chain complexes in characteristic~0: all properties, and every quasi-isomorphism is symmetric flat.
\li Simplicial modules: all properties.  In characteristic~0 every weak equivalence is symmetric flat.
\li Symmetric simplicial spectra: all properties, weak equivalences are symmetric flat.

\proclaim Theorem.
^^={operads in enriched presheaves}
Suppose $@V$ is a left proper combinatorial model category
that is a tractable (meaning it admits a set of generating cofibrations with cofibrant domains)
symmetric monoidal model category whose weak equivalences are closed under filtered colimits.
In the case of symmetric operads, we assume $@V$ to be symmetric h-monoidal and in the case of nonsymmetric operads, we assume $@V$ to be h-monoidal.
All operads are colored.
All statements below are formulated for $@PreSm_@V$, and an analogous version for $@Sm_@V$ also holds.
\typesetli
\li (Compare ^!{admissibility of operads in presheaves}.)
The category of algebras over any operad~$O$ admits a model structure transferred along the forgetful functor that extracts underlying objects.
\li If $f:O→O'$ is a weak equivalence of operads, then it induces a Quillen equivalence of model categories of algebras over $O$ and $O'$
if and only if $f$ is a (symmetric) flat map.
(In the nonsymmetric case, flat maps coincide with weak equivalences.)
\li (Compare ^!{quasicategorical algebras in presheaves}.)
For every operad~$O$ in $@PreSm_@V$,
the canonical comparison functor
\<@Alg_O(@PreSm_@V)^c[W_O^{-1}]→@Alg_O(@PreSm_@V[W^{-1}])\>
is an equivalence of quasicategories.
\li (Compare ^!{transport of algebras in presheaves}.)
There are Quillen equivalences
\<L⊣R:@Oper_{@V_@Δ} ⇄ @Oper_{@PreSm_@V}, \qquad L'⊣R':@Oper_{@V_@Δ} ⇄ @Oper_@V\>
of model categories of operads in $@V_@Δ$, $@PreSm_@V$, and $@V$.
\li For any cofibrant operad~$O∈@Oper_{@V_@Δ}$, there are Quillen equivalences
\<L_O⊣R_O:@Alg_O(@V_@Δ) ⇄ @Alg_{LO}(@PreSm_@V), \qquad @Alg_O(@V_@Δ) ⇄ @Alg_{L'O}(@V).\>
\li For any fibrant operad $P∈@Oper_{@PreSm_@V}$ (respectively $P'∈@Oper_@V$), there are Quillen equivalences
\<L_P⊣R_P:@Alg_{RP}(@V_@Δ) ⇄ @Alg_P(@PreSm_@V), \qquad @Alg_{R'P'}(@V_@Δ) ⇄ @Alg_{P'}(@PreSm_@V).\>

\proof Proof.
Combine ^!{existence of model structures on enriched presheaves} with Pavlov–Scholbach [\PSa, Theorems 5.11, 7.5, 7.11],
Haugseng [\Haug, Theorem 4.10].
For the last three parts, combine ^!{Quillen equivalence for enriched presheaves} with Pavlov–Scholbach [\PSa, Theorem~8.10].

The following result is implicit in Morel–Voevodsky [\AHTS, Proposition 3.3.3 and Corollary 2.3.5]
and is proved explicitly in Berwick-Evans–Boavida–Pavlov [\BBP, Proposition~1.3].
We give a short self-contained proof.

\proclaim Proposition.
^^={extended simplices are initial}
The functor $`Δ:@Δ→@Cart$ (^!{extended smooth simplex})
is an initial functor and a homotopy initial functor.

\proof Proof.
To show that $`Δ$ is a homotopy initial functor (and hence an initial functor), we verify that for every $V∈@Cart$,
the comma category $`Δ/V$ has a weakly contractible nerve.
Objects of $`Δ/V$ are pairs $([m],`Δ^m→V)$
and morphisms $([m],`Δ^m→V)→([n],`Δ^n→V)$ are maps of simplices $f:[m]→[n]$
that make the triangle with vertices $`Δ^m$, $`Δ^n$, and $V$ commute.
By construction, $`Δ/V$ is the category of simplices of the ^{smooth singular simplicial set} of~$V$.
Hence, the nerve of $`Δ/V$ is isomorphic to the subdivision of the ^{smooth singular simplicial set} of~$V$.
Therefore, the nerve of $`Δ/V$ is weakly equivalent to the ^{smooth singular simplicial set} of~$V$, which is contractible.

\proclaim Theorem.
^^={Quillen equivalence for enriched presheaves}
The right adjoint functors \<"Sing_@V:@PreSm_@V→@V_@Δ, \qquad "Sing_@V:@Sm_@V→@V_@Δ\>
are right Quillen equivalences,
in fact, ^{weak monoidal Quillen equivalences} (^!{weak monoidal Quillen equivalences}).
Here $@PreSm_@V$ and $@Sm_@V$ are equipped with the model structure of ^!{existence of model structures on enriched presheaves}
and $@V_@Δ$ is equipped with the model structure of ^!{model structure on simplicial objects}.

\proof Proof.
We prove the claim for $@PreSm_@V$ first.
We denote the left adjoint of~$"Sing_@V$ by~$‖{-}‖$ (omitting $@V$).
The functor $‖{-}‖$ sends a generating (acyclic) cofibration $i⊗Δ^n$ in $@V_@Δ$ to an (acyclic) cofibration $i⊠‖Δ^n‖$ in $@PreSm_@V$.
Furthermore, the left derived functor of $‖{-}‖$ sends each morphism $X⊗Δ^n→X⊗Δ^0$ to a weak equivalence $X'⊗‖Δ^n‖→X'⊗‖Δ^0‖$,
where $X'→X$ is a cofibrant replacement of~$X$.
Therefore, the functor $‖{-}‖$ is a left Quillen functor by the universal property of left Bousfield localizations.
\ppar
For $@PreSm_@V$, the functor $"Sing_@V$ preserves colimits.
Thus, the derived unit natural transformation $X→"Sing_@V‖X‖$ is cocontinuous in~$X∈@V_@Δ$.
Since weak equivalences in $@V_@Δ$ are closed under filtered colimits,
we can present $X$ as a transfinite composition of cobase changes of morphisms $i⊠δ_n:A→B$, where $δ_n$ is a ^{boundary inclusion} (^!{boundary inclusion})
and $i:P→Q$ is a generating cofibration of~$@V$
and reduce the problem to the following elementary step:
if $X→Y$ is a cobase change of $i⊠δ_n$
such that $X→"SmSing‖X‖$ is a weak equivalence,
then so is $Y→"SmSing‖Y‖$.
Indeed, we have a natural transformation
$$\cd{
A&\mapright{}&X\cr
\lmapdown{i⊠δ_n}&&\mapdown{}\cr
B&\mapright{}&Y\cr
}
\quad⟹\qquad\qquad\quad
\sqcd{
"Sing_@V‖A‖&\mapright{}&"Sing_@V‖X‖\cr
\lmapdown{"Sing_@V ‖i⊠δ_n‖}\hskip.4in&&\hskip.3in\mapdown{}\cr
"Sing_@V‖B‖&\mapright{}&"Sing_@V‖Y‖\cr
}
$$
of corresponding pushout squares.
The component \<X→"Sing_@V‖X‖\> is a weak equivalence by assumption.
The component
\<B→"Sing_@V‖B‖\>
is isomorphic to the map
\<Q⊗Δ^n=Q⊗"SmSing‖Δ^n‖,\>
which is itself a weak equivalence because the map $Δ^n→"SmSing‖Δ^n‖$ has contractible source and target.
By inductive assumption (prove the claim by induction on~$n$),
the component \<A→"Sing_@V‖A‖\> is a weak equivalence.
The maps $i⊠δ_n$ and $"Sing_@V‖i⊠δ_n‖≅i⊠‖δ_n‖$ are cofibrations in $@V_@Δ$, hence h-cofibrations, hence both squares are homotopy pushout squares in $@V_@Δ$
and the component \<Y→"Sing_@V‖Y‖\> is a weak equivalence.
\ppar
For $@Sm_@V$, we combine the previous argument for $@PreSm_@V$ with ^!{infinity-sheafification is a weak equivalence}.
\typesetpar
Finally, to show that the established Quillen equivalences
are weak monoidal Quillen equivalences in the sense of Schwede–Shipley [\EMMC, Definition~3.6],
we observe that the left adjoint functor $‖{-}‖$ (respectively $|{-}|$)
preserves small colimits and commutes with tensoring over~$@V$.
This allows us to prove that the comonoidal maps $L(A⊗B)→LA⊗LB$
are weak equivalences for all cofibrant objects $A,B∈@PreSm_@V$
by induction on~$A$.
If $A=∅$, then the comonoidal map is identity on~$∅$.
Suppose the claim is true for~$A$
and the map $A→A'$ is given by the cobase change of a generating cofibration $i⊠‖δ_n‖$.
The natural transformation of left Quillen functors $L(-⊗B)→L(-)⊗LB$
induces a natural transformation of the resulting cobase change squares.
Since $@V_@Δ$ and $@PreSm_@V$ are left proper, cofibrations are h-cofibrations and the two cobase changes squares are homotopy cobase change squares.
This reduces the problem to showing that the three components of the natural transformation between squares are weak equivalences.
This is true for~$A$ by assumption, holds for the domain of $i⊠‖δ_n‖$ by induction,
and holds for the codomain of $i⊠‖δ_n‖$ by the following argument.
After performing a symmetric reduction for~$B$,
we reduce the problem to the case $A=P⊗‖Δ^m‖$, $B=Q⊗‖Δ^n‖$.
The comonoidal map is $P⊗Q⊗‖Δ^m⨯Δ^n‖→P⊗Q⊗‖Δ^m‖⨯‖Δ^n‖$,
which is a weak equivalence because $‖Δ^m⨯Δ^n‖→‖Δ^m‖⨯‖Δ^n‖$ is a map between weakly contractible objects in $@PreSmSet$.
Thus, the cube lemma (Hovey [\MC, Lemma~5.2.6]) implies that $L(A'⊗B)→LA'⊗LB$ is also a weak equivalence.
\typesetpar
To show weak monoidality in the case of $@Sm_@V$ we combine the previous argument with ^!{infinity-sheafification is a weak equivalence}.

\section The smooth Oka principle for enriched presheaves

The following result improves on the usual way of computing derived internal homs in cartesian model categories by eliminating the fibrant replacement functor.

\proclaim Proposition.
^^={smooth Oka principle for simplicial presheaves}
(The smooth Oka principle for simplicial smooth sets.  Berwick-Evans–Boavida–Pavlov [\BBP, Theorem~1.1].)
If $X$ is a smooth manifold, the functor
\<"Hom(X,-):@PreSm_@sSet→@PreSm_@sSet\>
preserves weak equivalences (^!{weak equivalence of smooth simplicial sets})
and computes the derived internal hom in the model structure of ^!{existence of model structures on enriched presheaves},
yielding a natural simplicial weak equivalence
\<"sSmSing "Hom(X,F) ≃ \rdf"Hom("sSmSing X,"sSmSing F).\>
Here the functor $"=sSmSing$ takes the diagonal of the bisimplicial set $"Sing_@sSet(-)$ (^!{smooth simplicial set}).

\proclaim Remark.
Berwick-Evans–Boavida–Pavlov [\BBP, Theorem~1.1] use simplicial presheaves on the site of smooth manifolds $@Man$,
whereas we used simplicial presheaves on the cartesian site $@Cart$ to formulate ^!{smooth Oka principle for simplicial presheaves}.
However, the version for $@Cart$ is equivalent to the version for $@Man$ in loc.~cit.,
since the product of a manifold with a cartesian space is cofibrant by ^!{manifolds are cofibrant},
so for a Čech-local fibrant simplicial presheaf on~$@Man$,
the internal hom over $@Cart$ is weakly equivalent to the restriction of the internal hom over~$@Man$.

\proclaim Definition.
A ^={model variet[y|ies]} is a combinatorial model category~$@C$
that admits a set~$G$ of objects
such that for every $X∈G$ the functor $"Map(X,-):@C→@sSet$
preserves homotopy sifted homotopy colimits
and every object in~$@C$ is a homotopy sifted homotopy colimit of objects from~$G$.
Here $"=Map(-,-)$ denotes the mapping simplicial set given by the Dwyer–Kan hammock localization of~$@C$.

\proclaim Remark.
^^={left proper model variety}
By Rezk [\Proper, Theorem~B], every ^{model variety} is Quillen equivalent to a left proper simplicial ^{model variety}.
The following equivalent definitions of a ^{model variety} are found in the literature and can be shown to be equivalent to ^!{model variety} using Pavlov [\Comb, Theorem~1.1].
\typesetli
\li
A combinatorial model category
whose underlying quasicategory is equivalent to a ^={projectively generated ∞-categor[y|ies]}
in the sense of Lurie [\HTT, Definition 5.5.8.23].
\li
A combinatorial model category
whose underlying fibrant simplicial category
(e.g., the fibrant replacement of the hammock localization)
is a ^={homotopy variet[y|ies]} in the sense of Rosický [\HVar, Definition~4.10].
\li
A combinatorial model category connected by a chain of Quillen equivalences
to the model category of algebras over a simplicial algebraic theory
(Rosický [\HVar, Theorem~4.15, Corollary~4.18] and [\HVarErr]).

\proclaim Examples.
The following model categories are examples of ^{model varieties}:
\typesetli
\li Simplicial sets, simplicial monoids, simplicial groups, simplicial rings, simplicial objects in any variety of algebras.
\li Many models for connective spectra, e.g., $Γ$-spaces or connective simplicial symmetric spectra.
\li Nonnegatively graded chain complexes with quasi-isomorphisms.
\li $\E_n$-spaces ($0≤n≤∞$) and group-like $\E_n$-spaces ($1≤n≤∞$) in simplicial sets.
\li Many models for connective $\E_n$-ring spectra ($0≤n≤∞$) in simplicial sets.

\proclaim Definition.
Suppose $@V$ is a left proper simplicial ^{model variety} (^!{model variety}, see also ^!{left proper model variety}).
The functor \<\csp:@PreSm_@V→@V, \qquad F↦"=hocolim_{n∈@Δ^@op}F(`Δ^n)\>
is known as the {\it path ∞-groupoid\/} functor,
or, abusing the language, simply as the ^={shape[|s]} functor.
(The shape modality of~$F$ is the locally constant sheaf on the path ∞-groupoid of~$F$.)

\proclaim Theorem.
^^={smooth Oka principle for varieties}
(The smooth Oka principle for model varieties.)
Suppose $@V$ is a left proper simplicial ^{model variety} (^!{model variety}, ^!{left proper model variety}).
\typesetli
\li If $X$ is a smooth manifold, the endofunctor
\<"Hom(X,-):@PreSm_@V→@PreSm_@V, \qquad F↦(M↦F(M⨯X))\>
preserves weak equivalences and therefore computes the derived internal hom.
In particular, $X↦"Hom(X,F)$ is an ∞-sheaf of the form $@Man^@op→@PreSm_@V$.
\li Given a Čech-local object $F∈@PreSm_@V$,
the functor \<\csp("Hom(-,F)):@Man^@op→@V\> is connected by a zigzag of natural weak equivalences to the functor
\<"Hom("SmSing(-),\csp F):@Man^@op→@V,\> where the latter $"Hom$ denotes the powering of $@V$ over simplicial sets.

\proof Proof.
Since $@V$ is a ^{model variety}, we can find a generating set~$G$ of objects in~$@V$ as in ^!{model variety}.
In particular, for any $X∈G$ the functors $"Map(X,-):@PreSm_@V→@sSet$
jointly reflect weak equivalences: if $"Map(X,f)$ is a weak equivalence for all $X∈G$, then $f$ is a weak equivalence.
Furthermore, they preserve all small homotopy limits and homotopy sifted homotopy colimits,
in particular, they preserve the homotopy limits used for Čech descent objects and the homotopy colimit used in the definition of the functor~$\csp$.
Together, these properties allow us to reduce the case of arbitrary ^{model variety}~$@V$ to the case of $@sSet$, which holds by ^!{smooth Oka principle for simplicial presheaves}.

The following result answers a question posed to the author by Kiran Luecke.

\proclaim Proposition.
^^={shapes of presheaves}
(See ^!{infinity-sheafification is a weak equivalence}.)
If $@V$ is a left proper ^{model variety} (^!{model variety}, ^!{left proper model variety}),
$F∈@PreSm_@V$, and $G∈@PreSm_@V$ is the associated Čech-local object (i.e., the associated ∞-sheaf of~$F$),
with the localization map $F→G$,
then the induced map on ^{shapes} (^!{shape}) $\csp F→\csp G$ is a weak equivalence.

\section Applications: classifying spaces
^^={end of part 3}

We revisit the classical theorems on classifications of differential geometric objects
such as closed differential forms (classified by real cohomology),
bundle $(d-1)$-gerbes with connection (classified by integral cohomology),
principal $G$-bundles with connection (classified by the classifying space of~$G$).
In addition to recovering the classical versions of these results for smooth manifolds (^!{manifolds are cofibrant}),
we also establish them in much larger generality for arbitrary cofibrant ^{smooth sets} or ^{simplicial smooth sets}.

\proclaim Example.
Consider the internal hom object \<K=\Hom(M,Ω^n_\closed)\> in $@SmSet$, where $n≥0$ and $M$ is a smooth manifold, or, more generally, any cofibrant ^{smooth set}.
This internal hom computes the derived internal hom because the source $M$ is cofibrant and the target $Ω^n_\closed$ is fibrant.
Thus, the ^{shape} of $K$ can be computed
as the derived mapping simplicial set from the ^{shape} $\csp M$ of~$M$ to the ^{shape} of $Ω^n_\closed$.
The latter is simply $\EM(`R,n)$, the $n$th Eilenberg–MacLane space of the reals.
Thus, the smooth set $K$ can be seen as a smooth refinement
of the simplicial set \<\Hom("SmSing M,\EM(`R,n)).\>
In particular, connected components of $K$ are in bijection with $\HH^n(M,`R)$, the $n$th de Rham (or real singular) cohomology of~$M$.
This is well known when $M$ is a manifold, but appears to be new when $M$ is a cofibrant ^{smooth set}.
In concrete terms, the chain complex
\<Ω^n_\closed(M)←Ω^n_\closed(M⨯`Δ^1)←Ω^n_\closed(M⨯`Δ^2)←⋯,\>
where the differential in degree~$m$ is given by alternating sums of face maps of~$`Δ^m$,
is quasi-isomorphic to the chain complex
\<Ω^n_\closed(M)←Ω^{n-1}(M)←Ω^{n-2}(M)←⋯,\>
where the quasi-isomorphism can be implemented by fiberwise integration over the maps $M⨯`Δ^m→M$.

\proclaim Example.
Consider $\Hom(M,D_n=(Ω^n←⋯←Ω^0←`Z))$ the internal hom object in $@Sm_{@Ch_{≥0}}$, where $n≥0$ and $M$ is a cofibrant ^{simplicial smooth set}.
(Here we convert simplicial sets into chain complexes using the normalized chains functor.)
The target~$D_n$ is also known as the {\it Deligne complex}.
This internal hom computes the derived internal hom because the source $M$ is cofibrant and the target is a fibrant object in $@Sm_{@Ch_{≥0}}$.
Thus, the ^{shape} of $\Hom(M,D_n)$ can be computed as the derived mapping chain complex from the ^{shape} $\csp M$ of~$M$ to the ^{shape} of $D_n$.
The latter is simply $\EM(`Z,n+1)$, the $(n+1)$st Eilenberg–MacLane space of the integers.
In particular, this proves that concordance classes of bundle $(n-1)$-gerbes with connections over~$M$
are classified by the group $\HH^{n+1}(\csp M,`Z)$.
This is well known when $M$ is a manifold, but appears to be new when $M$ is a cofibrant ^{simplicial smooth set}.

\proclaim Example.
Consider the internal hom $\Hom(M,\B G)$ in $@Sm_@sSet$, where $G$ is a Lie group and $M$ is a cofibrant ^{simplicial smooth set}.
The target~$\B G$ is the delooping of the representable presheaf of~$G$.
This internal hom computes the derived internal hom because the source $M$ is cofibrant and the target is a fibrant object in $@Sm_@sSet$.
Thus, the ^{shape} of $\Hom(M,\B G)$ can be computed as the derived mapping chain complex from the ^{shape} $\csp M$ of~$M$ to the ^{shape} of $\B G$.
The latter ^{shape} is simply $\tB G$, the classifying space of~$G$ as a topological group,
i.e., the delooping of the singular complex of~$G$.
In particular, this proves that concordance classes of principal $G$-bundles over~$M$
are classified by the set $[\csp M,\tB G]$.
This is well known when $M$ is a manifold,
but appears to be new when $M$ is a cofibrant ^{simplicial smooth set}.

\unsection References

\refs


\inewif\ifmybib
\ifx\documentclass\undefined 
\mybibtrue 
\else 
\mybibfalse
\fi 

\bib\SHDM[1947]
\ifmybib
Samuel Eilenberg. 
Singular homology in differentiable manifolds.
Annals of Mathematics 48:3 (1947), 670–681.
\doi:10.2307/1969134.
\else
S.~Eilenberg. 
Singular homology in differentiable manifolds.
Ann.~Math.~48:3 (1947), 670–681.
\fi

\bib\TDR[1952]
\ifmybib
André Weil.
Sur les théorèmes de de Rham.
Commentarii Mathematici Helvetici 26 (1952), 119–145.
\doi:10.1007/bf02564296, \eudml:139040.
\else
A.~Weil.
Sur les théorèmes de de Rham.
Comm.~Math.~Helv.~26 (1952), 119–145.
\fi

\bib\Chen[1973]
Kuo-Tsai Chen.
Iterated integrals of differential forms and loop space homology.
Annals of Mathematics 97:2 (1973), 217–246.
\doi:10.2307/1970846.

\bib\PenonI[1973]
Jacques Penon.
Quasi-topos.
Comptes Rendus \(Hebdomadaires des Séances de l'Académie des Sciences\) 276 (1973), 237–240.
\(\https://gallica.bnf.fr/ark:/12148/bpt6k6217213f/f251.image.\)

\bib\PenonII[1977]
\ifmybib
Jacques Penon.
Sur le quasi-topos.
Cahiers de topologie et géométrie dif\-férentielle catégoriques 18:2 (1977), 181–218.
\(\numdam CTGDC_1977__18_2_181_0.\)
\else
J.~Penon.
Sur le quasi-topos.
Cahiers top.~géom.~diff.~cat.~18:2 (1977), 181–218.
\fi

\bib\Dubuc[1979]
\ifmybib
Eduardo~J.~Dubuc.
Concrete quasitopoi.
Lecture Notes in Mathematics 753 (1979), 239–254.
\doi:10.1007/BFb0061821.
\else
E.~J.~Dubuc.
Concrete quasitopoi.
Lecture Notes in Math.~753 (1979), 239–254.
\fi

\bib\Souriau[1980]
\ifmybib
Jean-Marie Souriau.
Groupes différentiels.
\(In: Differential Geometrical Methods in Mathematical Physics.\)
Lecture Notes in Mathematics 836 (1980), 91–128.
\(\doi:10.1007/bfb0089728.\)
\else
J.~M.~Souriau.
Groupes différentiels.
Lecture Notes in Math.~836 (1980), 91–128.
\fi

\bib\Conv[1986]
Anders Kock.
Convenient vector spaces embed into the Cahiers topos.
Cahiers de topologie et géométrie différentielle catégoriques 27:1 (1986), 3–17.
\(\numdam CTGDC_1986__27_1_3_0.\)

\bib\AccCat[1989]
Michael Makkai, Robert Paré.
Accessible categories: the foundations of categorical model theory.
Contemporary Mathematics 104 (1989).
\(\doi:10.1090/conm/104.\)

\bib\Losik[1992]
\ifmybib
M.~V.~Losik.
Fréchet manifolds as diffeologic spaces.
Russian Mathematics 36:5 (1992), 36–42. 
\(\https://dmitripavlov.org/scans/losik-frechet-manifolds-as-diffeologic-spaces.pdf.\)
\else
M.~Losik.
Fréchet manifolds as diffeologic spaces.
Russ.~Math.~36:5 (1992), 36–42. 
\fi

\bib\QCMSS[1993]
Sjoerd~E.~Crans.
Quillen closed model structures for sheaves.
Journal of Pure and Applied Algebra 101:1 (1995), 35–57.
\(\doi:10.1016/0022-4049(94)00033-f.\)

\bib\MC[1999]
Mark Hovey.
Model categories.
\ifmybib
Mathematical Surveys and Monographs 63 (1999).
\doi:10.1090/surv/063.
\else
Math.~Surveys Monogr.~63 (1999).
\fi

\bib\AHTS[1999]
Fabien Morel, Vladimir Voevodsky.
${\bf A}^1$-homotopy theory of schemes.
Publications mathématiques de l'I.H.É.S.\ 90:1 (1999), 45--143.
\doi:10.1007/bf02698831.

\bib\Hovey[1999]
Mark Hovey.
Algebraic Topology Problem List.
Model categories.
March \(6,\) 1999.
\bgroup\ifmybib\font\small=cmr7 \small\else\tiny\fi 
\https://web.archive.org/web/20000830081819/http://claude.math.wesleyan.edu/~mhovey/problems/model.html.
\egroup

\bib\SHMC[2000]
Tibor Beke.
Sheafifiable homotopy model categories.
Mathematical Proceedings of the Cambridge Philosophical Society 129:3 (2000), 447--475.
\(\arXiv:math/0102087v1, \doi:10.1017/s0305004100004722.\)

\bib\Proper[2000]
Charles Rezk.
Every homotopy theory of simplicial algebras admits a proper model.
Topology and its Applications 119:1 (2002), 65–94.
\(\arXiv:math/0003065v1, \doi:10.1016/s0166-8641(01)00057-8.\)

\bib\UHT[2000]
\ifmybib
Daniel Dugger.
Universal homotopy theories.
Advances in Mathematics 164:1 (2001), 144–176.
\(\arXiv:math/0007070v1, \doi:10.1006/aima.2001.2014.\)
\else
D.~Dugger.
Universal homotopy theories.
Adv.~Math.~164:1 (2001), 144–176.
\fi

\bib\LPMS[2001]
Benjamin~A.~Blander.
Local projective model structures on simplicial presheaves.
K-Theory 24:3 (2001), 283–301.
\doi:10.1023/a:1013302313123.

\bib\SoE[2002]
\ifmybib
Peter~T.~Johnstone.
Sketches of an Elephant.  A Topos Theory Compendium.  II.
Oxford Logic Guides 44 (2002).
\else
P.~Johnstone.
Sketches of an Elephant II.
Oxford Logic Guides 44 (2002).
\fi

\bib\THT[2002]
Denis-Charles Cisinski.
Théories homotopiques dans les topos.
Journal of Pure and Applied Algebra 174:1 (2002), 43–82.
\(\doi:10.1016/s0022-4049(01)00176-1.\)

\bib\EMMC[2002]
Stefan Schwede, Brooke Shipley.
Equivalences of monoidal model categories.
Algebraic \& Geometric Topology 3 (2003), 287–334.
\(\arXiv:math/0209342v2, \doi:10.2140/agt.2003.3.287.\)

\bib\MCL[2003]
\ifmybib
Philip~S.~Hirschhorn.
Model categories and their localizations.
Mathematical Surveys and Monographs 99 (2003).
\doi:10.1090/surv/099.
\else
P.~S.~Hirschhorn.
Model categories and their localizations.
Math.~Surv.~Monogr.~99 (2003).
\fi

\bib\FLA[2003]
Denis-Charles Cisinski.
Faisceaux localement asphériques.
January 26, 2003.
\https://cisinski.app.uni-regensburg.de/mtest2.pdf.

\bib\HVar[2005]
Jiří Rosický.
On homotopy varieties.
Advances in Mathematics 214 (2007), 525–550.
\arXiv:math/0509655v2, \doi:10.1016/j.aim.2007.02.011.

\bib\LR[2007]
Clark Barwick.
On left and right model categories and left and right Bousfield localizations.
Homology, Homotopy and Applications 12:2 (2010), 245--320.
\(\arXiv:0708.2067v2, \doi:10.4310/hha.2010.v12.n2.a9.\)

\bib\BH[2008]
John~C.~Baez, Alexander~E.~Hoffnung.
Convenient categories of smooth spaces.
Transactions of the American Mathematical Society 363:11 (2011), 5789–5825.
\(\arXiv:0807.1704v4, \doi:10.1090/s0002-9947-2011-05107-x.\)

\bib\Stacey[2008]
\ifmybib
Andrew Stacey.
Comparative smootheology.
Theory and Applications of Categories 25:4 (2011), 64–117.
\(\arXiv:0802.2225v2, \http://tac.mta.ca/tac/volumes/25/4/25-04abs.html.\)
\else
A.~Stacey.
Comparative smootheology.
Theo.~Appl.~Categ.~25:4 (2011), 64–117.
\fi

\bib\Five[2009]
\ifmybib
Hisham Sati, Urs Schreiber, Jim Stasheff.
Twisted differential String and Fivebrane structures.
Communications in Mathematical Physics 315:1 (2012), 169–213.
\arXiv:0910.4001v2, \doi:10.1007/s00220-012-1510-3.
\else
H.~Sati, U.~Schreiber, J.~Stasheff.
Twisted differential String and Fivebrane structures.
Comm.~Math.~Phys.~315:1 (2012), 169–213.
\fi

\bib\Cech[2010]
\ifmybib
Domenico Fiorenza, Urs Schreiber, Jim Stasheff.
Čech cocycles for differential characteristic classes: an ∞-Lie theoretic construction.
Advances in Theoretical and Mathematical Physics 16:1 (2012), 149--250.
\(\arXiv:1011.4735v2, \doi:10.4310/atmp.2012.v16.n1.a5.\)
\else
D.~Fiorenza, U.~Schreiber, J.~Stasheff.
Čech cocycles for differential characteristic classes: an ∞-Lie theoretic construction.
Adv.~Theor.~Math.~Phys.~16:1 (2012), 149--250.
\fi

\bib\IZ[2013]
\ifmybib
Patrick Iglesias-Zemmour.
Diffeology.
Mathematical Surveys and Monographs 185 (2013).
\doi:10.1090/surv/185.
\else
P.~Iglesias-Zemmour.
Diffeology.
Math.~Surveys Monogr.~185 (2013).
\fi

\bib\HTAPM[2013]
Michael Batanin, Clemens Berger.
Homotopy theory for algebras over polynomial monads.
Theory and Applications of Categories 32:6 (2017), 148–253.
\(\arXiv:1305.0086v7, \http://www.tac.mta.ca/tac/volumes/32/6/32-06abs.html.\)

\bib\DCCT[2013]
\ifmybib
Urs Schreiber.
Differential cohomology in a cohesive infinity-topos.
\arXiv:1310.7930v1.
\else
U.~Schreiber.
Differential cohomology in a cohesive infinity-topos.
\arXiv:1310.7930v1.
\fi

\bib\CW[2013]
\ifmybib
J.~Daniel Christensen, Enxin Wu.
The homotopy theory of diffeological spaces.
New York Journal of Mathematics 20 (2014), 1269–1303.
\(\arXiv:1311.6394v4, \https://nyjm.albany.edu/j/2014/20-59.html.\)
\else
J.~D.~Christensen, E.~Wu.
The homotopy theory of diffeological spaces.
New York J.~Math.~20 (2014), 1269–1303.
\fi

\bib\HVarErr[2014]
\ifmybib
Jiří Rosický.
Corrigendum to “On homotopy varieties”.
Advances in Mathematics 259 (2014), 841–842.
\doi:10.1016/j.aim.2014.02.032.
\else
J.~Rosický.
Corrigendum to “On homotopy varieties”.
Adv.~Math.~259 (2014), 841–842.
\fi

\bib\PSa[2014]
Dmitri Pavlov, Jakob Scholbach.
Admissibility and rectification of colored symmetric operads.
Journal of Topology 11:3 (2018), 559–601.
\(\arXiv:1410.5675v4, \doi:10.1112/topo.12008.\)

\bib\AFR[2015]
David Ayala, John Francis, Nick Rozenblyum.
A stratified homotopy hypothesis.
Journal of the European Mathematical Society 21:4 (2019), 1071–1178.
\(\arXiv:1502.01713v4, \doi:10.4171/jems/856.\)

\bib\PS[2015]
Dmitri Pavlov, Jakob Scholbach.
Homotopy theory of symmetric powers.
Homology, Homotopy and Applications 20:1 (2018), 359–397.
\(\arXiv:1510.04969v3, \doi:10.4310/hha.2018.v20.n1.a20.\)

\bib\Ki[2016]
Hiroshi Kihara.
Model category of diffeological spaces.
Journal of Homotopy and Related Structures 14:1 (2019), 51–90.
\(\arXiv:1605.06794v3, \doi:10.1007/s40062-018-0209-3.\)

\bib\KiQ[2017]
Hiroshi Kihara.
Quillen equivalences between the model categories of smooth spaces, simplicial sets, and arc-generated spaces.
\arXiv:1702.04070v1.

\bib\HTT[2017]
\ifmybib
Jacob Lurie.
Higher Topos Theory.
April 9, 2017.
\https://www.math.ias.edu/~lurie/papers/HTT.pdf.
\else
J.~Lurie.
Higher Topos Theory.
Ann.~of Math.~Stud.~170 (2009).
\fi

\bib\Haug[2019]
Rune Haugseng.
Algebras for enriched ∞-operads.
\arXiv:1909.10042v1.

\bib\BBP[2019]
Daniel Berwick-Evans, Pedro Boavida de Brito, Dmitri Pavlov.
Classifying spaces of infinity-sheaves.
Algebraic \& Geometric Topology (accepted).
\arXiv:1912.10544v2.

\bib\KBig[2020]
Hiroshi Kihara.
Smooth homotopy of infinite-dimensional $C^∞$-manifolds.
Memoirs of the American Mathematical Society 289:1436 (2023).
\arXiv:2002.03618v1, \doi:10.1090/memo/1436.

\bib\Rlocal[2020]
Severin Bunk.
The $`R$-local homotopy theory of smooth spaces.
Journal of Homotopy and Related Structures 17:4 (2022), 593–650.
\(\arXiv:2007.06039v3, \doi:10.1007/s40062-022-00318-7.\)

\bib\POC[2020]
Hisham Sati, Urs Schreiber.
Proper orbifold cohomology.
\arXiv:2008.01101v2.

\bib\TS[2020]
Tadayuki Haraguchi, Kazuhisa Shimakawa.
A model structure on the category of diffeological spaces, I.
\arXiv:2011.12842v1.

\bib\CCGT[2021]
Christopher Adrian Clough.
A convenient category for geometric topology.
Ph.D.~dissertation, The University of Texas at Austin, 2021.
\doi:2152/114981.

\bib\DC[2021]
Araminta Amabel, Arun Debray, Peter~J.~Haine.
Differential Cohomology: Categories, Characteristic Classes, and Connections.
\arXiv:2109.12250v2.

\bib\Comb[2021]
Dmitri Pavlov.
Combinatorial model categories are equivalent to presentable quasicategories.
Journal of Pure and Applied Algebra (accepted).
\arXiv:2110.04679v2.

\bib\SS[2021]
\ifmybib
Hisham Sati, Urs Schreiber.
Equivariant principal ∞-bundles.
\gbreak
\arXiv:2112.13654v3.
\else
H.~Sati, U.~Schreiber.
Equivariant principal ∞-bundles.
\arXiv:2112.13654v3.
\fi

\bib\HTDS[2023]
\ifmybib
Adrian Clough.
The homotopy theory of differentiable sheaves.
\arXiv:2309.01757v1.
\else
A.~Clough.
The homotopy theory of differentiable sheaves.
\arXiv:2309.01757v1.
\fi


%% file: dpmac.tex
\edef\inewcount{\noexpand\csname newcount\endcsname}
\edef\inewdimen{\noexpand\csname newdimen\endcsname}
\edef\inewskip{\noexpand\csname newskip\endcsname}
\edef\inewmuskip{\noexpand\csname newmuskip\endcsname}
\edef\inewbox{\noexpand\csname newbox\endcsname}
\edef\inewhelp{\noexpand\csname newhelp\endcsname}
\edef\inewtoks{\noexpand\csname newtoks\endcsname}
\edef\inewread{\noexpand\csname newread\endcsname}
\edef\inewwrite{\noexpand\csname newwrite\endcsname}
\edef\inewfam{\noexpand\csname newfam\endcsname}
\edef\inewlanguage{\noexpand\csname newlanguage\endcsname}
\edef\inewinsert{\noexpand\csname newinsert\endcsname}
\edef\inewif{\noexpand\csname newif\endcsname}

\countdef\ch=253
\ch="80 \loop\ifnum\ch<"100 \lccode\ch=0 \uccode\ch=0 \advance\ch1 \repeat 
\ch="80 \loop\ifnum\ch<"C0 \catcode\ch=11 \advance\ch1 \repeat 
\ch="C0 \loop\ifnum\ch<"100 \catcode\ch=15 \advance\ch1 \repeat 
\ch="C2 \loop\ifnum\ch<"F5 \catcode\ch=\active \advance\ch1 \repeat 
\ch="C2 \loop\ifnum\ch<"E2 \uccode\ch=1 \advance\ch1 \repeat 
\let\ch\underfined

\catcode0=12 \def\cdef#1#2{\begingroup\lccode0=`#2 \lowercase{\endgroup \def#1{^^@}}} \catcode0=9

\catcode`@=\active
\def@#1{\cdef\chr#1 \edef#1##1{\noexpand\csname\chr##1\noexpand\endcsname}} 
@^^c2@^^c3@^^c4@^^c5@^^c6@^^c7@^^c8@^^c9@^^ca@^^cb@^^cc@^^cd@^^ce@^^cf@^^d0@^^d1@^^d2@^^d3@^^d4@^^d5@^^d6@^^d7@^^d8@^^d9@^^da@^^db@^^dc@^^dd@^^de@^^df
\def@#1{\cdef\chr#1 \edef#1##1##2{\noexpand\csname\chr##1##2\noexpand\endcsname}} 
@^^e0@^^e1@^^e2@^^e3@^^e4@^^e5@^^e6@^^e7@^^e8@^^e9@^^ea@^^eb@^^ec@^^ed@^^ee@^^ef
\def@#1{\cdef\chr#1 \edef#1##1##2##3{\noexpand\csname\chr##1##2##3\noexpand\endcsname}} 
@^^f0@^^f1@^^f2@^^f3@^^f4
\let\chr\undefined
\let\cdef\undefined

\def\grabfuturelet{\futurelet\next\grabexamine}
\def\grabexamine{\ifx\next\csname\expandafter\grab\fi}
\obeylines \def\grab\csname#1\endcsname#2^^M{\expandafter\def\csname#1\endcsname{#2}\expandafter\grabfuturelet} \expandafter\grabfuturelet%
 ~
¢{\hbox{\rm\rlap/c}}
£{\it\$}
«\leftguillemet
­\-
»\rightguillemet
À{\`A}
Á{\'A}
Â{\^A}
Ã{\~A}
Ä{\"A}
Ç{\c C}
È{\`E}
É{\'E}
Ê{\^E}
Ë{\"E}
Ì{\`I}
Í{\'I}
Î{\^I}
Ï{\"I}
Ð\ETH
Ñ{\~N}
Ò{\`O}
Ó{\'O}
Ô{\^O}
Õ{\~O}
Ö{\"O}
Ù{\`U}
Ú{\'U}
Û{\^U}
Ü{\"U}
Ý{\'Y}
Þ\THORN
à{\`a}
á{\'a}
â{\^a}
ã{\~a}
ä{\"a}
ç{\c c}
è{\`e}
é{\'e}
ê{\^e}
ë{\"e}
ì{\`\i}
í{\'\i}
î{\^\i}
ï{\"\i}
ð\eth
ñ{\~n}
ò{\`o}
ó{\'o}
ô{\^o}
õ{\~o}
ö{\"o}
ù{\`u}
ú{\'u}
û{\^u}
ü{\"u}
ý{\'y}
þ\thorn
ÿ{\"y}
Ā{\=A}
ā{\=a}
Ă{\u A}
ă{\u a}
Ć{\'C}
ć{\'c}
Ĉ{\^C}
ĉ{\^c}
Ċ{\.C}
ċ{\.c}
Č{\v C}
č{\v c}
Ď{\v D}
ď{\v d}
Ē{\=E}
ē{\=e}
Ĕ{\u E}
ĕ{\u e}
Ė{\.E}
ė{\.e}
Ě{\v E}
ě{\v e}
Ĝ{\^G}
ĝ{\^g}
Ğ{\u G}
ğ{\u g}
Ġ{\.G}
ġ{\.g}
Ģ{\c G}
ģ{\c g}
Ĥ{\^H}
ĥ{\^h}
Ĩ{\~I}
ĩ{\~\i}
Ī{\=I}
ī{\=\i}
Ĭ{\u I}
ĭ{\u\i}
İ{\.I}
ĲIJ
ĳij
Ĵ{\^J}
ĵ{\^\j}
Ķ{\c K}
ķ{\c k}
Ĺ{\'L}
ĺ{\'l}
Ļ{\c L}
ļ{\c l}
Ľ{\v L}
ľ{\v l}
Ń{\'N}
ń{\'n}
Ņ{\c N}
ņ{\c n}
Ň{\v N}
ň{\v n}
Ō{\=O}
ō{\=o}
Ŏ{\u O}
ŏ{\u o}
Ő{\"O}
ő{\"o}
Ŕ{\'R}
ŕ{\'r}
Ŗ{\c R}
ŗ{\c r}
Ř{\v R}
ř{\v r}
Ś{\'S}
ś{\'s}
Ŝ{\^S}
ŝ{\^s}
Ş{\c S}
ş{\c s}
Š{\v S}
š{\v s}
Ţ{\c T}
ţ{\c t}
Ť{\v T}
ť{\v t}
Ũ{\~U}
ũ{\~u}
Ū{\=U}
ū{\=u}
Ŭ{\u U}
ŭ{\u u}
Ű{\H U}
ű{\H u}
Ŵ{\^W}
ŵ{\^w}
Ŷ{\^Y}
ŷ{\^y}
Ÿ{\"Y}
Ź{\'Z}
ź{\'z}
Ż{\.Z}
ż{\.z}
Ž{\v Z}
ž{\v z}
’'
‘`
”{''}
“{``}
‐-
–{--}
—{---}
¡{!`}
¿{?`}
−-
′'
ß\ss
æ\ae
Æ\AE
œ\oe
Œ\OE
ø\o
Ø\O
å\aa
Å\AA
ł\l
Ł\L
†\dag
‡\ddag
§\S
¶\P
©\copyright
…\dots
ı\relax\ifmmode\imath\else\i\fi
ȷ\relax\ifmmode\jmath\else\j\fi
¬\neg
△\bigtriangleup
÷\div
±\pm
◯\bigcirc
×\times
‖\|
←\leftarrow
→\rightarrow
∥\Vert
⟩\rangle
⟨\langle
∣|
∮\oint
≐\dot=
↑\uparrow
↓\downarrow
α\alpha
β\beta
γ\gamma
δ\delta
ϵ\epsilon
ζ\zeta
η\eta
θ\theta
ι\iota
κ\kappa
λ\lambda
μ\mu
ν\nu
ξ\xi
οo
π\pi
ρ\rho
σ\sigma
τ\tau
υ\upsilon
ϕ\phi
χ\chi
ψ\psi
ω\omega
ε\varepsilon
ϑ\vartheta
ϖ\varpi
ϱ\varrho
ς\varsigma
φ\varphi
Γ\Gamma
Δ\Delta
Θ\Theta
Λ\Lambda
Ξ\Xi
Π\Pi
Σ\Sigma
Υ\Upsilon
Φ\Phi
Ψ\Psi
Ω\Omega
ℵ\aleph
ℏ\hbar
ℓ\ell
℘\wp
ℜ\Re
ℑ\Im
∂\partial
∞\infty
∅\emptyset
∇\nabla
√\surd
⊤\top
⊥\bot
∠\angle
△\triangle
∀\forall
∃\exists
♭\flat
♮\natural
♯\sharp
♣\clubsuit
♢\diamondsuit
♡\heartsuit
♠\spadesuit
∐\coprod
⋁\bigvee
⋀\bigwedge
⨄\biguplus
⋂\bigcap
⋃\bigcup
∫\int
∏\prod
∑\sum
⨂\bigotimes
⨁\bigoplus
⨀\bigodot
⨆\bigsqcup
◁\triangleleft
▷\triangleright
▽\bigtriangledown
∧\wedge
∨\vee
∩\cap
∪\cup
⊓\sqcap
⊔\sqcup
⊎\uplus
⨿\amalg
⋄\diamond
∙\bullet
≀\wr
⊙\odot
⊘\oslash
⊗\otimes
⊖\ominus
⊕\oplus
∓\mp
∘\circ
○\Orb
∖\setminus
⋅\cdot
∗\ast
⨯\times
⋆\star
∝\propto
⊑\sqsubseteq
⊒\sqsupseteq
∥\parallel
∣\divides
⊣\dashv
⊢\vdash
↗\nearrow
↘\searrow
↖\nwarrow
↙\swarrow
⇔\Leftrightarrow
⇐\Leftarrow
⇒\Rightarrow
≠\neq
≤\leq
≥\geq
≻\succ
≺\prec
≈\approx
≽\succeq
≼\preceq
⊃\supset
⊂\subset
⊇\supseteq
⊆\subseteq
∈\in
∋\ni
≫\gg
≪\ll
↔\leftrightarrow
↦\mapsto
∼\sim
≃\simeq
⟂\perp
≡\equiv
≍\asymp
⌣\smile
⌢\frown
↼\leftharpoonup
↽\leftharpoondown
⇀\rightharpoonup
⇁\rightharpoondown
↪\hookrightarrow
↩\hookleftarrow
⋈\bowtie
⊨\models
⟹\Longrightarrow
⟶\longrightarrow
⟵\longleftarrow
⟸\Longleftarrow
⟼\longmapsto
⟷\longleftrightarrow
⟺\Longleftrightarrow
⋯\cdots
⋮\vdots
⋱\ddots
↕\updownarrow
⇑\Uparrow
⇓\Downarrow
⇕\Updownarrow
⌉\rceil
⌈\lceil
⌋\rfloor
⌊\lfloor
≅\cong
∉\notin
⇌\rightleftharpoons
≐\doteq
∄\not\exists
∌\not\ni
∔\dot+
∕/
∤\not|
∦\not\|
∬\int\!\!\!\int
∭\int\!\!\!\int\!\!\!\int
∮\oint
∸\dot-
≁\not\sim
≄\not\simeq
≆\not\cong
≇\not\cong
≉\not\approx
≔:=
≕=:
≢\not\equiv
≭\not\asump
≮\not<
≯\not<
≰\not\le
≱\not\ge
⊀\not\prec
⊁\not\succ
⊄\not\subset
⊅\not\supset
⊈\not\subseteq
⊉\not\supseteq
⊦\vdash
⊧\models
⊬\not\vdash
⊭\not\models
⊲\triangleleft
⊳\triangleright
⋠\not\preceq
⋡\not\succeq
⋤\not\sqsubseteq
⋥\not\sqsupseteq
⋪\not\triangleleft
⋫\not\triangleright
◻\square

\catcode`\^^M=5 %
\let\grabfuturelet\undefined \let\grabexamine\undefined \let\grab\undefined

\let\xcsname=\csname
\let\xendcsname=\endcsname
\def@#1{\def#1##1{\expandafter\ifx\csname\string#1##1\endcsname\relax\errmessage{Undefined UTF-8 sequence \string#1##1}\else\xcsname\string#1##1\xendcsname\fi}}
@^^c2@^^c3@^^c4@^^c5@^^c6@^^c7@^^c8@^^c9@^^ca@^^cb@^^cc@^^cd@^^ce@^^cf@^^d0@^^d1@^^d2@^^d3@^^d4@^^d5@^^d6@^^d7@^^d8@^^d9@^^da@^^db@^^dc@^^dd@^^de@^^df
\def@#1{\def#1##1##2{\expandafter\ifx\csname\string#1##1##2\endcsname\relax\errmessage{Undefined UTF-8 sequence \string#1##1##2}\else\xcsname\string#1##1##2\xendcsname\fi}}
@^^e0@^^e1@^^e2@^^e3@^^e4@^^e5@^^e6@^^e7@^^e8@^^e9@^^ea@^^eb@^^ec@^^ed@^^ee@^^ef
\def@#1{\def#1##1##2##3{\expandafter\ifx\csname\string#1##1##2##3\endcsname\relax\errmessage{Undefined UTF-8 sequence \string#1##1##2##3}\else\xcsname\string#1##1##2##3\xendcsname\fi}}
@^^f0@^^f1@^^f2@^^f3@^^f4
\let@\undefined
\catcode`@=12

\newif\ifscroll 
\newif\ifsuppressunusedbib 

\def\printerr#1{\immediate\write17{#1}}
\def\warningline#1#2{\printerr{! #2}\printerr{l.#1}\printerr{}}
\def\ewarningline#1#2#3{\printerr{! #2}\printerr{l.#1 #3}\printerr{}}
\def\warning{\warningline{\the\inputlineno}}

\long\def\gobble#1{}
\ifx\gobbleinit\undefined{\long\gdef\gobbleinit#1\par{}}\fi
\def\expand#1{\edef\expandmacro{#1}\expandmacro\let\expandmacro\undefined}
\def\setetok#1#2{\expand{\noexpand#1{#2}}}
\def\expandtoks#1{\expandafter\edef\expandafter\expandmacro\expandafter{\the#1}#1\expandafter{\expandmacro}}
\def\appendexpand#1#2{\setetok#1{\the#1#2}}
\long\def\append#1#2{#1\expandafter{\the#1#2}}
\long\def\appendtoksexpand#1#2{#1\expandafter\expandafter\expandafter{\expandafter\the\expandafter#1\the#2}}
\long\def\appendonceexpand#1#2{#1\expandafter\expandafter\expandafter{\expandafter\the\expandafter#1#2}}

\def\link#1#2{\lhighlight{#2}}
\def\llink#1{\printlink{llink #1}\link{\ohash#1}}
\catcode`\#=11 \def\ohash{#}\catcode`\#=6
\catcode`\&=11 \def\ampersand{&}\catcode`\&=4
\def\anchor#1#2{\printlink{anchor #1 #2}#2}

\def\setpapersize#1#2{} 
\def\dumpbox#1#2#3{\shipout\vbox{\setpapersize{#1}{#2}\unvbox#3}}
\def\mps#1{\epsfbox{#1}}
\def\metadata#1#2{}
\def\src{} 

\newread\epsffilein
\newif\ifepsfbbfound\inewif\ifepsffilecont
\newdimen\epsfxsize\inewdimen\epsfysize
\newdimen\pspoints\pspoints1bp
\let\runmp\errmessage 
\def\epsfbox#1{\openin\epsffilein=#1 \ifeof\epsffilein\runmp{Could not open file #1}\else
	{\def\do##1{\catcode`##1=12}\dospecials\catcode`\ =10\epsffileconttrue
		\epsfbbfoundfalse
		\loop\read\epsffilein to\epsffileline \ifeof\epsffilein\epsffilecontfalse\else\expandafter\epsfaux\epsffileline :. \\\fi\ifepsffilecont\repeat
		\ifepsfbbfound\else\errmessage{No HiResBoundingBox comment found in file #1}\fi}%
	\closein\epsffilein
	\epsfysize\epsfury\pspoints \advance\epsfysize-\epsflly\pspoints
	\epsfxsize\epsfurx\pspoints \advance\epsfxsize-\epsfllx\pspoints
	\setbox0\hbox{\vbox to\epsfury\pspoints{\vfil\hbox to\epsfxsize{\dimen0=\epsfllx\pspoints \kern-\dimen0 \includegraphics{#1}\hfil}}}%
	\dp0=\epsflly\pspoints \dp0=-\dp0
	\box0 \fi}
\catcode`\%=12 \def\epsfbblit{
\let\dummybrace=} 
\def\epsfaux#1:#2\\{\def\testit{#1}\ifx\testit\epsfbblit \epsfgrab #2 . . . \\\epsffilecontfalse\epsfbbfoundtrue\fi}
\def\empty{}
\def\epsfgrab #1 #2 #3 #4 #5\\{\gdef\epsfllx{#1}\ifx\epsfllx\empty\epsfgrab #2 #3 #4 #5 .\\\else\gdef\epsflly{#2}\gdef\epsfurx{#3}\gdef\epsfury{#4}\fi} 

\newif\ifpdf \pdffalse \ifx\pdfoutput\undefined\else\ifx\pdfoutput\relax\else\ifnum\pdfoutput<1 \else\pdftrue\fi\fi\fi
\ifpdf
\def\pdflink#1#2{\leavevmode \lhighlight{\pdfstartlink user { /Subtype /Link /Border [0 0 0] /A << /S #1 >> }#2\pdfendlink}}
\def\llink#1{\pdflink{/GoTo /D (#1)}}
\def\link#1{\pdflink{/URI /URI (#1)}}
\def\anchor#1#2{\pdfdest name {#1} xyz #2}

\def\setpapersize#1#2{\pdfpagewidth#1 \pdfpageheight#2 }
\def\dumpbox#1#2#3{\setpapersize{#1}{#2}\shipout\box#3}
\def\metadata#1#2{\pdfinfo{/Title (#1) /Author (#2)}}
\input supp-pdf.mkii 
\def\mps#1{\convertMPtoPDF{#1}{1}{1}}
\fi

\newtoks\buffertoks
\def\addcode{\immediate\write\mpout}
\def\addunexpandedcode#1{{\toks0={#1}\addcode{\the\toks0}}}
\def\addcodebuffer#1{\edef\tmp{#1}\buffertoks\expandafter\expandafter\expandafter{\expandafter\the\expandafter\buffertoks\tmp}}
\def\addunexpandedcodebuffer#1{\buffertoks\expandafter{\the\buffertoks#1}}
\def\inlinemp{\begininlinemp\grabcode}
\def\grabcode{\catcode`\#=12 \endlinechar=10 
	\afterassignment\dumpcode\outtoks} 
\def\dumpcode{\addcode{\the\outtoks}\endinlinemp\gobble} 
\def\mpsfilename{\filestem.\the\figno}
\def\begininlinemp{\inimp\begingroup\catcode`\^=7 \iftypesetting\mps{\mpsfilename}\let\addcode\gobble\fi \addcode{beginfig(\the\figno);}}
\def\endinlinemp{\addcode{endfig;}\addcode{}\endgroup\global\advance\figno1\relax}
\ifx\endprolog\undefined\let\endprolog\relax\fi

\def\plainfmtname{plain}\ifx\fmtname\plainfmtname\else
\edef\plainoutput{\the\output}
\global\chardef\itfam=4
\def\_{\leavevmode \kern.06em \vbox{\hrule width.3em}} 

\font\tensy=cmsy10

\catcode"18=12
\catcode`@=11
{
\global\let\end\@@end
\global\let\input\@@input
}
\def\eqalign#1{\null\,\vcenter{\openup\jot\m@th
  \ialign{\strut\hfil$\displaystyle{##}$&$\displaystyle{{}##}$\hfil
      \crcr#1\crcr}}\,}
\catcode`@=12
\fi

\font\tenmsa=msam10 \font\sevenmsa=msam7 \font\fivemsa=msam5 \newfam\msafam \textfont\msafam=\tenmsa \scriptfont\msafam=\sevenmsa \scriptscriptfont\msafam=\fivemsa 
\font\teneufm=eufm10 \font\seveneufm=eufm7 \font\fiveeufm=eufm5 \newfam\eufmfam \textfont\eufmfam=\teneufm \scriptfont\eufmfam=\seveneufm \scriptscriptfont\eufmfam=\fiveeufm 
\def\frak{\fam\eufmfam}
\font\teneufb=eufb10 \font\seveneufb=eufb7 \font\fiveeufb=eufb5 \newfam\eufbfam \textfont\eufbfam=\teneufb \scriptfont\eufbfam=\seveneufb \scriptscriptfont\eufbfam=\fiveeufb 

\font\teneurm=eurm10 \font\seveneurm=eurm7 \font\fiveeurm=eurm5 \newfam\eurmfam \textfont\eurmfam=\teneurm \scriptfont\eurmfam=\seveneurm \scriptscriptfont\eurmfam=\fiveeurm 
\def\eurm{\fam\eurmfam}
\font\teneurb=eurb10 \font\seveneurb=eurb7 \font\fiveeurb=eurb5 \newfam\eurbfam \textfont\eurbfam=\teneurb \scriptfont\eurbfam=\seveneurb \scriptscriptfont\eurbfam=\fiveeurb 

\font\teneusm=eusm10 \font\seveneusm=eusm7 \font\fiveeusm=eusm5 \newfam\eusmfam \textfont\eusmfam=\teneusm \scriptfont\eusmfam=\seveneusm \scriptscriptfont\eusmfam=\fiveeusm 

\font\teneusb=eusb10 \font\seveneusb=eusb7 \font\fiveeusb=eusb5 \newfam\eusbfam \textfont\eusbfam=\teneusb \scriptfont\eusbfam=\seveneusb \scriptscriptfont\eusbfam=\fiveeusb 

\font\tenss=cmss10 \font\sevenss=cmss7 \font\fivess=cmss5 \newfam\ssfam \textfont\ssfam\tenss \scriptfont\ssfam\sevenss \scriptscriptfont\ssfam\fivess 
\def\sf{\fam\ssfam}

\font\sevenit=cmti7 \scriptfont\itfam=\sevenit

\let\articletitle\seventeenss
\let\chaptertitle\twelvebf
\let\sectiontitle\tenbf
\let\subsectiontitle\tenbfit
\let\subsubsectiontitle\tenit
\let\contchaptertitle\tenbf 
\let\contsectiontitle\tenrm 
\let\contsubsectiontitle\sevenrm 
\let\contsubsubsectiontitle\fiverm 
\let\parnumfont\tenrm 
\let\parbackreffont\fiverm 
\let\proclaimfont\tenbf

\let\mainfont\tenrm

\def\lhighlight{} 
\def\suppresscomments#1#2#3{} 
\def\prohibitcomments#1#2#3{\errmessage{Draft comments are not allowed in the final version}} 

\ifx\format\undefined\else\format\fi 
\ifx\comment\undefined\def\comment{\prohibitcomments}\fi

\newdimen\hmargin
\newdimen\vmargin
\newdimen\plaintextwidth
\newdimen\plaintextheight
\newdimen\totalht
\newdimen\totalwd
\def\shipbox#1{%
	\totalht\ht#1
	\advance\totalht2\vmargin
	\totalwd\wd#1
	\advance\totalwd2\hmargin
	\hoffset-1in
	\advance\hoffset\hmargin
	\voffset-1in
	\advance\voffset\vmargin
	\dumpbox\totalwd\totalht{#1}}

\newtoks\cont
\def\contents{\begingroup\suppressbackreftrue\let\anchor\gobble\the\cont\endgroup}
\let\printcont\gobble
\def\contlinechapt#1#2{\printcont{#2}\smallskip\everypar{}\noindent{\contchaptertitle\llink{chapter.#1}{#2}}\par} 
\def\contlinesect#1#2{\printcont{#1.  #2}\everypar{}\indent{\contsectiontitle\llap{#1.\enskip}\llink{section.#1}{#2}}\par} 
\def\contlineunsect#1#2{\printcont{#2}\everypar{}\noindent{\contsectiontitle\llink{section.#1}{#2}}\par} 
\def\contlinesubsect#1#2{\printcont{#1.  #2}\everypar{}\indent{\contsubsectiontitle{#1.\enskip}\llink{paragraph.#1}{#2}}\par} 
\def\contlinesubsubsect#1#2{\printcont{#1.  #2}\everypar{}\indent\indent\indent{\contsubsubsectiontitle\llap{#1.\enskip}\llink{paragraph.#1}{#2}}\par} 
\def\addcont#1#2#3{\append\cont{#1}\appendexpand\cont{{#2}}\append\cont{{#3}}}

\newif\ifpresec \presecfalse 
\def\chapter#1\par{%
	\def\chapname{#1}%
	\parn0
	\subparn0
	\presectrue
	\numbfalse
	\addcont\contlinechapt\chapname{#1}%
	\curverb{Chapter~}%
	\assignlabel{chapter}{\chapname}%
	\tchapter#1\par}
\def\tchapter#1\par{
	\everypar{\let\beforesect\beforesection}
	\chapbreak\bigbreak
	\centerline{\plabel\chaptertitle\anchor{chapter.#1}{#1}}
	\nobreak\medskip
	\let\beforesect\relax 
}

\let\chapbreak\relax 

\newcount\secn \secn0
\def\section#1\par{%
	\ifx\sectionid\undefined\advance\secn1 \edef\sectionid{\the\secn}\fi
	\everypar{\numpar}
	\parn0
	\subparn0
	\presecfalse
	\numbfalse
	\addcont\contlinesect\sectionid{#1}%
	\curverb{\S}%
	\assignlabel{section}{\sectionid}%
	\tsection#1\par
}
\def\unsection#1\par{%
	\everypar{\numpar}
	\parn0
	\subparn0
	\presecfalse
	\numbfalse
	\addcont\contlineunsect{#1}{#1}%
	\curverb{\S}%
	\assignlabel{section}{#1}%
	\tsection#1\par
}
\def\tsection#1\par{
	\beforesect\let\beforesect\beforesection
	\typesetsection{#1}%
	\aftersection	
	\let\sectionid\undefined
}
\def\beforesection{\vskip0pt plus.3\vsize \penalty-250 \vskip0pt plus-.3\vsize \bigskip \vskip\parskip}
\def\aftersection{\nobreak\smallskip}
\def\typesetsection#1{\leftline{\sectiontitle
	\ifx\sectionid\undefined
		\plabel\anchor{section.#1}{}%
	\else
		\hbox to \parindent{\hss\plabel\anchor{section.\sectionid}{\sectionid}\enspace\hfill}%
	\fi#1}} 
\let\beforesect\beforesection

\def\subsection#1\par{\bigbreak\numbtrue\curverb{\S}\subsectiontitle\noindent#1\/\mainfont\par\nobreak\medskip\addcont\contlinesubsect{\the\secn.\the\parn}{#1}}
\def\subsubsection#1\par{\bigbreak\numbtrue\curverb{\S}\subsubsectiontitle\noindent#1\/\mainfont\par\nobreak\medskip\addcont\contlinesubsubsect{\the\secn.\the\parn}{#1}}

\def\sskip#1{\ifdim\lastskip<\medskipamount \removelastskip\penalty#1\medskip\fi}
\def\slug{\hbox{\kern1.5pt\vrule width2.5pt height6pt depth1.5pt\kern1.5pt}}
\newif\ifqed \newif\ifneedqed
\def\qed{\unskip\nobreak\ \slug\ifhmode\spacefactor3000 \fi\global\qedtrue}
\def\proclaim{\medbreak\atendpar{\sskip{55}}\numbtrue\gproclaim\proclaimfont} 
\def\proof{\medbreak\atendpar{\sskip{-55}}\needqedtrue\gproclaim\prooffont}
\def\xproclaim#1.{\medbreak\atendpar{\sskip{55}}{\everypar{}\noindent}{\proclaimfont#1.\enspace}\ignorespaces} 
\let\abstract\xproclaim 

\def\ppar{\endgraf{\everypar{}\indent}} 
\newtoks\atendpar 
\newtoks\atendbr 
\newif\ifnumb \numbfalse 
\def\finishpar{\ifhmode\ifneedqed\ifqed\else\qed\fi\qedfalse\needqedfalse\fi\iflist\endlist\fi\the\atendbr\atendbr{}\endgraf\the\atendpar\atendpar{}\numbfalse\fi}
\def\endlist{\iflist\listfalse\endgraf{\parskip\smallskipamount\everypar{}\noindent}\fi} 
\newcount\parn
\newcount\subparn
\newif\ifparbref
\def\nextpar{\ifnumb\advance\parn1 \printlabel{advancing paragraph number to \the\parn}\def\brt{}%
	\ifparbref\edef\cseq{\csname backreference-list.paragraph.\the\secn.\the\parn\endcsname}%
	\expandafter\ifx\cseq\relax\else\edef\brt{\cseq}\printbackref{back references for paragraph.\the\secn.\the\parn: \cseq}\fi\fi
	\assignlabel{paragraph}{\the\secn.\the\parn}\fi}
\def\numpar{\ifnumb\nextpar\expand{\noexpand\typesetparnum{\the\secn.\the\parn}\noexpand\typesetpbr{\brt}}\fi}
\newdimen\pindent \pindent\parindent

\ifx\draftnum\undefined 
\def\gproclaim#1#2.{\curverb{#2~}
	\ifnumb\nextpar\fi
	{\everypar{}\noindent}%
	\plabel#1#2%
	\ifnumb\ \anchor{paragraph.\the\secn.\the\parn}{}\the\secn.\the\parn\fi.\enspace\mainfont 
	\ifnumb\expand{\noexpand\typesetpbr{\brt}}\fi
	\ignorespaces}
\def\typesetparnum#1{\ifnumb{\plabel\parnumfont\anchor{paragraph.#1}{}#1.\enspace}\fi} 
\def\typesetpbr#1{\ifnumb\def\brtext{#1}\ifx\brtext\empty\else\setetok\atendbr{{\parbackreffont Used in \noexpand\stripcomma\brtext.}}\fi\fi}
\else 
\let\numbfalse\relax \numbtrue
\def\gproclaim#1#2.{\curverb{#2~}\medbreak\noindent#1#2.\enspace\mainfont\ignorespaces}
\inewdimen\brwidth \brwidth.6in
\def\typesetparnum#1{\ifnumb\llap{\plabel\anchor{paragraph.#1}{}\parnumfont#1\enspace}\fi} 
\def\typesetpbr#1{\ifnumb\def\brtext{#1}\ifx\brtext\empty\else
	\llap{\smash{\vtop{\everypar{}\raggedright\rightskip0pt plus 0pt \leftskip0pt plus 1fill \hsize\brwidth
	\parnumfont \strut \break 
	\parbackreffont\stripcomma#1}}\enspace}\fi\fi}
\fi

\def\hang{\hangindent\pindent}
\newif\iflist
\def\textindent#1{{\everypar{}\parindent\pindent\indent}\llap{#1\enspace}\listtrue\ignorespaces}

\def\li{\item{$\bullet$}}

\def\item{\endgraf\hang\textindent}
\def\filbreak{\endgraf\vfil\penalty-200\vfilneg}

\def\eject{\endgraf\break}
\def\supereject{\endgraf\penalty-20000}
\def\smallbreak{\endgraf\ifdim\lastskip<\smallskipamount
	\removelastskip\penalty-50\smallskip\fi}
\def\medbreak{\endgraf\ifdim\lastskip<\medskipamount
	\removelastskip\penalty-100\medskip\fi}
\def\bigbreak{\endgraf\ifdim\lastskip<\bigskipamount
	\removelastskip\penalty-200\bigskip\fi}

\let\printbackref\gobble
\def\predefbackref#1{%
	\printbackref{defining back reference list backreference-list.#1}%
	\expandafter\gdef\expandafter\cseq\expandafter{\csname backreference-list.#1\endcsname}
	\expandafter\ifx\cseq\relax\expandafter\gdef\cseq{}\else\printbackref{duplicate omitted}\fi
}
\newcount\backref \backref0
\newif\ifsuppressbackref \suppressbackreffalse
\def\firstletter#1#2\endletter{#1}
\newtoks\backreflist
\newif\ifaddbr \addbrfalse
\def\recordbackref#1{%
	\edef\params{{\ifpresec\expandafter\firstletter\chapname\endletter\else\the\secn\fi.\the\parn\ifnumb\else*\fi}{\the\backref}{\the\inputlineno}}%
	\printbackref{recording back reference \string#1 for future processing with params \params}%
	\edef\tmp{\the\backreflist\noexpand\processbackref\noexpand#1\params}%
	\global\backreflist\expandafter{\tmp}%
	\printbackref{new content of backreflist: \the\backreflist}} 
\def\processbackref#1#2#3#4{%
	\edef\key{\expandafter\gobble\string#1}%
	\edef\cseq{\csname id.\key\endcsname}%
	\expandafter\ifx\cseq\relax \ewarningline{#4}{Undefined reference \string#1.}{\string#1}\else
		\edef\lseq{\csname backreference-list.\cseq\endcsname}%
		\printbackref{processing back reference number #3 \string#1, originating from #2, at line #4; adding to \lseq}%
		\edef\lastnumber{\csname lastnumber.\cseq\endcsname}%
		\edef\newnumber{#2}%
		\addbrtrue
		\printbackref{lastnumber: \lastnumber; newnumber: \newnumber;}%
		\expandafter\ifx\csname lastnumber.\cseq\endcsname\relax 
		\else\ifx\lastnumber\newnumber 
			\printbackref{Suppressing duplicate back reference #2.}%
			\addbrfalse
		\fi\fi
		\expandafter\edef\csname lastnumber.\cseq\endcsname{#2}
		\ifaddbr
			\expandafter\expandafter\expandafter\gdef\expandafter\expandafter\csname backreference-list.\cseq\endcsname\expandafter{\lseq, \llink{backreference.#3}{#2}}
		\fi
	\fi
}
\newtoks\labelinitlist
\def\xxstripcomma, {}
\def\xstripcomma{\if\ntok,\let\xcont\xxstripcomma\else\let\xcont\relax\fi\xcont}
\def\stripcomma{\futurelet\ntok\xstripcomma}

\inewif\ifrecorddups \recorddupstrue
\def\checkduplicates#1#2{\edef\key{\expandafter\gobble\string#1}%
	\iftypesetting\else
	\expandafter\ifx\csname line:\key\endcsname\relax\printlabel{keydefline: relax}\ifrecorddups\expandafter\xdef\csname line:\key\endcsname{\the\inputlineno}\fi\else\edef\keydefline{\csname line:\key\endcsname}\errmessage{#2}\fi\fi}

\let\printlabel\gobble
\newtoks\curverb 
\ifx\draftlabel\undefined
\def\plabel{}
\else
\def\labeltext{} 
\def\plabel{\ifx\labeltext\empty\else\smash{\llap{\parbackreffont\labeltext\quad}}\gdef\labeltext{}\fi} 
\fi
\def\assignlabel#1#2{
	\ifx\lastlabel\undefined\else
	\edef\key{\expandafter\expandafter\expandafter\gobble\expandafter\string\lastlabel}
	\printlabel{label \key: id.\key\space = #1.#2, text.\key\space = #2}%
	\expandafter\xdef\csname id.\key\endcsname{#1.#2}
	\expandafter\xdef\csname text.\key\endcsname{#2}
	\edef\tmp{\the\curverb}%
	\ifx\tmp\empty\else 
	\printlabel{label v\key: id.v\key\space = #1.#2, text.v\key\space = \the\curverb#2}%
	\expandafter\xdef\csname id.v\key\endcsname{#1.#2}
	\expandafter\xdef\csname text.v\key\endcsname{\the\curverb#2}
	\fi
	\predefbackref{#1.#2}%
	\fi\let\lastlabel\undefined}
\newtoks\vlist 
\def\label#1{%
	\iftypesetting\else
	\checkduplicates#1{Label \string#1 was already defined at line \keydefline}%
	\addverunused#1\verifylabel
	\def\lastlabel{#1}%
	\fi
	\numbtrue 
}

\newif\ifyearkey 
\def\y{} 
\newdimen\bibindent 
\newtoks\bibt 
\def\tbib#1{
	\checkduplicates#1{Bibliographic reference \string#1 already defined at line \keydefline}%
	\addverunused#1\verifybib
	\edef\key{\expandafter\gobble\string#1}
	\printlabel{reference \key: id.\key\space = reference.\key, text.\key\space = \key}%
	\expandafter\edef\csname id.\key\endcsname{reference.\key}
	\expandafter\edef\csname text.\key\endcsname{\key}
	\predefbackref{reference.\key}%
	\ifyearkey\else\setbox0=\hbox{[\key]}\ifdim\bibindent<\wd0 \bibindent=\wd0 \fi \fi 
	\appendexpand\bibt{\noexpand\typesetbib\noexpand#1\src}%
	\ifyearkey\let\next\xxftbib\else\let\next\ftbibalpha\fi\next}
\def\ftbibalpha#1\par{\append\bibt{#1\par}}
\def\xxftbib{\futurelet\next\xftbib}
\def\xftbib{\if\next[\let\next\ftbibyear\else\let\next\ftbibnoyear\fi\next}
\def\ftbibyear[#1]{\edef\yearkey{#1}\expandafter\ftbibyearbis\ignorespaces}
\def\ftbibyearbis#1\par{\append\bibt{#1\par}\ftbibend}
\def\ftbibnoyear#1\par{\append\bibt{#1\par}\edef\yearkey{\extractyear#1\par}\ftbibend}
\def\ftbibend{\expandafter\ifx\csname year.\yearkey\endcsname\relax
		\edef\yearindex{0}%
	\else
		\edef\yearindex{\csname year.\yearkey\endcsname}%
	\fi
	{\count0=\yearindex
	\advance\count0 by 1 %
	\expandafter\xdef\csname year.\yearkey\endcsname{\the\count0 }%
	\xdef\alphakey{\ifcase \count0 ?\or a\or b\or c\or d\or e\or f\or g\or h\or i\or j\or k\or l\or m\or n\or o\or p\or q\or r\or s\or t\or u\or v\or w\or x\or y\or z\else .\the\count0 \fi}}%
	\printlabel{key \key, year key \yearkey, alpha key \alphakey}%
	\expandafter\edef\csname text.\key\endcsname{\noexpand\typesetyearalpha{\key}{\yearkey}{\alphakey}}%
	\printlabel{reference \key\space adjustment: id.\key\space = reference.\key, text.\key\space = \yearkey.\alphakey}%
	\setbox0=\hbox{[\yearkey.\alphakey]}\ifdim\bibindent<\wd0 \bibindent=\wd0 \fi
}
\def\typesetyearalpha#1#2#3{%
	\edef\yearindex{\csname year.#2\endcsname}%
	\ifnum\yearindex=1 #2\else#2.#3\fi}
\def\extractyear#1\y#2#3\par{#2}
\newif\iftype
\def\typesetbib#1#2\par{\edef\key{\expandafter\gobble\string#1}%
	\edef\bibbr{\csname backreference-list.reference.\key\endcsname}%
	\typetrue\ifsuppressunusedbib\ifx\bibbr\empty\typefalse
	\fi\fi
	\iftype
	\noindent\hbox to \bibindent{[\anchor{reference.\key}{\csname text.\key\endcsname}]\hfil}#2
	\ifx\bibbr\empty\else\expandafter\stripcomma\bibbr.\fi
	\hangindent\bibindent\filbreak
	\fi}

\let\printverify\gobble
\def\addverunused#1#2{\appendexpand\vlist{\noexpand#2\noexpand#1{\the\inputlineno}}}
\def\verifyref#1#2#3{\printverify{verifying for #3 \string#1 (line #2)}%
	\edef\key{\expandafter\gobble\string#1}%
	\edef\cseq{\csname id.\key\endcsname}%
	\edef\tmp{\csname backreference-list.\cseq\endcsname}%
        \ifx\tmp\empty\ewarningline{#2}{#3 \string#1.}{\string#1}\fi}
\def\verifylabel#1#2{\verifyref#1{#2}{Unused label}}
\def\verifybib#1#2{\verifyref#1{#2}{Unused reference}}

\let\printurl\gobble
\newtoks\urltext
\newtoks\urlt
\newif\ifpunct
\def\urldash{-}
\def\urltilde{{\tensy^^X}} 
\def\ndash{\def\urldash{--}}
\def\http://{\hfil\penalty900\hfilneg\urltext={http://}\urlt={http:/\negthinspace/}\punctfalse\urlgrab}
\def\https://{\hfil\penalty900\hfilneg\urltext={https://}\urlt={https:/\negthinspace/}\punctfalse\urlgrab}
\def\urlgrab{\catcode`\#=11 \catcode`\&=11 \futurelet\ntok\urldispatch}
\def\urldispatch{%
	\ifx\ntok~\let\proceed\urlcont\else
	\ifcat\noexpand\ntok\space\let\proceed\urlfinish\else
	\ifcat\noexpand\ntok\relax\let\proceed\urlfinish\else
	\let\proceed\urlcont
	\fi\fi\fi\proceed}
\def\urlcont#1{\ifpunct\appendexpand\urltext\punctc\appendexpand\urlt\punctc\punctfalse\fi
	\ifx\ntok~\appendexpand\urltext{\noexpand~}\appendexpand\urlt\urltilde
	\else\if\ntok\ampersand\appendexpand\urltext{&}\appendexpand\urlt{\&}%
	\else\if\ntok\ohash\appendexpand\urltext\ohash\appendexpand\urlt\#%
	\else\if\ntok_\appendexpand\urltext_\appendexpand\urlt\_%
	\else\if\ntok-\appendexpand\urltext-\appendexpand\urlt\urldash
	\else\if\ntok.\puncttrue\def\punctc{.}%
	\else\if\ntok,\puncttrue\def\punctc{,}%
	\else\if\ntok;\puncttrue\def\punctc{;}%
	\else\appendexpand\urltext{#1}\appendexpand\urlt{#1}%
	\fi\fi\fi\fi\fi\fi\fi\fi\urlgrab}
\def\urlfinish{\catcode`\#=6 \catcode`\&=4 \hbox{\printurl{\the\urltext}\link{\the\urltext}{\the\urlt}}\ifpunct\punctc\punctfalse\fi\def\urldash{-}}
\def\idgrab{\futurelet\ntok\iddispatch}
\def\iddispatch{\ifcat\noexpand\ntok\space\let\proceed\urlfinish
		\else\if\ntok,\let\proceed\urlfinish
		\else\let\proceed\idcont
		\fi\fi\proceed}
\def\idcont#1{\ifpunct\appendexpand\urltext.\appendexpand\urlt.\punctfalse\fi
	\if\ntok.\puncttrue\def\punctc{.}%
	\else\if\ntok_\appendexpand\urltext_\appendexpand\urlt\_%
	\else\appendexpand\urltext{#1}\appendexpand\urlt{#1}%
	\fi\fi\idgrab}

\let\printgrab\gobble
\newtoks\grabname
\newtoks\grabtoks 
\newtoks\grabcseq 
\newtoks\subsuptoks 
\newtoks\dtoks 
\newcount\grabsize
\newif\ifgrabsubscript 
\newif\ifgrabsupscript 
\def\grabsequence{\bgroup 
	\grabsubscripttrue\grabsupscripttrue 
	\grabstring}
\def\grabalpha{\bgroup 
	\grabsubscriptfalse\grabsupscriptfalse 
	\grabstring}
\def\grabingroup{\ifinfont\errmessage{Already inside a math token}\fi\append\grabtoks{\bgroup\grablink}\infonttrue}
\def\graboutgroup{\ifinfont\append\grabtoks{\endgrablink\egroup}\infontfalse\fi}
\def\grabstring#1#2#3{
	\let\specialhat^
	\catcode`\^=7
	\aftergroup#3 
	\ifx\specialaddon\undefined\else\expandafter\aftergroup\specialaddon\let\specialaddon\undefined\fi
	\inewif\ifdefine \inewif\ifinfont
	\printgrab{}\printgrab{grab a string of type #1, typeset using font \string#2, with postcommand \string#3}%
	\grabname{#1}\def\grabfont{#2}\grabsize0 \grabtoks={}\grabingroup \append\grabtoks{#2}\grabcseq={}%
	\futurelet\ntok\grabdeflookahead}
\def\grabdeflookahead{\if=\noexpand\ntok 
	\definetrue\printgrab{defining}\expandafter\grabgobblefuturelet
	\else\printgrab{referencing}\definefalse\expandafter\grablookahead\fi}
\def\grabgobblefuturelet#1{\futurelet\ntok\grabtestforsilent} 
\newif\ifsilentgrab
\def\grabtestforsilent{\if=\noexpand\ntok \silentgrabtrue \let\ncom\grabsilenteq \else \silentgrabfalse \let\ncom\grablookahead \fi \ncom}
\def\grabsilenteq={\grabfuturelet}
\def\grabfuturelet{\futurelet\ntok\grablookahead}
\def\grablookahead{\printgrab{futurelet token meaning: \meaning\ntok}%
	\let\ncom\grabfinish
	\if\bgroup\noexpand\ntok \printgrab{left brace, terminating}%
	\else \if\egroup\noexpand\ntok \printgrab{right brace, terminating}%
	\else \if\space\noexpand\ntok \printgrab{blank space, terminating}%
	\else \let\ncom\grabexamine \fi\fi\fi \ncom}
\def\grabexamine#1{\printgrab{grabexamine argument: \string#1, meaning \meaning#1}%
	\def\ncom{\grabfinish#1}%
	\ifcat$\ifcat*\string#1\fi$
		\ifcat _\noexpand#1 \ifgrabsubscript\printgrab{subscript, continuing}%
			\graboutgroup \append\grabtoks{#1}\subsuptoks{#1}\def\ncom{\grabsubsupfuturelet}%
						\else\printgrab{subscript, terminating}\fi
		\else \ifcat ^\noexpand#1 \ifgrabsupscript\printgrab{superscript, continuing}%
			\graboutgroup \append\grabtoks{#1}\subsuptoks{#1}\def\ncom{\grabsubsupfuturelet}%
						\else\printgrab{superscript, terminating}\fi
		\else \ifx \specialhat#1 \ifgrabsupscript\printgrab{specialhat superscript, continuing}%
			\graboutgroup \append\grabtoks{#1}\subsuptoks{#1}\def\ncom{\grabsubsupfuturelet}%
						\else\printgrab{superscript, terminating}\fi
		\else \ifcat\noexpand~\noexpand#1 \printgrab{active character \string#1, examining further}%
			\ifnum1=\uccode`#1 \printgrab{UTF-8 letter, continuing}%
				\advance\grabsize1 \append\grabtoks{#1}\appendexpand\grabcseq{\string#1}\def\ncom{\grabfuturelet}%
			\else
				\ifnum\the\grabsize=0 \printgrab{Nothing grabbed so far, continuing}%
					\advance\grabsize1 \append\grabtoks{#1}\appendexpand\grabcseq{\string#1}\def\ncom{\grabfuturelet}%
				\else\printgrab{Not a UTF-8 letter and not the first character in a string, terminating}%
				\fi
			\fi
		\else \ifcat a\noexpand#1 \printgrab{letter #1, continuing}%
			\advance\grabsize1 \append\grabtoks{#1}\append\grabcseq{#1}\def\ncom{\grabfuturelet}%
		\else\printgrab{nonactive character \string#1}%
			\ifnum\the\grabsize=0 \printgrab{sole argument, adding and terminating}%
				\advance\grabsize1 \append\grabtoks{#1}\append\grabcseq{#1}\def\ncom{\grabfinish}%
			\else\printgrab{terminating}\fi
		\fi\fi\fi\fi\fi
	\else \printgrab{command sequence \string#1, terminating}\fi
	\ncom}
\def\grabsubsupfuturelet{\futurelet\ntok\grabsubsuplookahead}
\newcount\dig
\newif\ifdigit
\def\grabsubsuplookahead{\printgrab{subsup futurelet token meaning: \meaning\ntok}%
	\if\bgroup\noexpand\ntok \printgrab{left brace, continuing}\let\ncom\grabentiresubsup%
	\else \if\egroup\noexpand\ntok \printgrab{right brace, continuing}\let\ncom\grabentiresubsup%
	\else \if\space\noexpand\ntok \errmessage{Blank space after \the\subsuptoks}%
	\else \let\ncom\grabsubsupexamine \fi\fi\fi \ncom}
\def\grabentiresubsup#1{\printgrab{subsup entire group added}\grabingroup\append\grabtoks{#1}\graboutgroup\grabfuturelet}
\def\grabsubsupexamine#1{\printgrab{examining subsup argument \string#1, meaning \meaning#1}%
	\ifcat$\ifcat*\string#1\fi$
		\ifcat\noexpand~\noexpand#1 \printgrab{active character \string#1, continuing}%
			\grabingroup\appendexpand\grabtoks{\grabfont\noexpand#1}\let\ncom\grabfuturelet
		\else\ifnum"8000=\the\mathcode`#1 \printgrab{math active character \string#1, continuing}%
			\let\specialaddon\egroup
			\def\ncom{#1}%
			\grabtypeset
		\else\ifcat a\noexpand#1 \printgrab{letter #1, checking whether single or not}%
			\dtoks{#1}\def\ncom{\futurelet\ntok\grabsubsupsecondletterlookahead}%
		\else\printgrab{something else, inserting a single-character sub/superscript, continuing}
			\appendexpand\grabtoks{\bgroup\grabfont\noexpand#1\egroup}%
			\advance\grabsize2 \appendexpand\grabcseq{\the\subsuptoks\string#1}
			\def\ncom{\grabfuturelet}\fi\fi\fi
	\else \printgrab{command sequence \string#1, continuing}%
		\append\grabtoks{#1}\let\ncom\grabfuturelet\fi
	\ncom}
\def\grabsubsupsecondletterlookahead{\def\ncom{\appendexpand\grabtoks{\the\dtoks}\grabfuturelet}%
	\ifcat a\noexpand\ntok \printgrab{not a single letter, grabbing the entire subsupscript}%
		\advance\grabsize1 \appendexpand\grabcseq{\expandafter\string\the\subsuptoks}
		\grabingroup
		\appendexpand\grabtoks{\grabfont\the\dtoks}%
		\appendexpand\grabcseq{\the\dtoks}%
		\def\ncom{\grabfuturelet}%
	\else \printgrab{single letter, continuing}\fi\ncom}
\def\grabfinish{\printgrab{grabfinish}\graboutgroup\grabtypeset\egroup}
\def\grabtypeset{\printgrab{grabtypeset grabsize=\the\grabsize, grabtoks=\the\grabtoks, grabcseq=\the\grabcseq}%
	\def\grablink##1\endgrablink{##1}%
	\ifnum\the\grabsize=0 \errmessage{No string to grab}\fi
	\ifnum\the\grabsize>1 
		\ifdefine 
			\expandafter\checkduplicates\csname\the\grabname.\the\grabcseq\endcsname{Mathematical identifier \key\space already defined at line \keydefline}%
			\iftypesetting 
				\ifsilentgrab
					\expandafter\gdef\csname silent:\the\grabname.\the\grabcseq\endcsname{}
				\else
					\anchor{\the\grabname.\the\grabcseq}{}
				\fi 
			\else
				\expandafter\xdef\csname id.\the\grabname.\the\grabcseq\endcsname{paragraph.\the\secn.\the\parn}
				\predefbackref{paragraph.\the\secn.\the\parn}%
  			\fi
		\else 
			\iftypesetting 
				\global\advance\backref1
				\expandafter\ifx\csname line:\the\grabname.\the\grabcseq\endcsname\relax 
					\warning{Undefined mathematical identifier \the\grabname.\the\grabcseq}%
					\expandafter\gdef\csname\the\grabname.\the\grabcseq\endcsname{\relax}
				\else 
					\expandafter\ifx\csname silent:\the\grabname.\the\grabcseq\endcsname\empty
						\edef\grablink##1\endgrablink{{##1}}%
					\else 
						\edef\grablink##1\endgrablink{\noexpand\llink{\the\grabname.\the\grabcseq}{##1}
							\noexpand\anchor{backreference.\the\backref}{}}
					\fi
				\fi
			\else 
				\ifsuppressbackref\else
					\global\advance\backref1
					\printbackref{math back reference to \the\grabname.\the\grabcseq: backref.\the\backref\space at line \the\inputlineno}%
					\expandafter\recordbackref\csname\the\grabname.\the\grabcseq\endcsname
				\fi
			\fi
		\fi
	\fi
	\ifsilentgrab\else\the\expandafter\grabtoks\fi}

\inewif\ifsilent \inewif\ifanchor
\newif\ifsuppresscs
\catcode`\^=7   \def\specialhat{\ifmmode\def\next{^}\else\let\next\beginxref\fi\next} \catcode`\^=\active \let^=\specialhat
\def\silentxref#1{\futurelet\next\silentxrefswitch}
\def\silentxrefswitch{\silenttrue\xref}
\def\beginxref{\futurelet\next\beginxrefswitch}
\def\beginxrefswitch{\ifx\next\specialhat\let\next\silentxref \else\silentfalse\let\next\xref\fi \next}
\def\xref{\leavevmode\futurelet\next\xrefswitch}
\def\xrefswitch{\ifx\next!\let\next\verbalxref \else \ifx\next=\let\next\anchorxref \else \anchorfalse \let\next\normalxref \fi \fi \next}
\newtoks\vtoksl
{\count0="C2 \loop\ifnum\count0<"F5 \catcode\count0=11 \advance\count0 by 1 \repeat 
\gdef\plainaccents{\suppresscstrue%
	\def\`##1{##1\empty ̀}%
	\def\'##1{##1\empty ́}%
	\def\^##1{##1\empty ̂}%
	\def\"##1{##1\empty ̈}%
	\def\~##1{##1\empty ̃}%
	\def\=##1{##1\empty ̄}%
	\def\.##1{##1\empty ̇}%
	\def\u##1{##1\empty ̆}%
	\def\v##1{##1\empty ̌}%
	\def\H##1{##1\empty ̋}%
	\def\t##1{##1\empty ͡}%
}}
\def\plainaccents{\let\xcsname=\empty \let\xendcsname=\empty}
\def\verbalxref!{\begingroup\plainaccents\verbalxrefaux}
\def\verbalxrefaux#1{%
	\lowercase{\vtoksl{#1}}%
	\expandtoks\vtoksl
	\iftypesetting
		\expandafter\gdef\expandafter\cseq\expandafter{\csname verbal.\the\vtoksl\endcsname}%
		\printlabel{verbal xref \the\vtoksl}%
		\expandafter\ifx\cseq\relax\errmessage{Undefined reference to \the\vtoksl}\else\cseq\fi
	\else
		\ifsuppressbackref\else
			\global\advance\backref1 %
			\printbackref{label verbal.#1: backref.\the\backref\space at line \the\inputlineno}%
			\blah
			\expandafter\recordbackref\csname verbal.\the\vtoksl\endcsname
		\fi
	\fi
	\endgroup}
\def\initverballabelcommand#1{%
	\printlabel{initializing verbal label #1 (\the\curverb\the\secn.\the\parn)}%
	\ifx\sectionid\undefined
		\expandafter\xdef\csname id.verbal.#1\endcsname{paragraph.\the\secn.\the\parn}%
	\else
		\expandafter\xdef\csname id.verbal.#1\endcsname{section.\the\secn}%
	\fi
	\ifnum\parn=0 %
	\expandafter\xdef\csname text.verbal.#1\endcsname{\the\curverb\the\secn}%
	\else
	\expandafter\xdef\csname text.verbal.#1\endcsname{\the\curverb\the\secn.\the\parn}%
	\fi
	\expandafter\initlabelcommand\csname verbal.#1\endcsname
}
\def\anchorxref={\anchortrue\futurelet\next\anchorxrefswitch}
\def\anchorxrefswitch{\ifx\next:\let\next\nonitalicanchor\else\italictrue\let\next\normalxref\fi \next}
\def\nonitalicanchor:{\italicfalse\normalxref}
\newtoks\firsttoks \newtoks\secondtoks \inewif\ifplural \inewif\ifitalic
\def\parseplural#1[#2|#3]{\let\next\parseplural\ifx\hfuzz#2\hfuzz\ifx\hfuzz#3\hfuzz\let\next\relax\else\pluraltrue\fi\else\pluraltrue\fi
	\append\firsttoks{#1#2}\append\secondtoks{#1#3}\next}
\newtoks\firsttoksl
\newtoks\secondtoksl
\newtoks\nexttoks
\def\normalxref{\begingroup\plainaccents\normalxrefaux}
\def\normalxrefaux#1{\firsttoks{}\secondtoks{}\pluralfalse\parseplural#1[|]%
	\lowercase\expandafter{\expandafter\firsttoksl\expandafter{\the\firsttoks}}%
	\expandtoks\firsttoksl
	\lowercase\expandafter{\expandafter\secondtoksl\expandafter{\the\secondtoks}}%
	\expandtoks\secondtoksl
	\ifanchor
		\iftypesetting\else
			\initverballabelcommand{\the\firsttoksl}%
			\ifplural\initverballabelcommand{\the\secondtoksl}\fi
			\predefbackref{paragraph.\the\secn.\the\parn}%
			\expandafter\checkduplicates\csname\the\firsttoksl\endcsname{Verbal label \the\firsttoksl\space was already defined at line \keydefline}
			\ifplural\expandafter\checkduplicates\csname\the\secondtoksl\endcsname{Verbal label \the\secondtoksl\space was already defined at line \keydefline}\fi
		\fi
		\anchor{verbal.\the\firsttoksl}{}
		\ifplural\anchor{verbal.\the\secondtoksl}{}\fi 
	\else
		\ifsuppressbackref\else
			\global\advance\backref1
			\printbackref{label normal.#1: backref.\the\backref\space at line \the\inputlineno}%
			\iftypesetting
			\else
				\expandafter\recordbackref\csname verbal.\the\secondtoksl\endcsname
			\fi
			\anchor{backreference.\the\backref}{}
		\fi
	\fi
	\nexttoks{}%
	\ifsilent
		\nexttoks{\ignorespaces}%
	\else
		\ifanchor
			\ifitalic\nexttoks\expandafter{\expandafter\bgroup\expandafter\it\the\firsttoks\italcorr}%
			\else\nexttoks\expandafter{\the\firsttoks}%
			\fi
		\else
			\edef\tmp{{verbal.\the\secondtoksl}}%
			\nexttoks\expandafter{\expandafter\llink\tmp}
			\nexttoks\expandafter\expandafter\expandafter{\expandafter\the\expandafter\nexttoks\expandafter{\the\firsttoks}}%
		\fi
	\fi
	\expandafter\endgroup\the\nexttoks}
\def\italcorr{\futurelet\next\italcorrtest}
\def\italcorrtest{\if,\noexpand\next\else\if.\noexpand\next\else\/\fi\fi\egroup}

\def\numdam{\urltext={http://www.numdam.org/item/?id=}\urlt={numdam:}\punctfalse\idgrab}

\def\gen:{\http://libgen.rs/book/index.php?md5=}
\def\jstor:{\https://www.jstor.org/stable/}
\def\eudml:{\https://eudml.org/doc/}
\def\repo#1#2{\urltext={#1}\urlt={#2}\punctfalse\idgrab}
\def\arXiv:{\urltext={https://arxiv.org/abs/}\urlt={arXiv:}\punctfalse\idgrab}

\def\Zbl:{\urltext={https://zbmath.org/?q=an:}\urlt={Zbl:}\punctfalse\idgrab}
\def\doi:{\ndash\urltext={https://doi.org/}\urlt={doi:}\punctfalse\urlgrab}

\def\sqr#1#2{{\thinspace\vbox{\hrule height.#2pt \hbox{\vrule width.#2pt height#1pt \kern#1pt \vrule width.#2pt} \hrule height0pt depth.#2pt}\thinspace}}
\def\square{\mathchoice\sqr64\sqr64\sqr{4.2}3\sqr33}

\def\ltoarr#1{\mathop{\count0=#1 \loop\ifnum\count0>0 \smash-\mkern-7mu \advance\count0 -1 \repeat \mathord\rightarrow}\limits} 
\def\lto#1#2{\mathrel{\ltoarr{#1}^{#2}}} 
\def\longto#1^#2_#3{\mathrel{\ltoarr{#1}^{#2}_{#3}}} 
\def\lgetsarr#1{\mathop{\mathord\leftarrow \count0=#1 \loop\ifnum\count0>0 \mkern-7mu\smash-\advance\count0 -1 \repeat}\limits} 
\def\longgets#1^#2_#3{\mathrel{\lgetsarr{#1}\limits^{#2}_{#3}}} 

\def\toto{\mathrel{\vcenter{\hbox{$\to$}\kern-1.5ex \hbox{$\to$}}}}
\def\prearrfill{\smash-\mkern-7mu}
\def\postarrfill{\mkern-7mu\smash-}
\def\midarrfill#1{\cleaders\hbox{$\mkern-2mu\smash-\mkern-2mu$}\hskip0pt plus #1fil}
\def\rightarrfill{\mkern-7mu\mathord\rightarrow}
\def\leftarrfill{\mathord\leftarrow\mkern-7mu\midarrfill1\postarrfill}
\def\ltoto#1#2#3{\ifinner
	\mathrel{\vcenter{\hbox to #1em{$\prearrfill\midarrfill1{\scriptstyle#2}\midarrfill3 \rightarrfill$}%
		\kern-1.5ex \hbox to #1em{$\prearrfill\midarrfill3{\scriptstyle#3}\midarrfill1 \rightarrfill$}}}%
	\else
	\mathrel{\mathop{\vcenter{\hbox to #1em{\rightarrowfill}%
		\kern-1.5ex \hbox to #1em{\rightarrowfill}}}\limits^{#2}_{#3}}%
	\fi}
\def\ltogets#1#2#3{\ifinner
	\mathrel{\vcenter{\hbox to #1em{$\prearrfill\midarrfill1{\scriptstyle#2}\midarrfill3 \rightarrfill$}%
		\kern-1.5ex \hbox to #1em{$\leftarrfill\midarrfill3{\scriptstyle#3}\midarrfill1 \postarrfill$}}}%
	\else
	\mathrel{\mathop{\vcenter{\hbox to #1em{\rightarrowfill}%
		\kern-1.5ex \hbox to #1em{\leftarrowfill}}}\limits^{#2}_{#3}}%
	\fi}

\def\ltogetscore#1#2{\dimen0=\fontdimen6 #1 2 \divide\dimen0 by 2 \multiply\dimen0 by #2 \vcenter{\hbox to \dimen0{\rightarrowfill}\kern-1.8ex \hbox to \dimen0{\leftarrowfill}}}
\def\ltogets#1#2#3{\mathrel{\mathop{\mathchoice{\ltogetscore\textfont{#1}}{\ltogetscore\textfont{#1}}{\ltogetscore\scriptfont{#1}}{\ltogetscore\scriptscriptfont{#1}}}^{#2}_{#3}}}

\def\rx#1#2{\rlap{\kern #1pt \raise#1pt \hbox{#2}}}
\def\dottednearrow{\rx{-8}. \rx{-6}. \rx{-4}. \rx{-2}. \rx0. \rx2. \rx4. \kern6pt \raise7.7pt \hbox{$\nearrow$}}

\def\gmatrix#1#2{\null\,\vcenter{\normalbaselines
	\ialign{#1\crcr
		\mathstrut\crcr\noalign{\kern-\baselineskip}
		#2\crcr\mathstrut\crcr\noalign{\kern-\baselineskip}}}\,}
\def\cdmatrix{\gmatrix{\hfil$##$\hfil&&\enspace\hfil$##$\hfil\enspace&\hfil$##$\hfil}}
\def\sqmatrix{\gmatrix{\hfil$##$&\enspace\hfil$##$\hfil\enspace&$##$\hfil}}
\def\cdbl{\def\normalbaselines{\baselineskip20pt \lineskip3pt \lineskiplimit3pt }}

\def\cd{\cdbl\cdmatrix}
\def\sqcd{\cdbl\let\vagap\;\sqmatrix}

\newcount\arrowsize \arrowsize3
\def\mapright#1{\smash{\lto\arrowsize{#1}}}

\def\rvagap{\vagap} \def\lvagap{\vagap} \def\rvaskip{\vaskip} \def\lvaskip{\vaskip} \def\vaskip{} \def\vagap{}
\def\mapdown#1{\rvagap\Big\downarrow\rlap{$\vcenter{\hbox{$\scriptstyle#1$}}$}\rvaskip}

\def\lmapdown#1{\lvaskip\llap{$\vcenter{\hbox{$\scriptstyle#1$}}$}\Big\downarrow\lvagap}

\newcount\forno \forno0
\def\arrno#1#2{\global\advance\arr1 \edef\eeqnno{\the\arr}%
	\global\advance\forno1 \edef\eforno{\the\forno}%
	\xdef#2{\noexpand\llink{equation.\eforno}{\eeqnno}}
	#1{(\anchor{equation.\eforno}{\eeqnno})}} 
\newbox\mdiag
\def\wrapdiagram{%
	\setbox\mdiag\vtop\bgroup
	\null 
	\vskip\baselineskip
	\inewcount\arr \arr0
	\baselineskip0pt
	\lineskip4pt
	\lineskiplimit4pt
	\let\par\cr
	\obeylines
	\halign\bgroup\hfil$\displaystyle##$\hfil\cr
	\ewrapdiagram}

\def\ewrapdiagram#1{#1
	\egroup
	\egroup
	\vskip0pt plus \dp\mdiag \penalty-250 \vskip0pt plus-\dp\mdiag 
	\hangafter-\dp\mdiag
	\divide\hangafter\baselineskip
	\advance\hangafter-2
	\hangindent-\wd\mdiag
	\advance\hangindent-2em
	\hbox to\hsize{\hfil\dp\mdiag0pt \box\mdiag}%
	\ignorespaces}

\inewtoks\preamble
{\endlinechar=10 \catcode`#=12 \global\preamble{
prologues := 3;

verbatimtex
\let\endprolog

\input dpmac/dpmac
\expandafter\gobbleinit\input }\global\appendonceexpand\preamble\jobname\global\append\preamble{
\catcode`\^=7
etex

input cmarrows
setup_cmarrows(arrow_name = "texarrow"; parameter_file = "cmr10.mf"; macro_name = "drawarrow");
setup_cmarrows(arrow_name = "doublearrow"; parameter_file = "cmr10.mf"; macro_name = "drawdarrow");
def drawmarrow expr p = _apth:=p; _finmarr enddef;
rule_thickness#:=.4pt#;    
def _finmarr text t_ =
  drawarrow subpath(0, 0.5 * length(_apth)) of _apth t_;
  draw subpath(0.5 * length(_apth), length(_apth)) of _apth withpen pencircle scaled rule_thickness# t_;
enddef;

def object(suffix O)(expr x,y)(expr l) =
  save O;       
  pair O;
  O := (x,y) * u;
  picture O.tx;    
  O.tx := thelabel(l,O);
  draw O.tx;                         
enddef;                                    

def smorphism(suffix A,B) =
  save ss, tt;
  ss := xpart ((A..B) intersectiontimes bbox A.tx);
  tt := xpart ((A..B) intersectiontimes bbox B.tx);
  drawarrow subpath(ss,tt) of (A..B);
enddef;

def morphism(suffix A,B)(expr l)(expr f)(suffix $) =
  smorphism(A, B);
  label.$(l, point f[ss, tt] of (A..B));
enddef;
}}%

\newcount\vertex                                      
\ifx\mathspecials\undefined\def\mathspecials{}\fi

\def\grabdiagramaux#1;{\vertex0 \dcom#1,*}                          
\def\dcom{\futurelet\next\dcomswitch}                                          
\def\dcomswitch{\ifcat a\noexpand\next \let\next\dcomalpha                                       
        \else\if*\noexpand\next \addcode{\the\buffertoks}\endinlinemp\let\next\gobble                                                 
        \else\if[\noexpand\next \let\next\grabscale                                                                
        \else\errmessage{Unrecognized diagram command \next}\let\next\relax\fi\fi\fi\next}
\def\grabscale[#1]{\dimen0=#1 \addcode{save u; u = \the\dimen0;}\dcom}
\def\dcomalpha#1{\def\objectname{#1}\futurelet\next\dcomalphaswitch}
\def\dcomalphaswitch{\ifcat a\noexpand\next \let\next\dcommorphism\else                                                          
        \expandafter\edef\csname vertex:\objectname\endcsname{\the\vertex}%
        \if:\noexpand\next\let\next\grabcoords                                              
        \else\if=\noexpand\next\let\next\grabobjectequ                                     
        \else\errmessage{Expected : or = while processing a diagram object, got \meaning\next}\let\next\relax\fi\fi\fi\next}
\def\grabcoords:#1,#2={\toks\vertex{#1,#2}\grabobjectlabel}
\def\grabobjectequ={\edef\tmp{\the\toks\vertex}\ifx\tmp\empty\errmessage{No coordinates specified for vertex \the\vertex: \objectname}\fi\grabobjectlabel}
\def\grabobjectlabel#1,{\addcode{object(\objectname, \the\toks\vertex, }\addunexpandedcode{btex #1 etex);}\advance\vertex1 \dcom}
\def\dcommorphism#1{\def\tobjectname{#1}%
        \def\labelpos{.5}%
        \futurelet\next\dcommorphismswitch}                     
\def\dcommorphismswitch{\if.\noexpand\next \expandafter\grabmorphismdir \else \setupmorphismlabel \expandafter\grabmorphismpos\fi}
\def\setupmorphismlabel{
        \edef\vlabeldira{direction.\objectname.\tobjectname}%
        \edef\vlabeldirb{direction.\tobjectname.\objectname}%
        \expandafter\ifx\csname\vlabeldira\endcsname\relax
        	\expandafter\ifx\csname\vlabeldirb\endcsname\relax
		        \edef\tmpa{\csname vertex:\objectname\endcsname}\expandafter\ifx\tmpa\relax\errmessage{No such vertex: \objectname}\fi
		        \edef\tmpb{\csname vertex:\tobjectname\endcsname}\expandafter\ifx\tmpb\relax\errmessage{No such vertex: \tobjectname}\fi
		        \ifnum \tmpa<\tmpb \edef\vlabeldir{direction.\tmpa.\tmpb}\else \edef\vlabeldir{direction.\tmpb.\tmpa}\fi
		        \expandafter\ifx\csname\vlabeldir\endcsname\relax\errmessage{No label direction specified for morphism \objectname->\tobjectname}\fi
		        \edef\labeldir{\csname\vlabeldir\endcsname}%
		\else
        		\edef\labeldir{\csname\vlabeldirb\endcsname}%
		\fi
	\else
        	\edef\labeldir{\csname\vlabeldira\endcsname}%
	\fi
}                                                                                        
\def\grabmorphismdir.{\def\labeldir{}\futurelet\next\grabmorphismdirec}
{\catcode`\@=13                                                               
\gdef\grabmorphismdirec{\ifx@\next \let\next\grabmorphismposition
        \else\if=\noexpand\next \let\next\grabmorphismequ               
        \else \let\next\grabmorphismdirect\fi\fi\next}                       
\gdef\grabmorphismdirect#1{\edef\labeldir{\labeldir#1}\futurelet\next\grabmorphismdirec}
\gdef\grabmorphismpos{\ifx@\next \expandafter\grabmorphismposition \else \expandafter\grabmorphismequ\fi} 
\gdef\grabmorphismposition@#1={\def\labelpos{#1}\grabmorphismlabel}
}
\def\grabmorphismequ={\grabmorphismlabel}
\def\grabmorphismlabel#1,{%
        \addcodebuffer{morphism(\objectname, \tobjectname, }%
        \addunexpandedcodebuffer{btex \everymath{\scriptstyle}#1 etex, }%
        \addcodebuffer{\labelpos, \labeldir);}%
        \dcom}

\let\printextract\gobble
\def\hinitlabelcommand#1{\printextract{initializing label command \string#1}%
	\gdef#1{\printlabel{invoked label \string#1}
		\ifsuppressbackref
			\edef\key{\expandafter\gobble\string#1}%
			\llink{\csname id.\key\endcsname}{\csname text.\key\endcsname}
		\else
			\global\advance\backref1 %
			\printbackref{label \expandafter\gobble\string#1: backref.\the\backref\space at line \the\inputlineno}%
			\iftypesetting
			\else
				\blah
				\recordbackref#1
			\fi
			\anchor{backreference.\the\backref}{}
			\edef\key{\expandafter\gobble\string#1}%
			\llink{\csname id.\key\endcsname}{\csname text.\key\endcsname}
		\fi
		}}
\let\initlabelcommand\hinitlabelcommand
\def\pinitlabelcommand#1{\printextract{initializing label command \string#1}%
	\gdef#1{\printlabel{invoked plain label \string#1}
		\iftypesetting
			\edef\key{\expandafter\gobble\string#1}%
			\csname text.\key\endcsname
		\else
			\blah
		\fi}}
\newread\labelin
\newif\iflabelcont
\let\terminate=\relax 
\long\def\labelauxaux#1\terminate
{}
\def\preprocesslabel#1{%
	\printextract{Processing label \string#1}%
	\edef\key{\expandafter\gobble\string#1}%
	\initlabelcommand#1
	\expandafter\initlabelcommand\csname v\key\endcsname
	\labelauxaux
}
\def\preprocessbib#1{%
	\printextract{Processing bib \string#1}%
	\initlabelcommand#1
	\labelauxaux
}
\long\def\labelaux#1{\ifx#1\label\let\next\preprocesslabel\else\ifx#1\bib\let\next\preprocessbib\else\let\next\labelauxaux\fi\fi\next}
\def\processoneline{\expandafter\labelaux\labelline\relax\relax\relax\terminate}
\def\preprocess#1{
	\openin\labelin=#1 \labelconttrue
	\loop
		\read\labelin to\labelline
		\ifeof\labelin\let\next\labelcontfalse\else\let\next\processoneline\fi
		\next
	\iflabelcont\repeat
	\expandafter\gobbleinit\input#1
	\ifhmode\par\fi\vfill\eject 
	\ifhmode\par\fi\vfill\supereject
}

\def\importlabels#1{ 
	\recorddupsfalse
	\let\initlabelcommand\pinitlabelcommand
	\preprocess{#1}%
	\let\initlabelcommand\hinitlabelcommand
	\recorddupstrue
	\cont={}
	\vlist={}
	\let\chapname\undefined \secn0 \parn0 \backref0 
}

\newif\iftypesetting 

\typesettingfalse 
\def\blah{blah} 
\output{\setbox0\box255 \setbox0\box\footins \deadcycles0 } 
\let\bib\tbib 
\let\refs\relax 
\everypar{\numpar}
\def\par{\finishpar} 
\newif\ifmetapost \metapostfalse
\edef\filestem{\jobname.gen}
\newcount\figno
\newwrite\mpout
\newtoks\outtoks
\def\inimp{\ifmetapost\else\global\metaposttrue\immediate\openout\mpout=\filestem.mp \addcode{\the\preamble}\fi}
\newdimen\oldhsize
\oldhsize=\hsize
\newdimen\oldvsize
\oldvsize=\vsize
\hsize\maxdimen
\vsize\maxdimen
\preprocess\jobname 
\hsize\oldhsize
\vsize\oldvsize
\ifx\fmtname\plainfmtname
\hsize\plaintextwidth
\vsize\plaintextheight
\fi
\printerr{(\jobname.tex}
\printbackref{list of all back references: \the\backreflist}%
\the\backreflist 
\ifsuppressunusedbib\def\verifybib#1#2{}\fi
\the\vlist 
\printerr{)}
\ifmetapost
\addcode{end}
\immediate\closeout\mpout
\immediate\write18{mpost -interaction nonstopmode \filestem.mp}
\ifeof18 \warning{Compile the METAPOST file \filestem.mp manually using mpost \filestem.mp}\let\runmp\warning\fi
\fi
\let\addcode\gobble 

\typesettingtrue
\def\importlabels#1{} 
\ifscroll \vsize\maxdimen\inewbox\abox\output{\setbox\abox\vbox{\unvbox255\unskip}\shipbox\abox}%
\else 
	\ifx\fmtname\plainfmtname
	\totalht\plaintextheight
	\advance\totalht.1in
	\advance\totalht2\vmargin
	\totalwd\plaintextwidth
	\advance\totalwd2\hmargin
	\setpapersize\totalwd\totalht
	\advance\hoffset\hmargin
	\advance\voffset\vmargin
	\fi
	\output{\plainoutput}
\fi
\let\blah\undefined 
\def\bib#1\par{} 
\def\refs{\raggedright\rightskip0em plus \maxdimen \advance\bibindent1em \everypar{}\the\bibt\ignorespaces} 
\def\addverunused#1#2{} 
\def\addcont#1#2#3{} 
\def\prepass#1{} 
\parbreftrue 
\everypar{\numpar}
\let\chapname\undefined \secn0 \backref0 
\expandafter\gobbleinit\input\jobname\relax
\ifhmode\par\fi\vfill\eject 
\ifhmode\par\fi\vfill\supereject
\end